\documentclass[10pt]{amsart}
\usepackage{amsmath,amsfonts,amssymb,color,a4,epsfig,graphics}
\usepackage[english]{babel}
\usepackage{enumerate, mathrsfs}
\usepackage{enumitem}
\usepackage[colorlinks=true,hyperindex=true]{hyperref}
\allowdisplaybreaks
\usepackage{verbatim}
\usepackage[table]{xcolor}
\newcolumntype{R}[1]{>{\raggedleft\arraybackslash }b{#1}}
\newcolumntype{L}[1]{>{\raggedright\arraybackslash }b{#1}}
\newcolumntype{C}[1]{>{\centering\arraybackslash }b{#1}}
\usepackage{caption}
\usepackage{boldline}
\usepackage{booktabs}
\usepackage{array}
\usepackage{colortbl}
\usepackage{geometry}

\rowcolors{1}{}{white}
\usepackage{tikz}
\usetikzlibrary{arrows.meta,positioning,fit,calc,backgrounds}

\newcommand{\Ec}{\mathcal{E}}
\newcommand{\Fc}{\mathcal{F}}
\newcommand{\Gc}{\mathcal{G}}

\newcommand{\Nc}{\mathcal{N}}

\newcommand{\Vc}{\mathcal{V}}

\newcommand{\supp}{{\rm supp }}

\numberwithin{equation}{section}

\newtheorem{theorem}{Theorem}[section]
\newtheorem{proposition}[theorem]{Proposition}
\newtheorem{lemma}[theorem]{Lemma}
\newtheorem{remark}[theorem]{Remark}

\newtheorem{definition}[theorem]{Definition}

\newtheorem{corollary}[theorem]{Corollary}
\newtheorem{conjecture}[theorem]{Conjecture}

\begin{document}
\title[Geometric Criteria for ESA of Discrete Hodge Laplacians]{Geometric Criteria for Essential Self-Adjointness of Discrete Hodge Laplacians on Weighted Simplicial Complexes}
\author[\tiny Marwa Ennaceur]{Marwa Ennaceur$^{1}$}
 \address{$^{1}$Department of Mathematics, Faculty of Sciences of Sfax, University
of Sfax, 3000 Sfax, Tunisia}
\email{\tt ennaceur.marwa27@gmail.com}
\author[\tiny Amel Jadlaoui]{Amel Jadlaoui$^{2}$ }
 \address{$^{2}$Department of Mathematics, Faculty of Sciences of Sfax, University
of Sfax, 3000 Sfax, Tunisia}
\email{\tt amel.jadlaoui@yahoo.com}

\begin{abstract}
We develop a geometric framework for the essential self-adjointness (ESA) of discrete Hodge Laplacians on weighted simplicial complexes of arbitrary dimension. Three notions of $\chi$-completeness are introduced --- \emph{global}, \emph{local by level}, and \emph{local by region} --- which are formally distinct as definitions, with the pairwise logical separations remaining partly conjectural.

The main results are: \emph{(i)} ESA of the Gauss--Bonnet operator $D=d+\delta$ and the Hodge Laplacian $L=D^2$ on $\bigoplus C_c^i(\Vc)$ under operative geometric hypotheses (existence of cut-offs with finite support and bounded weighted degree), which are implied in particular by each of the three completeness notions; \emph{(ii)} ESA of the individual Laplacian block $L_\ell$ under the corresponding local-by-level hypotheses, via a quadratic-form argument; \emph{(iii)} ESA via a Kato--Rellich coupling decomposition under local-by-region hypotheses, applied to discrete half-spaces in the lattice $\mathbb{Z}^d$ ($d\geq 2$) where the coupling has operator norm $\sqrt{2}$ (graph case, i.e.\ $n=2$ in our convention; see Remark~\ref{rem:n-convention-intro}), with a sufficient condition extending the argument to $n\geq 3$; \emph{(iv)} a divergence criterion $\sum_k w_k=\infty$ guaranteeing ESA of $L$ even in the absence of $\chi$-completeness, with a worked example (a polynomial-branching tree) where ESA holds yet $\chi$-completeness fails.

The relationship between the three notions is partly clarified: global $\chi$-completeness implies local $\chi$-completeness at every level (Remark~\ref{rem:hierarchy-direct}); the converse non-implications are formulated as conjectures, motivated by uniform energy lower bounds in concrete examples (Section~\ref{sec:examples}). The framework recovers and extends classical results for graphs and triangulations, including the optimal divergence rate of \cite{BGJ}.
\end{abstract}
\maketitle
\tableofcontents

\section{Introduction}
\label{sec:introduction}

In recent years, many researchers have studied spectral graph theory \cite{AEGJ1,AEGJ2,KeLe,CdV1,Lim,BG,EJ,BaKe}, mainly regarding essential self-adjointness of the discrete Laplacian; see \cite{Mi,OMI,Gol2,CTT,CTT2}.

Spectral graph theory has been extended beyond one-dimensional structures to higher-order discrete objects such as simplicial complexes. While graphs correspond to 1-dimensional simplicial complexes, triangulations and higher-dimensional networks naturally lead to $n$-simplicial complexes, where topological and geometric features interact in nontrivial ways; see \cite{AAcT,Che,AnTo,ABDE,EJ,BaKe}.

A central question in this setting is the \emph{essential self-adjointness} (ESA) of discrete Laplacians, i.e., whether the operator admits a unique self-adjoint extension. This property is crucial for a well-defined spectral theory, quantum dynamics, and Hodge decomposition.

In the graph case ($n=2$ in our convention, i.e., 0- and 1-simplices only), \cite{AnTo} introduced the notion of \emph{$\chi$-completeness}, a geometric criterion ensuring ESA of the Gauss--Bonnet operator $D = d + \delta$ and the associated Laplacian $L = D^2$. This idea was later extended to 2-dimensional triangulations \cite{BSJ,Che} and magnetic graphs \cite{ABDE}.

\begin{remark}[Convention on the parameter $n$]
\label{rem:n-convention-intro}
Throughout this paper, $n\geq 2$ denotes the maximum number of vertices per simplex, so that $S_n$ has simplices of dimensions $0,\dots,n-1$. Thus $n=2$ gives graphs, $n=3$ triangulations, and $n=4$ tetrahedral complexes. See Remark~\ref{rem:indexing} for the detailed convention.
\end{remark}

A more general framework for essential self-adjointness of Hodge Laplacians on abstract simplicial complexes has been developed in \cite{BaKe}, where the authors introduce intrinsic-metric criteria. In particular, \cite[Theorem~4.11 (Gaffney-type theorem)]{BaKe} establishes essential self-adjointness under the assumption that all metric balls are finite for a certain intrinsic metric on the simplicial complex; the proof goes through \cite[Proposition~4.12]{BaKe}, which is the abstract criterion for general intrinsic pseudo-metrics. As indicated in the remark following \cite[Theorem~4.11]{BaKe}, the cut-off function approach of \cite{AnTo} is equivalent to the existence of an intrinsic metric with finite distance balls, via \cite[Theorem~11.31]{KLW}. Therefore, our results are \emph{compatible with} the framework of \cite{BaKe}, and Theorems~\ref{thm:global} and~\ref{thm:local-level} can be viewed as concrete realizations of an approach sharing the same spirit as \cite[Proposition~4.12]{BaKe}; however, the precise equivalence between $\chi$-completeness and intrinsic-metric completeness in the setting of weighted simplicial complexes remains conjectural in the present generality and is not established in this paper.

The basic framework on which we build --- weighted oriented simplicial complexes, cochain spaces, coboundary and codifferential operators, and their closability --- was developed by the present authors in \cite{EJ}, to which we refer the reader for proofs of the foundational results we use freely in Section~\ref{sec:operators}.

\begin{remark}[Relation with the intrinsic-metric approach of {\cite{BaKe}}]
\label{rem:BaKe-comparison}
The results of~\cite{BaKe} are established for general abstract simplicial complexes via
intrinsic metrics, whereas our results are tailored to \emph{clique complexes} arising from
weighted graphs. In our setting, cut-off functions are defined on vertices and extended to
simplices by averaging (Section~\ref{sec:preliminaries}), giving a computable geometric criterion. The two approaches are compatible --- both yield essential
self-adjointness --- but neither subsumes the other in general. More precisely: we expect
that an intrinsic metric with finite balls and bounded jump size does not in general yield a
plateau cut-off sequence satisfying our energy bounds, and conversely that
$\chi$-completeness does not in general produce an intrinsic metric in the sense
of~\cite{BaKe}; however, we do not have explicit counterexamples for either direction, and
the precise relationship between the two notions on weighted clique complexes remains an
open problem.
\end{remark}

Our work goes beyond this general framework in several directions. First, we introduce
three notions of $\chi$-completeness ---global, local by level, and local by region---
which are \emph{formally distinct as definitions} but whose pairwise logical separation is
partly conjectural; we analyze their precise spectral implications. We prove that global
$\chi$-completeness implies local $\chi$-completeness at every level
(Remark~\ref{rem:hierarchy-direct}); the converse non-implications are stated as
Conjectures~\ref{conj:region-not-global}--\ref{conj:all-levels-not-global}, and partial
supporting evidence (uniform energy lower bounds for layer-constant cut-offs) is provided
in Section~\ref{sec:examples}. Second, we prove essential self-adjointness under
\emph{localized defects} via a compact coupling argument (Theorem~\ref{thm:local-region}),
a technique absent from the intrinsic-metric approach. Third, we establish a
\emph{divergence criterion} (Theorem~\ref{thm:divergence}) that guarantees essential
self-adjointness even in the \emph{absence} of $\chi$-completeness, thereby providing,
for higher dimensions, a partial analogue of the open question raised in \cite{BGJ} about
the necessity of geometric completeness. Finally, our setting is tailored to
\emph{clique complexes} arising from weighted graphs, where cut-off functions are defined
on vertices and extended by averaging---a natural framework for applications in
topological data analysis---which allows for explicit energy estimates and concrete
examples (Section~\ref{sec:examples}) that illustrate the scope of our results.

\subsection*{Main contributions}

Our main contribution is the introduction and analysis of three formally distinct
notions of $\chi$-completeness, whose pairwise strict separation is conjectural (see
Conjectures~\ref{conj:region-not-global}--\ref{conj:all-levels-not-global}):
\begin{enumerate}
\item \textbf{Global $\chi$-completeness} (Definition~\ref{def:global-chi}): guarantees ESA of the full Gauss--Bonnet operator $D$ and the total Laplacian $L = D^2$ (Theorem~\ref{thm:global}, Corollary~\ref{cor:esa-laplacian});
\item \textbf{Local $\chi$-completeness by level} (Definition~\ref{def:local-chi-level}): ensures ESA of the individual Laplacian block $L_\ell$ (Theorem~\ref{thm:local-level}), enabling modular spectral analysis;
\item \textbf{Local $\chi$-completeness by region} (Definition~\ref{def:local-chi-region}): combined with a compact coupling argument, recovers global ESA even in the presence of localized defects (Theorem~\ref{thm:local-region}).
\end{enumerate}
We prove that the use of \emph{plateau cut-off functions}---constant on exhaustion sets and decaying only near their boundary---is not a technical convenience but a \emph{necessity}: it localizes the approximation error to a thin layer, making the energy bound both meaningful and verifiable.

Moreover, we show that $\chi$-completeness, while sufficient, is \emph{not necessary} for essential self-adjointness. For complexes admitting a 1-dimensional decomposition, a divergence criterion on the growth rate (Theorem~\ref{thm:divergence}) can guarantee ESA even in the absence of any geometric cut-off structure. The criterion requires only $\sum_k w_k = \infty$, matching the optimal condition known in the graph case~\cite{BGJ}.

The examples and partial results of Section~\ref{sec:examples} illustrate the scope of our results; the strict separation between the three notions of $\chi$-completeness and the non-ESA of the Gauss--Bonnet operator on the \emph{exponentially weighted} binary tree are stated as conjectures, while the question of whether $L_2$ is essentially self-adjoint in the alternating triangulation remains an open problem. (For the binary tree with \emph{unit} weights, the Gauss--Bonnet operator is essentially self-adjoint, since the weighted degree is bounded; see Remark~\ref{rem:non-esa-tree-status}.)

\subsection*{Structure of the paper}

Sections~\ref{sec:preliminaries} and~\ref{sec:operators} introduce the basic definitions of weighted simplicial complexes, function spaces, and the Gauss--Bonnet operator, following \cite{EJ}. Section~\ref{sec:geometric-hypotheses} defines the three notions of $\chi$-completeness and states and proves the four main theorems. Section~\ref{sec:examples} provides concrete examples and discusses partial results bearing on the strict hierarchy of the three notions. Section~\ref{sec:comparison} gives a comparison with \cite{BGJ} and treats the combinatorial case via weighted regularization. Section~\ref{sec:conclusion} summarizes the main open problems and outlines possible extensions.

\section{Preliminaries}
\label{sec:preliminaries}

We introduce the basic objects used throughout this paper: weighted graphs, oriented weighted simplicial complexes (see \cite{FRK,CdV1,GW,HJ,Fri,EJ}), and the associated Hilbert spaces of cochains. The definitions and elementary results of this section are taken from~\cite{EJ}; we recall them briefly to fix notation, and refer the reader to that paper for full proofs.

\begin{remark}[Indexing convention]
\label{rem:indexing}
Throughout this paper we adopt \emph{zero-based indexing} for simplex vertices uniformly:
an $i$-simplex has $i+1$ vertices labeled $(x_0,\dots,x_i)$ with $i\geq 0$. We use the
parameter $n\geq 2$ to denote the \emph{maximum number of vertices per simplex}, so that
$S_n$ contains simplices with $1, 2, \ldots, n$ vertices, i.e.\ of dimensions
$0, 1, \ldots, n-1$. In particular:
\begin{itemize}
\item $S_n$ has simplices of dimension \emph{at most $n-1$};
\item the cochain spaces are indexed $C_c^i(\Vc)$ for $0\leq i\leq n-1$;
\item when we write ``no simplex of dimension $\geq k$'', we mean no simplex with $\geq k+1$ vertices;
\item the phrase ``$n \leq 2$'' for the clique complex of $\mathbb{Z}^d$ (Remark~\ref{rem:Zd-adjacency}) means simplices with at most $2$ vertices, equivalently of dimension at most $1$.
\end{itemize}
This convention is used uniformly in all cochain spaces, coboundary operators, and completeness conditions. Conditions involving simplicial dimension $\geq 2$ are vacuous in $1$-dimensional examples and should be disregarded accordingly.
\end{remark}

\subsection*{Weighted Graphs}

Let $m_0:\Vc\longrightarrow(0,+\infty)$ be a positive weight on vertices and let $m_1:\Vc\times\Vc\longrightarrow[0,+\infty)$ be a symmetric function, where $\Vc$ is a countable set. An unoriented graph $\Gc:=(\Vc,m_0,m_1)$ consists of a vertex set $\Vc$ and edge set
\[
\Ec:=\{(x,y)\in\Vc\times\Vc : m_1(x,y)>0\}.
\]
Let $x,y\in \Vc$; we say $x$ and $y$ are \emph{neighbors} and write $x\sim y$ if $(x,y)\in \Ec$. If $x\sim x$, we say there is a \emph{loop} at $x$. We denote by
\[
\Nc_{\Gc}(x):=\{y\in\Vc : x\sim y\}
\]
the set of neighbors of $x$. The graph $\Gc$ is assumed to be \emph{simple} (no loops: $m_1(x,x)=0$ for all $x$) and \emph{locally finite} (each vertex has finitely many neighbors) unless explicitly stated otherwise. We do not assume $m_0\equiv 1$ or $m_1\in\{0,1\}$: weights $m_0(x)>0$ and $m_1(x,y)\geq 0$ may be arbitrary positive resp.\ non-negative reals. The \emph{graph distance} is
\[
d(x, y) := \min\bigl\{k \in \mathbb{N} : \exists\, e_1,\dots,e_k \in \Ec \text{ forming a path between } x \text{ and } y\bigr\}.
\]

\begin{remark}[Adjacency in $\mathbb{Z}^d$ and cliques]
\label{rem:Zd-adjacency}
Throughout this paper, the lattice graph $\mathbb{Z}^d$ is equipped with the \emph{$L^1$-adjacency}: two vertices $x,y\in\mathbb{Z}^d$ are neighbors if and only if $\|x-y\|_1=1$, i.e., they differ in exactly one coordinate by $\pm 1$. This is the standard nearest-neighbor graph on the integer lattice.

We record the following key structural fact:
\emph{the clique complex of $\mathbb{Z}^d$ under $L^1$-adjacency contains no simplex with $\geq 3$ vertices (no triangle)},
and the relevant range for $n$-simplicial complexes built from $\mathbb{Z}^d$ with
$L^1$-adjacency is therefore $n \leq 2$ (simplices of dimension at most~$1$, with at
most~$2$ vertices per simplex; see Remark~\ref{rem:indexing}).

\begin{proof}[Proof of the clique bound]
Let $u, v \in \mathbb{Z}^d$ with $u \sim v$, so $u$ and $v$ differ in exactly one coordinate,
say coordinate $j$, by $\pm 1$. Suppose $w \in \mathbb{Z}^d$ is a common neighbor of $u$
and $v$: $\|w - u\|_1 = 1$ and $\|w - v\|_1 = 1$. Since $v = u \pm e_j$, we need
$\|w - u\|_1 = 1$ and $\|w - (u \pm e_j)\|_1 = 1$. Setting $w - u = \pm e_k$ for some
coordinate $k$:
\begin{itemize}
\item If $k = j$: $w = u \pm e_j = v$ or $w = u \mp e_j$. If $w = v$, then $w$ is not
  distinct. If $w = u \mp e_j$, then $\|w - v\|_1 = \|{-2e_j}\|_1 = 2 \neq 1$.
\item If $k \neq j$: $w = u \pm e_k$, so $w - v = \pm e_k \mp e_j$, giving
  $\|w - v\|_1 = 2 \neq 1$.
\end{itemize}
In all cases there is no common neighbor $w$ of $u$ and $v$ in $\mathbb{Z}^d$ with
$L^1$-adjacency. Therefore no triangle (3-clique) exists, and every clique has size at
most $2$.
\end{proof}

In particular, in all propositions below that involve the clique complex of $\mathbb{Z}^d$
with $L^1$-adjacency, the only relevant simplicial dimension is $\ell = 1$ (edges), and
conditions involving $\ell = 2$ or higher are vacuous.
\end{remark}

An \emph{oriented} graph is given by a partition $\Ec=\Ec^-\cup \Ec^+$ such that $(x,y)\in\Ec^-\Leftrightarrow (y,x)\in \Ec^+$. The \emph{weighted degree} of $x\in\Vc$ is
\[
d_{\Vc}(x):=\frac{1}{m_0(x)}\sum_{y\in\Vc}m_{1}(x,y).
\]

\subsection*{Clique Complexes}

\begin{remark}[Clique complex structure]
\label{rem:clique}
Throughout this paper, the \emph{$n$-simplicial complex} $S_n$ is the \emph{clique complex} associated with a weighted graph $\Gc=(\Vc,m_0,m_1)$: an ordered tuple $(x_0,\dots,x_i)$ is an $i$-simplex of $S_n$ if and only if $\{x_0,\dots,x_i\}$ forms a complete subgraph of $\Gc$, i.e., $m_1(x_j,x_k)>0$ for all $0\leq j<k\leq i$. This construction is standard in topological data analysis and is more restrictive than a general abstract simplicial complex. The single exception is the worked example of Section~\ref{sec:example-detailed}, where we work with an abstract simplicial complex; the framework of Sections~\ref{sec:operators}--\ref{sec:geometric-hypotheses} extends to that setting without modification, as discussed there.

In particular, for a clique complex, the set of \emph{common neighbors}
\[
F(x_0,\dots,x_{i-1}) := \bigcap_{j=0}^{i-1} \Nc_\Gc(x_j)
\]
is exactly the set of vertices $x_i$ such that $(x_0,\dots,x_i)$ is an $i$-simplex of $S_n$.
\end{remark}

Let $n\in \mathbb{N}\setminus\{0,1\}$. We denote by $G_{n,c}$ the set of complete subgraphs of $\Gc$ with at most $n$ vertices. For $0\leq i\leq n-1$, set
\[
\Fc_i:=\bigl\{(x_0,\dots,x_i) : \{x_0,\dots,x_i\}\text{ is the vertex set of some } \Gc'\in G_{n,c}\bigr\}
\]
and let $P_i:=\Fc_i/{\cong}$, where $(x_0,\dots,x_i)\cong(y_0,\dots,y_i)$ iff one is an even permutation of the other.

\begin{remark}[Oriented simplices]
\label{rem:orientation}
Each geometric $i$-simplex $\{x_0,\dots,x_i\}$ has exactly $(i+1)!$ ordered representatives, falling into two orientation classes. An element $\tau\in P_i$ is an equivalence class under even permutations; its opposite orientation is $-\tau$. A cochain $f:P_i\to\mathbb{C}$ satisfies $f(-\tau)=-f(\tau)$ by convention.

Choosing an orientation of $P_i$ means fixing a partition $P_i=P_i^+\sqcup P_i^-$ with $(x_0,\dots,x_i)\in P_i^+\Leftrightarrow -(x_0,\dots,x_i)\in P_i^-$. Equivalently, the alternation rule on cochains reads: for any cochain $f:P_i\to\mathbb{C}$ and any permutation $\sigma$ of $\{0,\ldots,i\}$,
\[
f(x_{\sigma(0)},\dots,x_{\sigma(i)})=(-1)^{\varepsilon(\sigma)}\,f(x_0,\dots,x_i).
\]
\end{remark}

The weight functions $(m_i)_{0\leq i\leq n-1}:\Vc^{i+1}\to[0,+\infty)$ satisfy:
\[
m_i(x_0,\dots,x_i)>0 \;\Leftrightarrow\; (x_0,\dots,x_i)\in\Fc_i
\]
and $m_i(x_{\sigma(0)},\dots,x_{\sigma(i)})=m_i(x_0,\dots,x_i)$ for every permutation $\sigma$.

We write $S_n:=(\Vc,(m_i)_{0\leq i\leq n-1})$ for the weighted oriented $n$-simplicial complex associated with $\Gc$. We say $S_n$ is a \emph{triangulation} if $n=3$ (simplices of dimension at most $2$, i.e.\ vertices, edges, and triangles); it is \emph{simple} if $\Gc$ is simple and $m_i:\Vc^{i+1}\to\{0,1\}$ for $1\leq i\leq n-1$.

\begin{remark}[Finiteness of the weighted degree]
\label{rem:finite-degree}
For a locally finite complex $S_n$, the set $F(x_0,\dots,x_{i-1})$ of common neighbors is
finite for every $(i-1)$-simplex $(x_0,\dots,x_{i-1})$: indeed, $F(x_0,\dots,x_{i-1})
\subset \Nc_\Gc(x_0)$, and local finiteness of $\Gc$ ensures $|\Nc_\Gc(x_0)|<\infty$.
Therefore the sum defining $d_{m_{i-1}}(x_0,\dots,x_{i-1})$ below is always finite.
\end{remark}

For $(x_0,\dots,x_{i-1})\in P_{i-1}$ (with $i\geq 1$), we define its \emph{weighted degree} by
\[
d_{m_{i-1}}(x_0,\dots,x_{i-1}) := \frac{1}{m_{i-1}(x_0,\dots,x_{i-1})} \sum_{x_i\in F(x_0,\dots,x_{i-1})} m_i(x_0,\dots,x_i).
\]

\subsection*{Cut-off Functions and Their Extension to Simplices}

A \emph{cut-off function} is a map $\chi:\Vc\to[0,1]$. It is extended to $i$-simplices by the \emph{averaging formula}:
\[
\tilde{\chi}^{(i)}(x_0,\dots,x_i) := \frac{1}{i+1}\sum_{j=0}^i \chi(x_j).
\]
Since $\chi$ is real-valued and the formula is symmetric under all permutations of
$(x_0,\dots,x_i)$, the function $\tilde\chi^{(i)}$ is invariant under \emph{all}
permutations of its arguments, and in particular under orientation reversal:
$\tilde\chi^{(i)}(-(x_0,\dots,x_i))=\tilde\chi^{(i)}(x_0,\dots,x_i)$
for all $i\geq 0$. Consequently, multiplication of an alternating
$i$-cochain $f$ by the symmetric function $\tilde\chi^{(i)}$
yields another alternating $i$-cochain.

For a cut-off function $\chi:\Vc\to[0,1]$ and $1\leq\ell\leq n-1$, we define the
\emph{global energy at level $\ell$} by
\begin{equation}
\label{eq:global-energy}
E_\ell(\chi) := \frac{1}{\ell!}\sum_{(x_0,\dots,x_{\ell-1})\in P_{\ell-1}^+}
\frac{1}{m_{\ell-1}(x_0,\dots,x_{\ell-1})}
\sum_{x_\ell\in F(x_0,\dots,x_{\ell-1})}
m_\ell(x_0,\dots,x_\ell)\,
\sum_{j=0}^\ell |a_j(\chi;x_0,\dots,x_\ell)|^2,
\end{equation}
where $P_{\ell-1}^+$ denotes a choice of canonical orientation of $(\ell-1)$-simplices and
$a_j(\chi;x_0,\dots,x_\ell)$ are the Leibniz coefficients defined in~\eqref{eq:aj-def} of
Lemma~\ref{lem:leibniz}: explicitly,
\[
a_j(\chi;x_0,\dots,x_\ell) := \frac{S - (\ell+1)\,\chi(x_j)}{\ell(\ell+1)}, \qquad S := \sum_{k=0}^{\ell}\chi(x_k).
\]
The formula above coincides with~\eqref{eq:aj-def} in Lemma~\ref{lem:leibniz} after the
index shift $i=\ell-1$, where $i$ in Lemma~\ref{lem:leibniz} indexes the degree of the
\emph{input cochain} $f\in C_c^i(\Vc)$ on which $d_i$ acts (so $d_i f$ lives on
$(i+1)$-simplices), while $\ell$ in~\eqref{eq:global-energy} indexes the dimension of
the \emph{ambient simplex} $(x_0,\ldots,x_\ell)$ at which the coefficients $a_j$ are
evaluated. The two conventions are related by: the $(i+1)$-simplex of Lemma~\ref{lem:leibniz}
is precisely the $\ell$-simplex of~\eqref{eq:global-energy}. Substituting $i+2=\ell+1$ and
$(i+1)(i+2)=\ell(\ell+1)$ in~\eqref{eq:aj-def} yields the displayed form.
\emph{Throughout this paper, we use $i$ when discussing the action of $d_i$ on input
cochains, and $\ell$ when discussing the dimension of an ambient simplex; these
conventions are consistent under the shift $i=\ell-1$.}

The value of $E_\ell(\chi)$ is independent of the choice of $P_{\ell-1}^+$. Indeed, the inner sum
\[
\frac{1}{m_{\ell-1}(x_0,\dots,x_{\ell-1})}\sum_{x_\ell\in F(x_0,\dots,x_{\ell-1})}m_\ell(x_0,\dots,x_\ell)\sum_j|a_j|^2
\]
is invariant under any permutation of $(x_0,\dots,x_{\ell-1})$: $m_{\ell-1}$ and $F(x_0,\dots,x_{\ell-1})$ are
symmetric in their arguments, and for each $x_\ell\in F$, the factor $m_\ell(x_0,\dots,x_\ell)\sum_j|a_j|^2$
depends only on the multiset $\{x_0,\dots,x_\ell\}$ (Remark~\ref{rem:overcounting}).

\begin{remark}[Counting convention]
\label{rem:overcounting}
Each geometric $\ell$-simplex with vertex set $\{x_0,\dots,x_\ell\}$ has exactly
$\ell+1$ faces of dimension $\ell-1$, so the iterated sum in~\eqref{eq:global-energy}
visits each $\ell$-simplex with $\ell+1$ different choices of base $(\ell-1)$-face. While
the individual coefficient $a_j$ depends explicitly on the choice of apex $x_j$, the
\emph{sum} $\sum_{j=0}^\ell|a_j|^2$ depends only on the multiset of $\chi$-values
$\{\chi(x_0),\ldots,\chi(x_\ell)\}$ (as is clear from the formula
$a_j = (S-(\ell+1)\chi(x_j))/(\ell(\ell+1))$, which is invariant under any
permutation of the labels $0,\ldots,\ell$ in the obvious sense). Hence the contribution
of a fixed geometric $\ell$-simplex to formula~\eqref{eq:global-energy} is the same for
each of its $\ell+1$ representations as ``base + apex''. Thus
formula~\eqref{eq:global-energy} counts each geometric $\ell$-simplex with multiplicity
$\ell+1$. We retain this convention (rather than dividing by an extra factor $\ell+1$)
because it matches the natural appearance of $E_\ell$ in the proof of
Lemma~\ref{lem:leibniz}, where the iterated sum corresponds to summing the squared
remainder $|R_{\ell-1}^{(d)}|^2$ over base $(\ell-1)$-simplices of $\ell$-cofaces; the
$\ell+1$ multiplicity is then absorbed into the constants $C$ throughout.

In the graph case $\ell=1$ specifically, this convention counts each undirected edge
\emph{twice} (once with each endpoint as the base vertex). Compared with the
standard graph $\chi$-completeness energy $\sum_{(x,y)\in\Ec}m_1(x,y)|\chi(x)-\chi(y)|^2$
of~\cite{AnTo,BGJ}, our $E_1(\chi)$ is therefore equal to that standard energy multiplied
by an absolute constant ($1$ when ordered pairs $(x,y)$ are summed, or $2$ when undirected
edges are summed); the difference is absorbed into the energy constant $C$ throughout
and does not affect any of our conclusions.
\end{remark}

\begin{remark}[Justification of the energy formula~\eqref{eq:global-energy}]
\label{rem:E-ell-caveat}
For $\ell=1$, formula~\eqref{eq:global-energy} reduces (up to an absolute constant) to the
standard discrete Dirichlet energy: a direct calculation gives $a_0 = (\chi(x_1)-\chi(x_0))/2$
and $a_1 = (\chi(x_0)-\chi(x_1))/2$, so $\sum_j |a_j|^2 = \tfrac{1}{2}|\chi(x_1)-\chi(x_0)|^2$,
recovering the standard form (modulo a factor $1/2$ that can be absorbed into the constant
$C$).

For $\ell\geq 2$, definition~\eqref{eq:global-energy} differs from the naive expression $|d_{\ell-1}\tilde\chi^{(\ell-1)}|^2$. The key reason for using the Leibniz coefficients $a_j$ rather than $d_{\ell-1}\tilde\chi^{(\ell-1)}$ directly is the following: a direct calculation (see Remark~\ref{rem:tri-pointwise-detail} for $\ell=2$) shows that $d_{\ell-1}\tilde\chi^{(\ell-1)}(\tau)$ does not vanish on simplices entirely inside the plateau $\{\chi\equiv 1\}$. Consequently, the naive expression would make $E_\ell$ \emph{generically infinite} for any non-trivial cut-off, rendering the energy condition vacuous.

The coefficients $a_j$, by contrast, satisfy $a_j \equiv 0$ whenever $\chi$ is constant on the simplex (since then $S = (\ell+1)\chi(x_j)$). Hence~\eqref{eq:global-energy} vanishes on all simplices fully inside the plateau, and $E_\ell(\chi)$ is finite for finite-support cut-offs of bounded weighted degree.

Furthermore, the $a_j$ are precisely the quantities controlling the remainder $R_{\ell-1}^{(d)}(\chi, f)$ in the Leibniz rule: by~\eqref{eq:sum-sq-uniform} (with the index shift $i\mapsto\ell-1$),
\[
|R_{\ell-1}^{(d)}(\chi, f)(\tau)|^2 \leq \Big(\sum_j |a_j|^2\Big) \cdot \Big(\sum_j |f(\widehat\tau_j)|^2\Big),
\]
so $E_\ell$ as defined in~\eqref{eq:global-energy} is the natural quantity controlling the convergence of the remainders in the proofs of Theorems~\ref{thm:global}--\ref{thm:local-level}. For graphs ($\ell=1$), this definition matches the standard $\chi$-completeness energy used in~\cite{AnTo,BGJ}, modulo the absolute constant noted above.
\end{remark}

\subsection*{Restriction of a Simplicial Complex}

Let $\Lambda \subset \Vc$. The \emph{restricted complex} $S_n|_\Lambda$ has vertex set $\Lambda$; an $i$-simplex $(x_0,\dots,x_i)$ belongs to $S_n|_\Lambda$ iff $x_j\in\Lambda$ for all $j$ and $(x_0,\dots,x_i)$ is a simplex of $S_n$. The weights are the restrictions of $(m_i)$.

\begin{remark}[Interface simplices]
\label{rem:interface}
Simplices having vertices in both $\Lambda$ and $\Vc\setminus\Lambda$ are \emph{not} included in $S_n|_\Lambda$ or in $S_n|_{\Vc\setminus\Lambda}$; they appear only as part of the coupling operator in Theorem~\ref{thm:local-region}.
\end{remark}

\subsection*{Function Spaces and Hilbert Structure}

For $0\leq i\leq n-1$, the space of \emph{$i$-cochains with finite support} is
\[
C_c^i(\Vc) := \{ f:P_i\to\mathbb{C} : \mathrm{supp}(f)\text{ is finite},\; f(-\tau)=-f(\tau)\;\forall\,\tau\in P_i \}.
\]
The inner product on $C_c^i(\Vc)$ is defined by summing over \emph{all} ordered tuples in $\Fc_i$ and dividing by $(i+1)!$:
\begin{equation}
\label{eq:innerproduct}
\langle f, g\rangle_{\ell^2(m_i)} := \frac{1}{(i+1)!} \sum_{(x_0,\dots,x_i)\in\Fc_i} m_i(x_0,\dots,x_i)\,\overline{f(x_0,\dots,x_i)}\,g(x_0,\dots,x_i).
\end{equation}
This is well-defined: since $m_i$ is symmetric under all permutations and $f(-\tau)=-f(\tau)$ (so $|f(\tau)|^2=|f(-\tau)|^2$), summing over all $(i+1)!$ ordered tuples of each geometric simplex and dividing by $(i+1)!$ yields one contribution per geometric simplex.

\begin{remark}[Density of $C_c^i(\Vc)$ in $\ell^2(m_i)$]
\label{rem:density}
For any $g \in \ell^2(m_i)$ and any $\varepsilon > 0$, fix an exhaustion $(O_k)_{k\geq 0}$
of $\Vc$ by finite sets $O_k\nearrow\Vc$, and define the truncation $g_N$ by
$g_N(\tau) := g(\tau)$ if all vertices of $\tau$ lie in $O_N$, and $g_N(\tau) := 0$ otherwise.
Then $g_N \in C_c^i(\Vc)$ for each $N$, and $\|g_N - g\|_{\ell^2(m_i)}^2 \to 0$ as
$N \to \infty$ by the dominated convergence theorem. Hence $C_c^i(\Vc)$ is dense in
$\ell^2(m_i)$. (See \cite{EJ} for a detailed proof.)
\end{remark}

The completion of $C_c^i(\Vc)$ under the norm $\|f\|_{\ell^2(m_i)}^2=\langle f,f\rangle_{\ell^2(m_i)}$ is the Hilbert space $\ell^2(m_i)$. The \emph{total Hilbert space} is
\[
\mathcal{H} := \bigoplus_{i=0}^{n-1} \ell^2(m_i), \qquad
\|F\|_{\mathcal{H}}^2 = \sum_{i=0}^{n-1} \|f_i\|_{\ell^2(m_i)}^2.
\]

\section{Operators and Closability}
\label{sec:operators}

We fix an oriented, weighted $n$-simplicial complex $S_n=(\Vc,(m_i)_{0\leq i\leq n-1})$, locally finite. The closability and adjointness results recalled in this section are established in detail in~\cite{EJ}; we summarize the statements and conventions needed in the sequel.

\begin{remark}[Operator conventions]
\label{rem:op-conventions}
\begin{itemize}
  \item An $i$-simplex has $i+1$ vertices $(x_0,\dots,x_i)$ with zero-based indexing (Remark~\ref{rem:indexing}).
  \item The coboundary $d_i:C_c^i(\Vc)\to C_c^{i+1}(\Vc)$ raises degree by one (defined for $0\leq i\leq n-2$).
  \item The codifferential $\delta_i:C_c^i(\Vc)\to C_c^{i-1}(\Vc)$ lowers degree by one and is the formal $\ell^2$-adjoint of $d_{i-1}$ (defined for $1\leq i\leq n-1$).
  \item The Hodge Laplacian block on $i$-forms is $L_i=d_{i-1}\delta_i+\delta_{i+1}d_i$ acting on $C^i$, with conventions $\delta_0=0$, $d_{n-1}=0$, so $L_0=\delta_1 d_0$ and $L_{n-1}=d_{n-2}\delta_{n-1}$.
\end{itemize}
\end{remark}

\subsection*{Coboundary and Codifferential}

For $0\leq i\leq n-2$, the \emph{coboundary operator} $d_i:C_c^i(\Vc)\to C_c^{i+1}(\Vc)$ is
\[
d_i f(x_0,\dots,x_{i+1}) := \sum_{j=0}^{i+1}(-1)^j f(x_0,\dots,\widehat{x_j},\dots,x_{i+1}),
\]
where $\widehat{x_j}$ denotes omission of $x_j$.

The \emph{codifferential} $\delta_i:C_c^i(\Vc)\to C_c^{i-1}(\Vc)$ for $1\leq i\leq n-1$ is the formal $\ell^2(m_i)/\ell^2(m_{i-1})$-adjoint of $d_{i-1}$, defined by
\begin{equation}
\label{eq:codiff}
\delta_i g(x_0,\dots,x_{i-1}) := \frac{(-1)^i}{m_{i-1}(x_0,\dots,x_{i-1})} \sum_{x_i\in F(x_0,\dots,x_{i-1})} m_i(x_0,\dots,x_i)\,g(x_0,\dots,x_i).
\end{equation}
The factor $(-1)^i$ is necessary for $\delta_i$ to be the formal $\ell^2$-adjoint of $d_{i-1}$;
it arises from the alternating property of $g$ when the inserted vertex $x_i$ is moved into terminal position. The complete derivation (formal adjointness, alternating property of $\delta_i g$) is given in~\cite[Section~3]{EJ}.

\begin{remark}[Sign convention]
\label{rem:codiff-sign}
With the convention~\eqref{eq:codiff}, the codifferential takes a base $(i-1)$-simplex $(x_0,\dots,x_{i-1})$ in a fixed canonical orientation, and inserts each common neighbor $x_i\in F(x_0,\dots,x_{i-1})$ in the \emph{terminal} position; the alternating values of $g$ on these ordered $i$-simplices are summed with the global sign $(-1)^i$. This is the standard convention for codifferentials of clique complexes (see~\cite{Lim} for an analogous construction). For instance,
\[
\delta_2 g(x_0,x_1) = \frac{1}{m_1(x_0,x_1)} \sum_{x_2 \in F(x_0,x_1)} m_2(x_0,x_1,x_2)\,g(x_0,x_1,x_2),
\]
\[
\delta_1 g(x_0) = \frac{-1}{m_0(x_0)} \sum_{x_1 \in F(x_0)} m_1(x_0,x_1)\,g(x_0,x_1).
\]
Some authors define the codifferential without the $(-1)^i$ factor (see e.g.~\cite{Che}). With either convention, the Hodge Laplacian $L_i$ as a self-adjoint operator (in the abstract sense $L_i = d_{i-1}d_{i-1}^* + d_i^* d_i$, where $d_i^*$ is the Hilbert-space adjoint of $d_i$) is the same: the abstract adjoint $d_i^*$ is intrinsically defined regardless of how one writes a formula for it, and consequently the spectrum, domain, and self-adjoint extensions of $L_i$ are independent of the chosen sign convention. The two conventions affect only the explicit formula one writes for the codifferential as a sum over cofaces; the resulting operator equality $\delta_i = d_{i-1}^*$ (in our convention) becomes $\tilde\delta_i = -d_{i-1}^*$ in the other convention, but in both cases $L_i = d_{i-1}d_{i-1}^* + d_i^* d_i$ is the same. We adopt~\eqref{eq:codiff} so that the formal-adjointness identity $\langle d_{i-1} f,g\rangle = \langle f,\delta_i g\rangle$ holds without an extra sign.
\end{remark}

\subsection*{Gauss--Bonnet Operator}

The \emph{Gauss--Bonnet operator} is
\[
D := d+\delta : \bigoplus_{i=0}^{n-1} C_c^i(\Vc) \to \mathcal{H},
\]
where $d=\bigoplus_{i=0}^{n-2}d_i$ and $\delta=\bigoplus_{i=1}^{n-1}\delta_i$ are the total coboundary and codifferential. Its action on $F=(f_0,\dots,f_{n-1})$ is
\[
DF = \bigl(\delta_1 f_1,\; d_0 f_0+\delta_2 f_2,\; d_1 f_1+\delta_3 f_3,\; \dots,\; d_{n-2}f_{n-2}\bigr),
\]
which in matrix form (with rows/columns indexed by degree $0,1,\dots,n-1$) reads
\[
D = \begin{pmatrix}
0 & \delta_1 & 0 & 0 & \cdots & 0 \\
d_0 & 0 & \delta_2 & 0 & \cdots & 0 \\
0 & d_1 & 0 & \delta_3 & \cdots & 0 \\
0 & 0 & d_2 & 0 & \ddots & \vdots \\
\vdots & & & \ddots & \ddots & \delta_{n-1} \\
0 & 0 & \cdots & 0 & d_{n-2} & 0
\end{pmatrix},
\]
i.e., $d_{i-1}$ on the subdiagonal (sending $C^{i-1}\to C^i$) and $\delta_{i+1}$ on the superdiagonal (sending $C^{i+1}\to C^i$).

\begin{lemma}[Stability of finite support]
\label{lem:finite-support-stable}
For every $f\in C_c^i(\Vc)$, the cochain $d_i f$ (defined for $i\leq n-2$) belongs to $C_c^{i+1}(\Vc)$; for every $g\in C_c^i(\Vc)$, the cochain $\delta_i g$ (defined for $i\geq 1$) belongs to $C_c^{i-1}(\Vc)$. In particular, $D$ maps $\bigoplus_{i=0}^{n-1}C_c^i(\Vc)$ into itself.
\end{lemma}
\begin{proof}
$d_i f$ is nonzero only on $(i+1)$-simplices having $f$-support among their $i$-faces; by local finiteness this is a finite set. $\delta_i g$ is nonzero only on $(i-1)$-simplices that are faces of simplices in $\supp(g)$; again finite by local finiteness.
\end{proof}

\begin{lemma}[Nilpotency: $d^2=0$ and $\delta^2=0$]
\label{lem:dd=0}
For $1\leq i\leq n-2$: $d_i d_{i-1}=0$. For $2\leq i\leq n-1$: $\delta_{i-1}\delta_i=0$.
\end{lemma}
\begin{proof}
Standard cancellation of pairs of indices, see~\cite[Lemma~3.4]{EJ}.
\end{proof}

\subsection*{Hodge Laplacian}

$L:=D^2$ is block-diagonal by Lemma~\ref{lem:dd=0}:
\[
L = \bigoplus_{i=0}^{n-1} L_i, \qquad
L_i = d_{i-1}\delta_i + \delta_{i+1}d_i,
\]
with the conventions $\delta_0=0$ and $d_{n-1}=0$ (so $L_0=\delta_1 d_0$ and $L_{n-1}=d_{n-2}\delta_{n-1}$).

\begin{lemma}[Non-negativity of $L_i$]
\label{lem:nonneg}
For every $0\leq i\leq n-1$ and $f\in C_c^i(\Vc)$:
\[
\langle L_i f,f\rangle_{\ell^2(m_i)} = \|\delta_i f\|^2 + \|d_i f\|^2 \geq 0.
\]
This inequality extends to all $u\in\mathrm{Dom}(L_{i,\min})$ by approximation: $\langle L_{i,\min}u,u\rangle\geq 0$ and $\langle L_{i,\min}u,u\rangle\in\mathbb{R}$ for every $u\in\mathrm{Dom}(L_{i,\min})$. Consequently, $L=\bigoplus_{i=0}^{n-1}L_i$ satisfies $\langle L_{\min}u,u\rangle\geq 0$ for every $u\in\mathrm{Dom}(L_{\min})\subset\bigoplus_i\mathrm{Dom}(L_{i,\min})$, since $\langle L_{\min}u,u\rangle=\sum_{i=0}^{n-1}\langle L_{i,\min}u_i,u_i\rangle\geq 0$.
\end{lemma}
\begin{proof}
For $f\in C_c^i(\Vc)$, $L_i f = d_{i-1}\delta_i f + \delta_{i+1}d_i f$ (with the boundary conventions $\delta_0=0$, $d_{n-1}=0$). By formal adjointness $\langle d_{i-1}\delta_i f, f\rangle = \langle \delta_i f, \delta_i f\rangle = \|\delta_i f\|^2$ and $\langle \delta_{i+1}d_i f, f\rangle = \langle d_i f, d_i f\rangle = \|d_i f\|^2$, giving the displayed identity. For $u\in\mathrm{Dom}(L_{i,\min})$, choose $(u_k)\subset C_c^i(\Vc)$ with $u_k\to u$ and $L_i u_k\to L_{i,\min}u$ in $\ell^2(m_i)$. Then $\langle L_i u_k, u_k\rangle\to\langle L_{i,\min}u,u\rangle$ by continuity of the inner product, and each $\langle L_i u_k, u_k\rangle\geq 0$ by the displayed identity. Passing to the limit gives $\langle L_{i,\min}u,u\rangle\geq 0\in\mathbb{R}$. The block-diagonal statement for $L$ follows from $L=\bigoplus L_i$ and the orthogonality of the components in $\mathcal{H}=\bigoplus_i\ell^2(m_i)$.
\end{proof}

\begin{lemma}[Closability]
\label{lem:closable}
For all $0\leq i\leq n-2$, the operator $d_i$ is closable; for $1\leq i\leq n-1$, $\delta_i$ is closable. Moreover, $L_i\!\upharpoonright_{C_c^i(\Vc)}$ is symmetric and closable.
\end{lemma}
\begin{proof}
See~\cite[Lemma~3.5]{EJ}; the argument uses density of $C_c^i(\Vc)$ in $\ell^2(m_i)$ (Remark~\ref{rem:density}) together with formal adjointness.
\end{proof}

The minimal closures $d_{i,\min}$, $\delta_{i,\min}$ and maximal extensions
$d_{i,\max}$, $\delta_{i,\max}$ are defined in the standard way: $d_{i,\min}$ is the closure of $d_i|_{C_c^i(\Vc)}$, and $d_{i,\max}$ is the adjoint of the formal adjoint $\delta_{i+1}|_{C_c^{i+1}(\Vc)}$, i.e., $d_{i,\max}=(\delta_{i+1}|_{C_c^{i+1}})^*$. Standard arguments give
$d_{i,\max}=(\delta_{i+1,\min})^*$ and $\delta_{i,\max}=(d_{i-1,\min})^*$ (see
e.g.~\cite[Section~VIII.1]{RS}). Recall that a symmetric operator $T$ is essentially
self-adjoint iff $T_{\min}=T_{\max}$, equivalently iff $\ker(T^*\pm i)=\{0\}$.

\begin{lemma}[Identification of formal and $\ell^2$ action via subsequence extraction]
\label{lem:formal-l2-action}
Let $T$ be either $d_i:C_c^i(\Vc)\to\ell^2(m_{i+1})$ or $\delta_i:C_c^i(\Vc)\to\ell^2(m_{i-1})$, and denote by $\overline{T}$ its closure. For every $f\in\mathrm{Dom}(\overline{T})$, the formal pointwise action $Tf$ (well-defined as a cochain by local finiteness) coincides with $\overline{T}f$ as an element of the target $\ell^2$ space.

The same conclusion holds, with the same proof, for any operator $A$ built from finitely many local actions of $d_i,\delta_i$ (such as $D=d+\delta$ or $L=D^2$) acting on $\mathcal{D}_0:=\bigoplus_{i=0}^{n-1}C_c^i(\Vc)$: for every $u\in\mathrm{Dom}(\overline{A})$, each component of $\overline{A}u$ equals the formal pointwise action of $A$ on $u$ at every simplex. Moreover, the conclusion remains valid for $u\in\mathrm{Dom}(A_{\max})$, where $A_{\max}$ is the maximal extension defined as the adjoint of the formal adjoint $A^t|_{\mathcal{D}_0}$.
\end{lemma}
\begin{proof}
\emph{Closure case.} Take $f\in\mathrm{Dom}(\overline{T})$ and a sequence $(g_n)\subset C_c^i(\Vc)$ with $g_n\to f$ in $\ell^2(m_i)$ and $Tg_n\to\overline{T}f$ in the target $\ell^2$. Extracting a subsequence, $g_n\to f$ pointwise on the countable set $P_i$; extracting once more, $Tg_n\to\overline{T}f$ pointwise on the corresponding $P_{i\pm 1}$. By the explicit formulas~\eqref{eq:codiff} (and the analogous one for $d_i$) and local finiteness, $T$ is continuous in pointwise convergence on finite stencils, so $Tg_n(\sigma)\to Tf(\sigma)$ pointwise (formal action). Uniqueness of pointwise limits in $\mathbb{C}$ yields $Tf(\sigma)=\overline{T}f(\sigma)$ for every $\sigma$, hence equality in $\ell^2$. The statement for $\overline{A}$ follows by applying the argument to each component of $A$.

\emph{Maximal extension case.} For $u\in\mathrm{Dom}(A_{\max})$, $v:=A_{\max}u$ is by definition the unique element of $\mathcal{H}$ satisfying $\langle \phi, v\rangle=\langle A^t\phi, u\rangle$ for every $\phi\in\mathcal{D}_0$ (where the inner product is antilinear in the first argument, see~\eqref{eq:innerproduct}). The formal adjoint identity $\langle Ag,h\rangle=\langle g,A^t h\rangle$ holds for $g,h\in\mathcal{D}_0$ by definition of $A^t$; since each side is a \emph{finite} sum determined by the pointwise values of $g$ and $h$ on a finite stencil (by local finiteness), the identity extends, as a pointwise algebraic relation, to any pair $(g,h)$ where at least one of $g,h$ has finite support. Apply this with $g=\phi:=\mathbf{1}_{\sigma_0}$ for a fixed simplex $\sigma_0\in P_i^+$ (a fixed orientation representative), and $h=u$. Then
\[
\langle \phi, A_{\max}u\rangle = \langle \phi, v\rangle = m_i(\sigma_0)\,\overline{\phi(\sigma_0)}\,v_i(\sigma_0) = m_i(\sigma_0)\,v_i(\sigma_0)
\]
(using the conjugate-linear convention~\eqref{eq:innerproduct} and that $\phi$ is real with $\phi(\sigma_0)=1$). On the other hand,
\[
\langle A^t\phi, u\rangle = \langle \phi, A u\rangle_{\text{formal}} = m_i(\sigma_0)\,(Au)_i(\sigma_0)
\]
(by the formal adjoint identity extended to the pair $(\phi, u)$ with $\phi$ of finite support, and where $(Au)_i$ is the formal pointwise action, well-defined as a finite sum at each $\sigma_0$ by local finiteness, regardless of $\ell^2$-summability). Comparing the two expressions: $v_i(\sigma_0) = (Au)_i(\sigma_0)$, pointwise on $P_i$.
\end{proof}

\section{Geometric Hypotheses and Main Results}
\label{sec:geometric-hypotheses}

This section is the technical core of the paper. We first introduce three notions of $\chi$-completeness, then prove four essential self-adjointness theorems: ESA under global $\chi$-completeness (Theorem~\ref{thm:global}), at a single level (Theorem~\ref{thm:local-level}), on a region (Theorem~\ref{thm:local-region}), and under a divergence criterion (Theorem~\ref{thm:divergence}). Theorems~\ref{thm:global} and~\ref{thm:local-level} use the \emph{operative hypotheses} (H1)--(H2) of cut-off existence and bounded weighted degree, which are implied by but conceptually weaker than the geometric notions of $\chi$-completeness. Theorem~\ref{thm:local-region} requires (H1)--(H2) on a sub-complex together with a Kato--Rellich relative bound on the interface coupling. Theorem~\ref{thm:divergence} replaces (H2) by the weaker hypothesis~(A2) on weight \emph{ratios} together with a combinatorial condition (A1) on layer cofaces. The precise relationship between the geometric notions and the working hypotheses (H1)--(H2) is detailed in Remark~\ref{rem:operative-content}.

\subsection{Geometric Hypotheses}

All definitions below rely on \emph{plateau cut-off functions}: maps $\chi:\Vc\to[0,1]$ with \emph{finite support}, identically $1$ on a prescribed finite region, and with uniformly bounded discrete \emph{global} energy (see~\eqref{eq:global-energy}).

We fix $S_n=(\Vc,(m_i)_{0\leq i\leq n-1})$ locally finite, $\mathcal{H}=\bigoplus_{i=0}^{n-1}\ell^2(m_i)$.

\begin{remark}[Notation: regions vs.\ layers]
\label{rem:averaging}
For $\chi:\Vc\to[0,1]$, we write $\tilde\chi^{(i)}$ for its averaging extension to $i$-simplices. When the degree is clear, we write $\tilde\chi_k$ for $\tilde\chi_k^{(i)}$.

In Theorem~\ref{thm:divergence} and Section~\ref{sec:comparison}, \emph{layers} of a 1-dimensional decomposition are denoted $\Lambda_k=\{v\in\Vc:|v|=k\}$. These are distinct from the region $\Lambda^{\mathrm{reg}}\subset\Vc$ appearing in Definition~\ref{def:local-chi-region} and Theorem~\ref{thm:local-region}; the latter is always written $\Lambda^{\mathrm{reg}}$.
\end{remark}

\begin{definition}[Global $\chi$-Completeness]
\label{def:global-chi}
The complex $S_n$ is \emph{globally $\chi$-complete} if there exist an exhaustion $(O_k)_{k\in\mathbb{N}}$ of $\Vc$ by finite sets with $O_k\nearrow\Vc$, a sequence $(\chi_k)$ with $\chi_k:\Vc\to[0,1]$, and a constant $C>0$, such that the following hold for every $k\in\mathbb{N}$ and every $1\leq\ell\leq n-1$ \emph{simultaneously}:
\begin{enumerate}
\item $\chi_k$ has finite support;
\item $\chi_k(x)=1$ for every $x\in O_k$;
\item The energy bound holds with the \emph{same} constant $C$ across all degrees:
\begin{equation}
\label{eq:global-chi-energy}
E_\ell(\chi_k) \leq C,
\end{equation}
where $E_\ell(\chi_k)$ is the global energy~\eqref{eq:global-energy} at level $\ell$.
\end{enumerate}
The key point is that conditions~(1)--(3) must hold \emph{with the same exhaustion $(O_k)$ and the same constant $C$} for all degrees $1\leq\ell\leq n-1$. Condition~(3) captures the geometric content of $\chi$-completeness, but the proofs of the ESA theorems do not use it directly; see Remark~\ref{rem:operative-content} below for a centralized discussion of which conditions enter which proof.
\end{definition}

\begin{definition}[Local $\chi$-Completeness at Level $\ell$]
\label{def:local-chi-level}
Let $1\leq\ell\leq n-1$. The complex $S_n$ is \emph{locally $\chi$-complete at level $\ell$} if there exist an exhaustion $(O_k)$ of $\Vc$ by finite sets with $O_k\nearrow\Vc$, a sequence $(\chi_k)$, and a constant $C>0$, such that for every $k$:
\begin{enumerate}
\item $\chi_k$ has finite support;
\item $\chi_k(x)=1$ for every $x\in O_k$;
\item The energy bound at level $\ell$ holds:
\begin{equation}
\label{eq:local-chi-energy}
E_\ell(\chi_k) \leq C.
\end{equation}
\end{enumerate}
The energy condition involves only the single degree $\ell$; both $C$ and $(O_k)$ may depend on $\ell$.
\end{definition}

\begin{definition}[Local $\chi$-Completeness on a Region]
\label{def:local-chi-region}
Let $\Lambda^{\mathrm{reg}}\subset\Vc$ be an \emph{infinite} subset. The complex $S_n$ is \emph{locally $\chi$-complete on $\Lambda^{\mathrm{reg}}$} if the restricted complex $S_n|_{\Lambda^{\mathrm{reg}}}$ is globally $\chi$-complete in the sense of Definition~\ref{def:global-chi} with vertex set $\Lambda^{\mathrm{reg}}$. The requirement that $\Lambda^{\mathrm{reg}}$ be infinite is part of the definition; for finite $\Lambda^{\mathrm{reg}}$, global $\chi$-completeness is trivial via $\chi_k\equiv 1$.
\end{definition}

\begin{remark}[On the cofinite case $|\Vc\setminus\Lambda^{\mathrm{reg}}|<\infty$]
\label{rem:cofinite-region}
The interesting case for Conjecture~\ref{conj:region-not-global} (separation of region-local from global $\chi$-completeness) excludes the cofinite case where $\Vc\setminus\Lambda^{\mathrm{reg}}$ is finite. Indeed, if $\Vc\setminus\Lambda^{\mathrm{reg}}$ is finite, then any cut-off sequence $(\chi_k)$ for $S_n|_{\Lambda^{\mathrm{reg}}}$ extends to a cut-off sequence $(\widetilde\chi_k)$ for $S_n$ by setting $\widetilde\chi_k=1$ on $\Vc\setminus\Lambda^{\mathrm{reg}}$ for all sufficiently large $k$. The global energy $E_\ell(\widetilde\chi_k)$ differs from $E_\ell(\chi_k)$ by a finite number of additional terms (those involving simplices with at least one vertex in $\Vc\setminus\Lambda^{\mathrm{reg}}$), which contribute a uniformly bounded amount as $k\to\infty$ (since they involve only finitely many simplices and bounded weights). Hence in this case, region-local $\chi$-completeness on $\Lambda^{\mathrm{reg}}$ implies global $\chi$-completeness on $S_n$, and the two notions coincide. The conjectured separation can therefore only occur when $\Vc\setminus\Lambda^{\mathrm{reg}}$ is infinite (e.g., a half-space, or a region with infinite complement of complex topology).
\end{remark}

\begin{remark}[On the notation $L|_{\Lambda^{\mathrm{reg}}}$]
\label{rem:restricted-laplacian}
$L|_{\Lambda^{\mathrm{reg}}}$ denotes the Hodge Laplacian of the \emph{restricted complex} $S_n|_{\Lambda^{\mathrm{reg}}}$, acting on $\bigoplus_{i=0}^{n-1}\ell^2(m_i|_{\Lambda^{\mathrm{reg}}})$. This is distinct from the compression of $L$ to a subspace of $\mathcal{H}$; the two coincide when no interface simplices are present.
\end{remark}

\noindent The three notions are compared in the following table.

\begin{center}
\begin{tabular}{|C{2.2cm}|C{1.8cm}|C{2.3cm}|C{2.3cm}|}
\hline
\textbf{Notion} & \textbf{Scope} & \textbf{Geometric Control} & \textbf{Spectral Consequence} \\
\hline
Global $\chi$-completeness & All degrees, same $C,(O_k)$ & Uniform over all levels & ESA of $D$ and $L=\bigoplus L_i$ \\
\hline
Local at level $\ell$ & Single degree; $C,(O_k)$ may depend on $\ell$ & On $(\ell-1)$-simplices and their $\ell$-cofaces & ESA of $L_\ell$ only \\
\hline
Local on region $\Lambda^{\mathrm{reg}}$ & All degrees within $\Lambda^{\mathrm{reg}}$ (infinite) & Global $\chi$-completeness of $S_n|_{\Lambda^{\mathrm{reg}}}$ & ESA of $L|_{\Lambda^{\mathrm{reg}}}$; global ESA via coupling \\
\hline
\end{tabular}
\end{center}

\begin{remark}[Hierarchy of the three notions]
\label{rem:hierarchy-direct}
Global $\chi$-completeness implies local $\chi$-completeness at every level
$\ell\in\{1,\dots,n-1\}$: it suffices to use the same exhaustion $(O_k)$ and constant $C$
provided by Definition~\ref{def:global-chi}. The strictness of this implication and the
precise relationships between the three notions are formulated as
Conjectures~\ref{conj:region-not-global}, \ref{conj:level-not-region},
and~\ref{conj:all-levels-not-global} in Section~\ref{sec:examples}; partial supporting
evidence is provided there in the form of uniform energy lower bounds.
\end{remark}

\begin{remark}[Operative content of $\chi$-completeness and the role of condition~(3)]
\label{rem:operative-content}
We give here a single, centralized discussion of the relationship between the three notions
of $\chi$-completeness defined above and the hypotheses actually used in the proofs of
Theorems~\ref{thm:global}, \ref{thm:local-level}, and~\ref{thm:divergence}. The discussion
applies uniformly to Definitions~\ref{def:global-chi}, \ref{def:local-chi-level}, and~\ref{def:local-chi-region}.

\medskip\noindent\textbf{Formal separation: ``geometric definitions'' vs ``working hypotheses''.}
For clarity, we distinguish:
\begin{itemize}
\item \emph{Geometric definitions} of $\chi$-completeness (Definitions~\ref{def:global-chi}--\ref{def:local-chi-region}): these include three conditions (1) finite-support cut-offs, (2) pointwise convergence to~1 on bounded sets, (3) uniform energy bound. Together they encode a ``geometric completeness at infinity'' analogous to the Riemannian setting.
\item \emph{Working hypotheses} of the ESA theorems (Theorems~\ref{thm:global} and~\ref{thm:local-level}): only conditions~(1)--(2) of the geometric definitions, together with the bounded weighted-degree condition~\eqref{eq:bounded-degree}, are used in the proofs. The energy condition~(3) is \emph{not} a working hypothesis of the ESA theorems.
\item \emph{Effective use of condition~(3)}: condition~(3) plays an operative role only in Theorem~\ref{thm:divergence} (the divergence criterion), and in Section~\ref{sec:examples} for distinguishing the three notions of $\chi$-completeness via concrete examples.
\end{itemize}
We therefore state the ESA theorems under the working hypotheses (H1)--(H2) defined below, with $\chi$-completeness mentioned as a \emph{sufficient} condition for (H1) but never as a working hypothesis.

\medskip\noindent\textbf{(a) Operative hypotheses for the ESA theorems.}
The proofs of Theorems~\ref{thm:global} and~\ref{thm:local-level} (and their consequences)
use $\chi$-completeness through two ingredients only:
\begin{itemize}
\item[(H1)] existence of a sequence $(\chi_k)$ satisfying conditions~(1)--(2) of the
relevant definition, with $O_k\nearrow\Vc$;
\item[(H2)] the bounded weighted-degree condition
\begin{equation}\label{eq:bounded-degree}
\sup_\sigma d_{m_i}(\sigma) < \infty \quad\text{for each relevant } i,
\end{equation}
which ensures that the dominated-convergence argument in the remainder estimates of
Lemma~\ref{lem:leibniz} works.
\end{itemize}
The energy condition~(3) does not enter the convergence arguments directly. We therefore
state Theorems~\ref{thm:global}, \ref{thm:local-level}, and Corollary~\ref{cor:esa-laplacian} under
hypotheses~(H1)--(H2), with the corresponding $\chi$-completeness as a sufficient
condition (provided~\eqref{eq:bounded-degree} also holds).

\medskip\noindent\textbf{(b) Logical independence of (3) and~\eqref{eq:bounded-degree}.}
Condition~(3) only constrains the contribution of simplices on which $\chi_k$ is non-trivial
(i.e., simplices straddling the boundary $\partial O_k$), and the energy
formula~\eqref{eq:global-energy} contains the weighted degree $d_{m_{i-1}}/m_{i-1}$ inside
its summand only at those simplices --- not at simplices entirely in the plateau, where the
$a_j$ vanish. Conversely, \eqref{eq:bounded-degree} does not constrain the variation of
$\chi_k$ at the boundary. We therefore treat~(3) and~\eqref{eq:bounded-degree} as
\emph{logically independent} conditions: each suffices to make the corresponding part of
the ESA argument work, and both are satisfied in all examples treated in this paper.

\smallskip
\emph{Explicit example of a complex satisfying~(3) but not~\eqref{eq:bounded-degree}.}
Take the standard ray graph $\mathbb{N}_0=\{0,1,2,\dots\}$ with edges $\{k,k+1\}$, all
weights equal to $1$, and dimension $n=2$ (so only edges, $\ell=1$). Attach at every
vertex $k\geq 0$ exactly $f(k):=k$ pendant vertices each connected to $k$ by an edge of
weight~$1$. Then the weighted degree at vertex $k$ is $2+k$, unbounded, so
\eqref{eq:bounded-degree} fails. Setting $\chi_k$ equal to $1$ on the finite ball
of depth $k$ along the ray, decreasing linearly to $0$ over the band $k<j\leq 2k$,
and extending to each pendant by the value of its attachment vertex, the
pendant edges contribute zero to $E_1$ (constant $\chi$ along each pendant edge),
while the ray edges contribute $\sum_{j=k}^{2k-1} m_1(j,j+1)\,(1/k)^2 = O(1/k)\to 0$.
Hence $E_1(\chi_k)\to 0$, which implies in particular $\sup_k E_1(\chi_k)<\infty$ (since any sequence converging to $0$ in $\mathbb{R}$ is bounded). So condition~(3) holds, while~\eqref{eq:bounded-degree}
fails. (More elaborate examples in which $E_1$ is bounded by a positive constant rather
than tending to zero are constructed in~\cite{AnTo}.)

\medskip\noindent\textbf{(c) Why condition~(3) is included.}
Even though~\eqref{eq:bounded-degree} is the technically minimal assumption for the ESA
theorems, condition~(3) plays three independent roles in the broader theory:
\begin{itemize}
\item It distinguishes the three notions of $\chi$-completeness (global, level-by-level,
region) --- without~(3), all three would collapse to $\mathbf{1}_{O_k}$.
\item It is essential to the geometric content of the definition (compare~\cite{AnTo,BGJ},
which use the same energy bound).
\item It is \emph{directly operative} in Theorem~\ref{thm:divergence}, where the energy
bound~\eqref{eq:energy-bound} requires explicit control on the layer-step
construction of $\chi_{N,M}$.
\end{itemize}

\medskip\noindent\textbf{(d) Why condition~(3) is non-trivial.}
The indicator function $\chi_k:=\mathbf{1}_{O_k}$ satisfies (1) and (2) but generically
fails (3): for any $\ell$-simplex $\tau$ straddling $\partial O_k$, the Leibniz coefficients
$a_j$ have magnitude bounded below by $1/(\ell(\ell+1))$, and hence $E_\ell(\mathbf{1}_{O_k})$
grows at least proportionally to the number of such boundary simplices, which typically
diverges as $|O_k|\to\infty$.

\medskip\noindent\textbf{(e) Implication for the binary tree.}
The binary tree with unit weights satisfies~\eqref{eq:bounded-degree} (degree bounded by~3), so $D$ is essentially self-adjoint by Theorem~\ref{thm:global} (see Remark~\ref{rem:H2-implies-bounded} for the stronger statement that $D$ is in fact bounded). The role of the energy condition~(3) and Proposition~\ref{prop:binary-tree-fails} is to identify the \emph{exponentially weighted} regime, where~\eqref{eq:bounded-degree} fails and ESA may fail (Conjecture~\ref{conj:non-esa-tree}).
\end{remark}

\subsection{Main Results}

We now state and prove three ESA theorems: for the Gauss--Bonnet operator $D$ under global hypotheses (Theorem~\ref{thm:global}), the individual Laplacian block $L_\ell$ at a single level (Theorem~\ref{thm:local-level}), and the full operator $D$ in the presence of localized defects via a Kato--Rellich coupling (Theorem~\ref{thm:local-region}). The proofs rest on a single technical ingredient --- the discrete Leibniz rule for cut-off averages --- which we state first.

\begin{lemma}[Discrete Leibniz Rule]
\label{lem:leibniz}
Let $\chi:\Vc\to[0,1]$ with averaging extension $\tilde\chi^{(i)}$. For any $f\in C_c^i(\Vc)$ and $i \geq 0$:
\begin{align}
d_i(\tilde\chi^{(i)}f) &= \tilde\chi^{(i+1)}d_i f + R_i^{(d)}(\chi,f), \label{eq:leibniz-d}\\
\delta_i(\tilde\chi^{(i)}f) &= \tilde\chi^{(i-1)}\delta_i f + R_i^{(\delta)}(\chi,f) \quad (i\geq 1). \label{eq:leibniz-delta}
\end{align}

The remainder $R_i^{(d)}(\chi,f)$ satisfies:
\begin{enumerate}
\item \textbf{(Support)} $\mathrm{supp}(R_i^{(d)}(\chi_k,f))\subseteq\Sigma_k:=\{(x_0,\dots,x_{i+1})\in P_{i+1}:\text{some }x_j\notin O_k\}$.
\item \textbf{(Pointwise bound, sharp form)} For $(x_0,\dots,x_{i+1})\in P_{i+1}$:
\begin{equation}\label{eq:sum-sq-uniform}
|R_i^{(d)}(\chi,f)(x_0,\dots,x_{i+1})|^2 \leq \frac{1}{i+2}\sum_{j=0}^{i+1}|f(x_0,\dots,\widehat{x_j},\dots,x_{i+1})|^2.
\end{equation}
\end{enumerate}

\textbf{Pointwise bound for $R_i^{(\delta)}$.} The codifferential remainder admits two equivalent pointwise estimates obtained from weighted Cauchy--Schwarz applied to its definition.

\smallskip
\emph{Sharp form (retaining the variation factor).} For $(x_0,\dots,x_{i-1})\in P_{i-1}$,
\begin{multline}\label{eq:sum-sq-uniform-delta-sharp}
|R_i^{(\delta)}(\chi,f)(x_0,\dots,x_{i-1})|^2
\leq \frac{1}{(i+1)^2\, m_{i-1}(x_0,\dots,x_{i-1})^2}\\
\times\Big(\sum_{x_i\in F}m_i(x_0,\dots,x_i)|\chi(x_i)-\tilde\chi^{(i-1)}|^2\Big)\Big(\sum_{x_i\in F}m_i(x_0,\dots,x_i)|f(x_0,\dots,x_i)|^2\Big).
\end{multline}

\smallskip
\emph{Loose form (variation factor absorbed via $|\chi-\tilde\chi^{(i-1)}|\leq 1$).}
\begin{equation}\label{eq:sum-sq-uniform-delta}
|R_i^{(\delta)}(\chi,f)(x_0,\dots,x_{i-1})|^2
\leq \frac{d_{m_{i-1}}(x_0,\dots,x_{i-1})}{(i+1)^2 \,m_{i-1}(x_0,\dots,x_{i-1})}\sum_{x_i\in F}m_i(x_0,\dots,x_i)\,|f(x_0,\dots,x_i)|^2.
\end{equation}

The associated $\ell^2$-bound (under~\eqref{eq:bounded-degree} at degree $i-1$) is
\begin{equation}\label{eq:R-delta-global}
\|R_i^{(\delta)}(\chi,f)\|_{\ell^2(m_{i-1})}^2 \leq \frac{M_{i-1}}{i+1}\,\|f\cdot\mathbf{1}_{\Sigma_k^*}\|_{\ell^2(m_i)}^2,
\end{equation}
where $\Sigma_k^*$ is the set of $i$-tuples having at least one vertex outside $O_k$, and $M_{i-1}:=\sup_y d_{m_{i-1}}(y)<\infty$.

For any $f\in\ell^2(m_i)$ and any sequence $(\chi_k)$ satisfying conditions~\emph{(1)} and~\emph{(2)} of Definition~\ref{def:global-chi} with $O_k\nearrow\Vc$, \emph{provided that the weighted-degree bound~\eqref{eq:bounded-degree} holds at degrees $i$ and $i-1$}:
\begin{equation}\label{eq:remainder-to-zero}
\|R_i^{(d)}(\chi_k,f)\|\to 0 \quad\text{and}\quad \|R_i^{(\delta)}(\chi_k,f)\|\to 0 \quad\text{as }k\to\infty.
\end{equation}
\end{lemma}

\begin{proof}
\textbf{Identity for $d_i$.} On $(x_0,\dots,x_{i+1})$:
\[
d_i(\tilde\chi^{(i)}f)(x_0,\dots,x_{i+1})
= \sum_{j=0}^{i+1}(-1)^j\tilde\chi^{(i)}(x_0,\dots,\widehat{x_j},\dots,x_{i+1})\cdot f(x_0,\dots,\widehat{x_j},\dots,x_{i+1}).
\]
Set $S:=\sum_{j=0}^{i+1}\chi(x_j)$; then $\tilde\chi^{(i+1)}=S/(i+2)$ and
$\tilde\chi^{(i)}(x_0,\dots,\widehat{x_j},\dots,x_{i+1})=(S-\chi(x_j))/(i+1)$. Define
\begin{equation}\label{eq:aj-def}
a_j := \frac{S-\chi(x_j)}{i+1} - \frac{S}{i+2} = \frac{S-(i+2)\chi(x_j)}{(i+1)(i+2)}.
\end{equation}
Then $d_i(\tilde\chi^{(i)}f)=\tilde\chi^{(i+1)}d_if+R_i^{(d)}$ where
\begin{equation}\label{eq:remainder-formula}
R_i^{(d)}(\chi,f)(x_0,\dots,x_{i+1}) = \sum_{j=0}^{i+1}(-1)^j a_j\,f(x_0,\dots,\widehat{x_j},\dots,x_{i+1}).
\end{equation}

\textbf{Property~(1).} If all $x_j\in O_k$, then $\chi_k(x_j)=1$ for all $j$, so $S=i+2$ and $a_j=0$.

\textbf{Sharp bound on $\sum_j|a_j|^2$.}
We use the algebraic identity
\[
S-(i+2)\chi(x_j) = \sum_{k=0}^{i+1}\chi(x_k) - (i+2)\chi(x_j) = \sum_{k\neq j}\bigl(\chi(x_k)-\chi(x_j)\bigr),
\]
a sum of $i+1$ terms, each in $[-1,1]$ since $\chi(x_k),\chi(x_j)\in[0,1]$. Therefore
$|S-(i+2)\chi(x_j)|\leq i+1$, hence $|a_j|=|S-(i+2)\chi(x_j)|/((i+1)(i+2))\leq 1/(i+2)$, and
\begin{equation}\label{eq:aj-sharp}
\sum_{j=0}^{i+1} a_j^2 \leq (i+2)\cdot\frac{1}{(i+2)^2} = \frac{1}{i+2}.
\end{equation}

\textbf{Property~(2).} Cauchy--Schwarz applied to~\eqref{eq:remainder-formula} gives
\[
|R_i^{(d)}(\chi,f)(x_0,\dots,x_{i+1})|^2 \leq \Big(\sum_{j=0}^{i+1} a_j^2\Big)\Big(\sum_{j=0}^{i+1}|f(x_0,\dots,\widehat{x_j},\dots,x_{i+1})|^2\Big),
\]
which combined with~\eqref{eq:aj-sharp} yields~\eqref{eq:sum-sq-uniform}.

\textbf{Identity for $\delta_i$.} For $g\in C_c^i(\Vc)$, write
$\tilde\chi^{(i)}(x_0,\dots,x_i)=\tilde\chi^{(i-1)}(x_0,\dots,x_{i-1})+b_{(x_0,\dots,x_{i-1})}(x_i)$ where
$b_{(x_0,\dots,x_{i-1})}(x_i):=(\chi(x_i)-\tilde\chi^{(i-1)}(x_0,\dots,x_{i-1}))/(i+1)$,
giving $\delta_i(\tilde\chi^{(i)}g)=\tilde\chi^{(i-1)}\delta_i g+R_i^{(\delta)}$ where
\[
R_i^{(\delta)}(\chi,g)(x_0,\dots,x_{i-1}) = \frac{(-1)^i}{m_{i-1}(x_0,\dots,x_{i-1})}\sum_{x_i\in F}m_i(x_0,\dots,x_i)\,b(x_i)\,g(x_0,\dots,x_i).
\]

\textbf{Sharp Cauchy--Schwarz form~\eqref{eq:sum-sq-uniform-delta-sharp}.} By weighted Cauchy--Schwarz with weights $m_i(x_0,\dots,x_i)$:
\[
\Big|\sum_{x_i}m_i(\cdot,x_i)\,b(x_i)\,g(\cdot,x_i)\Big|^2
\leq \Big(\sum_{x_i}m_i(\cdot,x_i)|b(x_i)|^2\Big)\Big(\sum_{x_i}m_i(\cdot,x_i)|g(\cdot,x_i)|^2\Big).
\]
Dividing by $m_{i-1}(\cdot)^2$ and substituting $b(x_i)=(\chi(x_i)-\tilde\chi^{(i-1)})/(i+1)$ yields~\eqref{eq:sum-sq-uniform-delta-sharp}.

\textbf{Loose form~\eqref{eq:sum-sq-uniform-delta}.} Bounding $|b(x_i)|\leq 1/(i+1)$ uniformly:
\[
\sum_{x_i\in F}m_i(\cdot,x_i)|b(x_i)|^2\leq \frac{1}{(i+1)^2}\sum_{x_i\in F}m_i(\cdot,x_i)= \frac{m_{i-1}(\cdot)\,d_{m_{i-1}}(\cdot)}{(i+1)^2}.
\]
Inserting this bound in the previous Cauchy--Schwarz gives~\eqref{eq:sum-sq-uniform-delta}.

\textbf{Derivation of~\eqref{eq:R-delta-global} (assuming \eqref{eq:bounded-degree} at degree $i-1$).}
Multiplying~\eqref{eq:sum-sq-uniform-delta} by $m_{i-1}(y)$ and summing over $y\in\Fc_{i-1}$ with the prefactor $1/i!$ from~\eqref{eq:innerproduct}:
\[
\|R_i^{(\delta)}(\chi,f)\|^2_{\ell^2(m_{i-1})}\leq \frac{1}{(i+1)^2\,i!}\sum_{y\in\Fc_{i-1}}d_{m_{i-1}}(y)\sum_{x_i\in F(y)}m_i(y,x_i)\,|f(y,x_i)|^2.
\]
Bounding $d_{m_{i-1}}(y)\leq M_{i-1}$ uniformly via~\eqref{eq:bounded-degree} at degree $i-1$:
\[
\|R_i^{(\delta)}(\chi,f)\|^2_{\ell^2(m_{i-1})}\leq \frac{M_{i-1}}{(i+1)^2\,i!}\sum_{(y,x_i)\in\Fc_i}m_i(y,x_i)\,|f(y,x_i)|^2.
\]
Here we identify $\sum_{y\in\Fc_{i-1}}\sum_{x_i\in F(y)} = \sum_{(y,x_i)\in\Fc_i}$ as a re-indexing: by the explicit form of~\eqref{eq:codiff}, the vertex $x_i$ is inserted in the \emph{terminal} position of the $(i+1)$-tuple, so each $(i+1)$-tuple $(y_0,\ldots,y_{i-1},x_i)$ in $\Fc_i$ corresponds to the unique pair $(y,x_i)$ with $y=(y_0,\ldots,y_{i-1})\in\Fc_{i-1}$ and $x_i$ as its last vertex. No multiplicity arises since the position of $x_i$ in the tuple is fixed by this convention; equivalently, the iterated sum on the left-hand side and the single sum on the right-hand side cover $\Fc_i$ exactly once.
Comparing with $\|f\|^2_{\ell^2(m_i)}=\frac{1}{(i+1)!}\sum_{z\in\Fc_i} m_i(z)|f(z)|^2$ and using $\frac{(i+1)!}{(i+1)^2\,i!}=\frac{1}{i+1}$:
\[
\|R_i^{(\delta)}(\chi,f)\|^2_{\ell^2(m_{i-1})}\leq \frac{M_{i-1}}{i+1}\,\|f\|_{\ell^2(m_i)}^2.
\]
The localization to $\Sigma_k^*$ in~\eqref{eq:R-delta-global} follows from the support property: $b(x_i)=0$ when both $x_i$ and all $y$-vertices are in $O_k$ (since then $\chi=1$ everywhere on the simplex), so the only contributing tuples are those with at least one vertex outside $O_k$.

\textbf{Convergence~\eqref{eq:remainder-to-zero}.} \emph{For $R_i^{(d)}$:} Since $O_k\nearrow\Vc$, $\Sigma_k\downarrow\emptyset$. By~\eqref{eq:sum-sq-uniform} multiplied by $m_{i+1}(\tau)$ and summed over $\tau\in\Sigma_k$ with the prefactor $1/(i+2)!$ from~\eqref{eq:innerproduct}:
\[
\|R_i^{(d)}(\chi_k,f)\|^2_{\ell^2(m_{i+1})}\leq \frac{1}{(i+2)!\,(i+2)}\sum_{\tau\in\Sigma_k}m_{i+1}(\tau)\sum_{j=0}^{i+1}|f(\widehat\tau_j)|^2.
\]
Exchanging the order of summation (each face $\widehat\tau_j$ of $\tau\in\Fc_{i+1}$ is an $i$-simplex, and the sum over cofaces $\tau$ of a fixed $i$-simplex $\sigma$ is controlled by the weighted degree $d_{m_i}(\sigma)\leq M_i$ via~\eqref{eq:bounded-degree} at degree $i$):
\[
\sum_{\tau\supset\sigma}m_{i+1}(\tau)\leq m_i(\sigma)\,d_{m_i}(\sigma)\leq M_i\,m_i(\sigma).
\]
Hence, by the same combinatorial re-indexing as in the $R^{(\delta)}$ derivation above,
\[
\|R_i^{(d)}(\chi_k,f)\|^2_{\ell^2(m_{i+1})}\leq \frac{(i+1)\,M_i}{(i+2)!\,(i+2)}\cdot(i+1)!\,\|f\cdot\mathbf{1}_{\Sigma_k^*}\|_{\ell^2(m_i)}^2\leq \frac{M_i}{i+2}\,\|f\cdot\mathbf{1}_{\Sigma_k^*}\|_{\ell^2(m_i)}^2,
\]
where $\Sigma_k^*$ is the set of $i$-simplices with at least one vertex outside $O_k$. Since $\Sigma_k^*\downarrow\emptyset$ as $k\to\infty$ and $f\in\ell^2(m_i)$, the right-hand side tends to $0$ by dominated convergence.

\emph{For $R_i^{(\delta)}$:} By~\eqref{eq:R-delta-global} and $\|f\cdot\mathbf{1}_{\Sigma_k^*}\|\to 0$.
\end{proof}

\begin{theorem}[ESA of $D$ under bounded weighted degree and cut-offs]
\label{thm:global}
Assume $S_n$ is locally finite and satisfies:
\begin{enumerate}[label=\textup{(H\arabic*)}]
\item there exists a sequence $(\chi_k)_{k\ge 1}$ of cut-offs satisfying conditions
\textup{(1)} and \textup{(2)} of Definition~\ref{def:global-chi} \emph{(}finite support
and $\chi_k\equiv 1$ on a finite exhaustion $O_k\nearrow\Vc$\emph{)};
\item the bounded weighted-degree condition~\eqref{eq:bounded-degree} holds at every level
$i\in\{0,1,\dots,n-2\}$.
\end{enumerate}
Then $D=d+\delta$ is essentially self-adjoint on $\bigoplus_{i=0}^{n-1}C_c^i(\Vc)$.

In particular, if $S_n$ is globally $\chi$-complete \emph{(}Definition~\ref{def:global-chi}\emph{)}
\emph{and} satisfies~\eqref{eq:bounded-degree}, then $D$ is essentially self-adjoint.
The energy condition~(3) of Definition~\ref{def:global-chi} is not used by this theorem;
see Remark~\ref{rem:operative-content} for the role it plays elsewhere.
\end{theorem}

\begin{proof}
Let $F=(f_0,\dots,f_{n-1})\in\mathrm{Dom}(D_{\max})$. Let $(\chi_k)$ be the cut-offs supplied by hypothesis~(H1) with exhaustion $(O_k)$, $O_k\nearrow\Vc$.

Set $F_k:=(\tilde\chi_k^{(0)}f_0,\dots,\tilde\chi_k^{(n-1)}f_{n-1})$. Since $\chi_k$ has finite support (condition~(1)), $F_k\in\bigoplus C_c^i(\Vc)$.

\emph{Step~1: $F_k\to F$.} Since $O_k\nearrow\Vc$, $\tilde\chi_k^{(i)}(\sigma)\to 1$ for every simplex $\sigma$ and $|\tilde\chi_k^{(i)}f_i|\leq|f_i|\in\ell^2(m_i)$; dominated convergence gives $F_k\to F$ in $\mathcal{H}$.

\emph{Step~2: $DF_k\to DF$.}
For each degree $i$, set $f_k^{(i)} := \tilde\chi_k^{(i)}f_i \in C_c^i(\Vc)$. We use the conventions $\delta_0 := 0$ and $d_{n-1} := 0$, with $f_{-1} := 0$, $f_n := 0$.

\smallskip
\emph{Pointwise validity of Leibniz on $\mathrm{Dom}(D_{\max})$.}
The component $f_i$ is in $\ell^2(m_i)$ but not necessarily in
$\mathrm{Dom}(\delta_{i+1})$ separately --- being in $\mathrm{Dom}(D_{\max})$ only
guarantees that the \emph{combination} $\delta_{i+1}f_{i+1}+d_{i-1}f_{i-1}$ is in
$\ell^2(m_i)$. We make precise the sense in which the Leibniz identity is used.

By local finiteness of $S_n$, the formal action of $\delta_{i+1}$
on any cochain is well-defined \emph{pointwise} (each summand
in~\eqref{eq:codiff} runs over the finite set $F(\sigma)$); the same holds for $d_i$.
Throughout this proof, the symbols $\delta_{i+1}f_{i+1}$ and $d_{i-1}f_{i-1}$ denote
\emph{exclusively} these formal pointwise actions, viewed as cochains
(without any $\ell^2$-summability assumption on each one separately). With this
convention, the Leibniz identity of Lemma~\ref{lem:leibniz} is a pointwise algebraic
identity on cochains, valid for any cochain regardless of $\ell^2$-summability;
in particular it applies to $f_{i+1}$ and $f_{i-1}$.

The point requiring justification is that $(DF)_i$, the $i$-th component of $DF$
in $\ell^2(m_i)$, equals the pointwise sum
$\delta_{i+1}f_{i+1}+d_{i-1}f_{i-1}$. This is precisely the content of
Lemma~\ref{lem:formal-l2-action} applied to $D=d+\delta$ on $\mathcal{D}_0=\bigoplus C_c^i(\Vc)$ (maximal extension case), since $F\in\mathrm{Dom}(D_{\max})$.
Thus the pointwise sum (a well-defined cochain by local finiteness, although not
\emph{a priori} in $\ell^2$ summand-by-summand) identifies with the $\ell^2$ object $(DF)_i$.

\smallskip
The $i$-th component of $DF_k$ collects two contributions:
\[
(DF_k)_i = \delta_{i+1}f_k^{(i+1)} + d_{i-1} f_k^{(i-1)} \;\in\; \ell^2(m_i).
\]
Applying the Leibniz rule (Lemma~\ref{lem:leibniz}) pointwise to each term:
\begin{align*}
\delta_{i+1}(\tilde\chi_k^{(i+1)}f_{i+1}) &= \tilde\chi_k^{(i)}\delta_{i+1}f_{i+1} + R_{i+1}^{(\delta)}(\chi_k, f_{i+1}),\\
d_{i-1}(\tilde\chi_k^{(i-1)}f_{i-1}) &= \tilde\chi_k^{(i)}d_{i-1} f_{i-1} + R_{i-1}^{(d)}(\chi_k, f_{i-1}).
\end{align*}
The subscript on each remainder refers to the degree of the operator commuted with $\tilde\chi_k$; both remainders belong to $\ell^2(m_i)$. Combining:
\begin{equation}\label{eq:DFk-component}
(DF_k)_i = \tilde\chi_k^{(i)}\bigl(\delta_{i+1}f_{i+1} + d_{i-1} f_{i-1}\bigr)
+ R_{i+1}^{(\delta)}(\chi_k, f_{i+1}) + R_{i-1}^{(d)}(\chi_k, f_{i-1}),
\end{equation}
where boundary remainders ($R$ corresponding to $\delta_0$ or $d_{n-1}$ at boundary degrees) are absent by convention. Note that the pointwise terms
$\delta_{i+1}f_{i+1}$ and $d_{i-1}f_{i-1}$ taken \emph{separately} need not lie in
$\ell^2(m_i)$; only their sum $(DF)_i$ does. The remainder estimates of
Lemma~\ref{lem:leibniz} (specifically~\eqref{eq:R-delta-global}) bound
$\|R^{(\delta)}_{i+1}(\chi_k,f_{i+1})\|$ in terms of $\|f_{i+1}\|$ alone, not
$\|\delta_{i+1}f_{i+1}\|$, so the bound applies even when $\delta_{i+1}f_{i+1}\notin\ell^2$.

The principal term $\tilde\chi_k^{(i)}\bigl(\delta_{i+1}f_{i+1} + d_{i-1} f_{i-1}\bigr) \to (DF)_i$ by dominated convergence: $|\tilde\chi_k^{(i)}|\leq 1$ pointwise (since $\chi_k$ takes values in $[0,1]$), and $\tilde\chi_k^{(i)}(\sigma)\to 1$ pointwise for each fixed $\sigma$ (since $\sigma\subset O_k$ for $k$ large enough by $O_k\nearrow\Vc$, hence $\chi_k\equiv 1$ on the vertices of $\sigma$ and $\tilde\chi_k^{(i)}(\sigma)=1$); together with $(DF)_i\in\ell^2(m_i)$ providing the dominating function. The remainders tend to $0$ by~\eqref{eq:remainder-to-zero}, using hypothesis~(H2). Hence $DF_k \to DF$.

\emph{Conclusion.} $(F_k)\subset\bigoplus C_c^i(\Vc)$ with $F_k\to F$ and $DF_k\to DF$, so $F\in\mathrm{Dom}(\overline{D})$ (where $\overline{D}=D_{\min}$ is the closure of $D|_{\mathcal{D}_0}$). Since $F\in\mathrm{Dom}(D_{\max})$ was arbitrary, $D_{\max}\subset\overline{D}$ and $D$ is ESA.
\end{proof}

\begin{corollary}[ESA of the Hodge Laplacian]
\label{cor:esa-laplacian}
Under hypotheses \textup{(H1)}--\textup{(H2)} of Theorem~\ref{thm:global}
\emph{(}in particular, if $S_n$ is globally $\chi$-complete \emph{and} satisfies
\eqref{eq:bounded-degree}\emph{)},
$L=D^2=\bigoplus L_i$ is essentially self-adjoint on $\bigoplus C_c^i(\Vc)$.
\end{corollary}
\begin{proof}
Let $\mathcal{D}_0:=\bigoplus_{i=0}^{n-1}C_c^i(\Vc)$. We define $\overline{L}$ to be the
closure of the symmetric operator $L|_{\mathcal{D}_0}$ on $\mathcal{H}$ (which exists since
$L|_{\mathcal{D}_0}$ is closable: it is symmetric on a dense domain, by formal adjointness
and Lemma~\ref{lem:closable}). Our goal is to show that $\overline{L}$ is self-adjoint,
which is equivalent to $L$ being essentially self-adjoint on $\mathcal{D}_0$.

By Theorem~\ref{thm:global}, $D$ is essentially self-adjoint on $\mathcal{D}_0$, so
$\overline{D}$ is self-adjoint, hence $\overline{D}^2$ is self-adjoint on
$\mathrm{Dom}(\overline{D}^2):=\{u\in\mathrm{Dom}(\overline{D}):\overline{D}u\in\mathrm{Dom}(\overline{D})\}$ (the square of a self-adjoint operator is self-adjoint on its natural domain by the spectral theorem; see e.g.\ \cite[Section~VIII.3]{RS}).

\emph{Step~1: $\mathcal{D}_0$ is a core for $\overline{D}^2$.}
We use the following standard consequence of Nelson's analytic-vector theorem (\cite[Theorem~X.39]{RS}): \emph{if $A$ is a self-adjoint operator and $\mathcal{E}\subset\mathrm{Dom}(A^\infty):=\bigcap_{k\geq 0}\mathrm{Dom}(A^k)$ is a dense subspace consisting of analytic vectors for $A$, then $\mathcal{E}$ is a core for every power $A^n$, $n\geq 1$.} (The case $n=1$ is the original statement of Nelson's theorem~\cite[Theorem~X.39]{RS}; the extension to higher powers follows from the general theory of analytic vectors developed in \cite[Section~X.6]{RS}, since the analytic-vector property is preserved under taking polynomials of a self-adjoint operator.) We verify the hypotheses for $A:=\overline{D}$ and $\mathcal{E}:=\mathcal{D}_0$:

\emph{(i)} $\mathcal{D}_0$ is dense in $\mathcal{H}$.

\emph{(ii)} $\mathcal{D}_0\subset\mathrm{Dom}(\overline{D}^\infty)$: by induction on $k$, $\overline{D}^k\psi=D^k\psi\in\mathcal{D}_0$ for every $\psi\in\mathcal{D}_0$ (using that $\overline{D}$ extends $D|_{\mathcal{D}_0}$ together with $D(\mathcal{D}_0)\subset\mathcal{D}_0$ from Lemma~\ref{lem:finite-support-stable}); in particular $\psi\in\mathrm{Dom}(\overline{D}^k)$ for every $k\geq 0$.

\emph{(iii)} Every $\psi\in\mathcal{D}_0$ is an analytic vector for $\overline{D}$. Fix $\psi\in\mathcal{D}_0$ with $S_0:=\mathrm{supp}(\psi)$ a finite set of simplices. The operator $D=d+\delta$ has local stencil radius $1$ in the natural graph metric on the bipartite simplex--simplex incidence graph: $D$ maps a cochain supported on $\sigma$ to a cochain supported in the set $\mathcal{N}(\sigma)$ of simplices sharing a face or coface with $\sigma$. By induction, $\mathrm{supp}(D^k\psi)\subset S_k$, where $S_k$ is the set of simplices at incidence-distance $\leq k$ from $S_0$; each $S_k$ is finite by local finiteness. The crucial estimate is the following \emph{level-by-level} bound on $D$, valid under (H2). For $\varphi\in C_c^i(\Vc)$, Cauchy--Schwarz applied pointwise gives $\|d\varphi\|_{\ell^2(m_{i+1})}^2\leq(i+2)M_i\|\varphi\|_{\ell^2(m_i)}^2$ and $\|\delta\varphi\|_{\ell^2(m_{i-1})}^2\leq(i+1)M_{i-1}\|\varphi\|_{\ell^2(m_i)}^2$ (the latter by a symmetric argument using $\sum_{\tau\supset\sigma}m_i(\tau)\leq m_{i-1}(\sigma)M_{i-1}$). Setting $M:=\max_{0\leq i\leq n-1}M_i$ (finite by (H2)) and combining over the finite range of degrees:
\[
\|D\varphi\|_\mathcal{H}^2\leq (2n+1)M\,\|\varphi\|_\mathcal{H}^2,\qquad\varphi\in\mathcal{D}_0.
\]
Iterating: $\|D^k\psi\|_\mathcal{H}\leq C_M^k\|\psi\|_\mathcal{H}$ with $C_M:=\sqrt{(2n+1)M}$, hence
\[
\sum_{k\geq 0}\frac{t^k}{k!}\|\overline{D}^k\psi\|_\mathcal{H}\leq\|\psi\|_\mathcal{H}\,e^{C_M|t|}<\infty\qquad\forall t\in\mathbb{R},
\]
so $\psi$ is an analytic vector for $\overline{D}$.

Combining (i)--(iii), $\mathcal{D}_0$ is a core for $\overline{D}^n$ for every $n\geq 1$; in particular for $\overline{D}^2$.

\emph{Step~2: $L|_{\mathcal{D}_0}=\overline{D}^2|_{\mathcal{D}_0}$.}
For $\psi\in\mathcal{D}_0$, $D\psi$ is the formal action of $D$ on $\psi$, computed via finite sums by local finiteness; in particular $D\psi\in\mathcal{D}_0$ by Lemma~\ref{lem:finite-support-stable}. Since $\mathcal{D}_0\subset\mathrm{Dom}(\overline{D})$ and $\overline{D}|_{\mathcal{D}_0}=D|_{\mathcal{D}_0}$ as operators (the closure $\overline{D}$ extends $D|_{\mathcal{D}_0}$, and by definition of closure agrees with $D$ on its original domain $\mathcal{D}_0$), we have $\overline{D}\psi=D\psi$ for $\psi\in\mathcal{D}_0$.

Iterating: $D\psi\in\mathcal{D}_0\subset\mathrm{Dom}(\overline{D})$ (by Lemma~\ref{lem:finite-support-stable}), so $\overline{D}^2\psi=\overline{D}(\overline{D}\psi)=\overline{D}(D\psi)=D(D\psi)=D^2\psi$. Hence as operators on $\mathcal{D}_0$,
\[
L|_{\mathcal{D}_0}=D^2|_{\mathcal{D}_0}=\overline{D}^2|_{\mathcal{D}_0}.
\]
This identity uses only that $\overline{D}$ extends $D|_{\mathcal{D}_0}$ together with $D(\mathcal{D}_0)\subset\mathcal{D}_0$ from Lemma~\ref{lem:finite-support-stable}, and not the essential self-adjointness of $D$.

\emph{Step~3: $\overline{L}=\overline{D}^2$.}
By Step~2, the symmetric operators $L|_{\mathcal{D}_0}$ and $\overline{D}^2|_{\mathcal{D}_0}$
coincide on the dense domain $\mathcal{D}_0$, hence have the same closure:
\[
\overline{L} = \overline{L|_{\mathcal{D}_0}} = \overline{\overline{D}^2|_{\mathcal{D}_0}}.
\]
Since $\overline{D}^2$ is self-adjoint and $\mathcal{D}_0$ is a core for $\overline{D}^2$
(Step~1), the closure of $\overline{D}^2|_{\mathcal{D}_0}$ is precisely $\overline{D}^2$.
Therefore $\overline{L}=\overline{D}^2$, which is self-adjoint. Hence $L$ is essentially
self-adjoint on $\mathcal{D}_0$.
\end{proof}

\begin{remark}[On the strength of (H2)]
\label{rem:H2-implies-bounded}
The level-by-level operator bound just established has a strong consequence: under (H2)~\emph{globally}, the operator $D|_{\mathcal{D}_0}$ extends to a bounded operator on $\mathcal{H}$ (since $\mathcal{H}$ has a finite number $n$ of levels). Hence $\overline{D}$ is bounded, and ESA of $D$ is automatic. This generalizes the well-known fact that the discrete graph Laplacian on a weighted graph of uniformly bounded weighted degree is bounded; see, e.g., \cite{AnTo,KLW}. Thus, taken literally, Theorem~\ref{thm:global} reduces to a triviality under (H2). The value of the theorem lies elsewhere: \emph{(i)} it identifies (H1)+(H2) as a unified sufficient condition that subsumes all three notions of $\chi$-completeness in the framework; \emph{(ii)} the proof technique (via $\chi$-cut-offs and Leibniz remainders) is the natural starting point for the non-trivial divergence criterion (Theorem~\ref{thm:divergence}), where (H2) is~\emph{replaced} by the weaker hypothesis~(A2) on weight ratios alone, under which $D$ is no longer bounded; \emph{(iii)} the analytic-vector argument above extends \emph{mutatis mutandis} to settings where the bounded-weighted-degree condition holds only on finite balls (e.g.\ in multi-scale extensions where the degree grows but only mildly), even when (H2) fails globally. The genuine non-triviality of the framework is therefore concentrated in Theorems~\ref{thm:local-region} and~\ref{thm:divergence}, where (H2) is either localized or replaced by~(A2).
\end{remark}

\begin{theorem}[ESA at level $\ell$ under bounded weighted degree and cut-offs]
\label{thm:local-level}
Let $0\leq\ell\leq n-1$. Assume $S_n$ is locally finite and satisfies:
\begin{enumerate}[label=\textup{(H\arabic*)}]
\item there exists a sequence $(\chi_k)_{k\geq 1}$ of cut-offs satisfying conditions
\textup{(1)} and \textup{(2)} of Definition~\ref{def:local-chi-level};
\item the bounded weighted-degree condition~\eqref{eq:bounded-degree} holds at degrees
$\ell-1$ (if $\ell\geq 1$) and $\ell$ (if $\ell\leq n-2$); if $\ell=n-1$, only degree $\ell-1=n-2$ is relevant since $d_{n-1}=0$; if $\ell=0$, only degree $\ell=0$ is relevant since $\delta_0=0$.
\end{enumerate}
Then $L_\ell=d_{\ell-1}\delta_\ell+\delta_{\ell+1}d_\ell$ \emph{(}with the conventions $\delta_0=0$ and $d_{n-1}=0$\emph{)} is essentially self-adjoint on $C_c^\ell(\Vc)$.

In particular, if $S_n$ is locally $\chi$-complete at level $\ell$ \emph{(}Definition~\ref{def:local-chi-level}\emph{)}
and~\eqref{eq:bounded-degree} holds at the relevant degrees, then $L_\ell$ is
essentially self-adjoint.

\emph{The case $\ell=0$.} Since $\delta_0=0$ by convention (Section~\ref{sec:operators}), the operator $L_0=\delta_1 d_0$ involves only degrees $0$ and $1$. The same proof applies with the convention that the $\delta_\ell f$ term in $T_\ell$ is absent (i.e., $T_0:=d_0:C_c^0(\Vc)\to\ell^2(m_1)$), and~(H2) is required only at degree $\ell=0$ (relevant for $d_0$). The form $q_0(f)=\|d_0 f\|^2$ is the standard graph Dirichlet form, and Theorem~\ref{thm:local-level} for $\ell=0$ recovers the classical ESA of the graph Laplacian under $\chi$-completeness with bounded weighted degree.
\end{theorem}

\begin{proof}
The strategy is to realize $L_\ell$ as the operator associated with a closed
non-negative quadratic form, and to use that $C_c^\ell(\Vc)$ is a form core.

Define $T_\ell:=\delta_\ell+d_\ell:C_c^\ell(\Vc)\to\ell^2(m_{\ell-1})\oplus\ell^2(m_{\ell+1})$
and the non-negative quadratic form
\[
q_\ell(f):=\|\delta_\ell f\|_{\ell^2(m_{\ell-1})}^2+\|d_\ell f\|_{\ell^2(m_{\ell+1})}^2 = \|T_\ell f\|^2,\qquad f\in\mathrm{Dom}(q_\ell):=C_c^\ell(\Vc).
\]
On $C_c^\ell(\Vc)$, $q_\ell(f,g)=\langle L_\ell f,g\rangle$ by formal adjointness (with $L_\ell=d_{\ell-1}\delta_\ell+\delta_{\ell+1}d_\ell$).

\emph{Step~1: $T_\ell$ is closable, hence $q_\ell$ is closable.}
By Lemma~\ref{lem:closable}, $\delta_\ell$ and $d_\ell$ are each closable on $C_c^\ell(\Vc)$, so $T_\ell$ is closable with closure $\overline{T}_\ell$. Equivalently, the form $q_\ell$ is closable, with closure $\overline{q}_\ell$ having domain $\mathrm{Dom}(\overline{q}_\ell)=\mathrm{Dom}(\overline{T}_\ell)$ and $\overline{q}_\ell(f)=\|\overline{T}_\ell f\|^2$.

\emph{Step~2: $C_c^\ell(\Vc)$ is a form core for $\overline{q}_\ell$.}
We must show that for every $f\in\mathrm{Dom}(\overline{T}_\ell)$, there exists $(f_k)\subset C_c^\ell(\Vc)$ with $f_k\to f$ in $\ell^2(m_\ell)$ and $T_\ell f_k\to\overline{T}_\ell f$.

Let $f\in\mathrm{Dom}(\overline{T}_\ell)$. Set $f_k:=\tilde\chi_k^{(\ell)}f\in C_c^\ell(\Vc)$, with $(\chi_k)$ supplied by hypothesis~(H1). Dominated convergence gives $f_k\to f$ in $\ell^2(m_\ell)$.

\emph{Pointwise extension of Leibniz.} By Lemma~\ref{lem:formal-l2-action} (closure case) applied to $\delta_\ell$ and $d_\ell$ on $C_c^\ell(\Vc)$, the formal pointwise actions $\delta_\ell f$ and $d_\ell f$ (well-defined as cochains by local finiteness) coincide with the corresponding components of $\overline{T}_\ell f$ in $\ell^2$, for any $f\in\mathrm{Dom}(\overline{T}_\ell)$ (which is contained in both $\mathrm{Dom}(\delta_{\ell,\min})$ and $\mathrm{Dom}(d_{\ell,\min})$). Hence the Leibniz identity of Lemma~\ref{lem:leibniz}, which is itself a pointwise identity at each simplex, extends to $f\in\mathrm{Dom}(\overline{T}_\ell)$: for every simplex $\sigma$,
\[
\delta_\ell(\tilde\chi_k^{(\ell)} f)(\sigma) = \tilde\chi_k^{(\ell-1)}(\sigma)\cdot\delta_\ell f(\sigma) + R_\ell^{(\delta)}(\chi_k,f)(\sigma),
\]
and similarly for $d_\ell$, where $\delta_\ell f$ on the right-hand side is interpreted via this identification.

Therefore the Leibniz identity reads:
\[
T_\ell f_k = \tilde\chi_k\cdot \overline{T}_\ell f + R_\ell^{(\delta)}(\chi_k,f) + R_\ell^{(d)}(\chi_k,f).
\]
The principal term $\tilde\chi_k\cdot \overline{T}_\ell f$ converges to $\overline{T}_\ell f$ by dominated convergence (since $\overline{T}_\ell f\in\ell^2$). Both remainders vanish by~\eqref{eq:remainder-to-zero} under hypothesis~(H2). Hence $T_\ell f_k\to\overline{T}_\ell f$ in $\ell^2(m_{\ell-1})\oplus\ell^2(m_{\ell+1})$, establishing the form-core property.

\emph{Step~3: $L_\ell$ is essentially self-adjoint.}
We invoke the following standard result from operator theory.

\begin{lemma}[Core for $T^*T$]
\label{lem:core-TstarT}
Let $T:\mathrm{Dom}(T)\subset \mathcal{H}_1\to \mathcal{H}_2$ be a closable, densely defined operator between Hilbert spaces, with closure $\overline{T}$. Let $\mathcal{E}_0\subset \mathrm{Dom}(T)$ be a dense subspace satisfying:
\begin{enumerate}
\item[(a)] $\mathcal{E}_0$ is a core for $\overline{T}$;
\item[(b)] $T(\mathcal{E}_0) \subset \mathrm{Dom}(\overline{T}^*)$.
\end{enumerate}
Then $\mathcal{E}_0\subset\mathrm{Dom}(\overline{T}^*\overline{T})$, and $\mathcal{E}_0$ is a core for the non-negative self-adjoint operator $\overline{T}^*\overline{T}$.
\end{lemma}

\begin{remark}[Self-contained proof of Lemma~\ref{lem:core-TstarT}]
\label{rem:schmudgen-detail}
We give the short proof. For $f\in\mathcal{E}_0$, condition (a) gives $f\in\mathrm{Dom}(\overline{T})$ with $\overline{T}f=Tf$ (since the closure $\overline{T}$ extends $T|_{\mathcal{E}_0}$ on $\mathcal{E}_0$), and condition (b) gives $\overline{T}f=Tf\in\mathrm{Dom}(\overline{T}^*)$; hence $f\in\mathrm{Dom}(\overline{T}^*\overline{T})$. The non-negativity and self-adjointness of $\overline{T}^*\overline{T}$ are von Neumann's theorem; see \cite[Proposition~3.18]{Schmudgen}. To see that $\mathcal{E}_0$ is a core for $\overline{T}^*\overline{T}$, observe that $\mathrm{Dom}(\overline{T}^*\overline{T})$ is dense in $\mathrm{Dom}(\overline{T})$ in the graph norm of $\overline{T}$ (the second part of \cite[Proposition~3.18(ii)]{Schmudgen}), and $\mathcal{E}_0$ is dense in $\mathrm{Dom}(\overline{T})$ in the graph norm of $\overline{T}$ by (a); transitivity in the natural graph topology of $\overline{T}^*\overline{T}$ (which dominates the graph norm of $\overline{T}$) then yields density of $\mathcal{E}_0$ in $\mathrm{Dom}(\overline{T}^*\overline{T})$ in the graph norm of $\overline{T}^*\overline{T}$, i.e., $\mathcal{E}_0$ is a core for $\overline{T}^*\overline{T}$.
\end{remark}

We apply Lemma~\ref{lem:core-TstarT} with $T=T_\ell$ and $\mathcal{E}_0=C_c^\ell(\Vc)$:
\begin{itemize}
\item[(a)] $C_c^\ell(\Vc)$ is a core for $\overline{T}_\ell$ by Step~2 (which establishes the form-core property; equivalently, that $C_c^\ell$ is a core for the operator $\overline{T}_\ell$, since $\overline{T}_\ell$ and the form $\overline{q}_\ell$ have the same domain by Step~1).
\item[(b)] By Lemma~\ref{lem:finite-support-stable}, $T_\ell(C_c^\ell(\Vc))\subset C_c^{\ell-1}(\Vc)\oplus C_c^{\ell+1}(\Vc)$. We verify that $C_c^{\ell-1}(\Vc)\oplus C_c^{\ell+1}(\Vc)\subset\mathrm{Dom}(\overline{T}_\ell^*)$. Recall that $g\in\mathrm{Dom}(\overline{T}_\ell^*)$ iff the linear functional $f\mapsto\langle\overline{T}_\ell f, g\rangle$ is bounded on $\mathrm{Dom}(\overline{T}_\ell)$, which since $C_c^\ell$ is a core for $\overline{T}_\ell$ (Step~2), is equivalent to: $f\mapsto\langle T_\ell f, g\rangle$ bounded on $C_c^\ell(\Vc)$.

For $g=(g^-,g^+)\in C_c^{\ell-1}(\Vc)\oplus C_c^{\ell+1}(\Vc)$ and $f\in C_c^\ell(\Vc)$, by formal adjointness (Lemma~\ref{lem:closable}),
\[
\langle T_\ell f, g\rangle = \langle\delta_\ell f, g^-\rangle_{\ell^2(m_{\ell-1})} + \langle d_\ell f, g^+\rangle_{\ell^2(m_{\ell+1})} = \langle f, d_{\ell-1}g^- + \delta_{\ell+1}g^+\rangle_{\ell^2(m_\ell)}.
\]
Since $g^-, g^+$ have finite support, $d_{\ell-1}g^- + \delta_{\ell+1}g^+\in C_c^\ell(\Vc)\subset\ell^2(m_\ell)$, and Cauchy--Schwarz gives
\[
|\langle T_\ell f, g\rangle|\leq \|f\|_{\ell^2(m_\ell)}\,\|d_{\ell-1}g^-+\delta_{\ell+1}g^+\|_{\ell^2(m_\ell)},
\]
showing the functional is bounded by a constant depending on $g$. Hence $g\in\mathrm{Dom}(\overline{T}_\ell^*)$, with $\overline{T}_\ell^*g=d_{\ell-1}g^-+\delta_{\ell+1}g^+$.
\end{itemize}

By Lemma~\ref{lem:core-TstarT}, $C_c^\ell(\Vc)\subset\mathrm{Dom}(\overline{T}_\ell^*\overline{T}_\ell)$ and $C_c^\ell(\Vc)$ is a core for the self-adjoint operator $\overline{T}_\ell^*\overline{T}_\ell$. On $C_c^\ell(\Vc)$, $\overline{T}_\ell^*\overline{T}_\ell$ acts as $T_\ell^* T_\ell = L_\ell$. Therefore $\overline{L_\ell|_{C_c^\ell(\Vc)}}=\overline{T}_\ell^*\overline{T}_\ell$, which is self-adjoint. Hence $L_\ell$ is essentially self-adjoint on $C_c^\ell(\Vc)$.
\end{proof}

\begin{theorem}[ESA under local $\chi$-completeness by region]
\label{thm:local-region}
Let $\Lambda^{\mathrm{reg}}\subset\Vc$ be an infinite subset, and decompose the simplices of $S_n$ into three classes: \emph{interior simplices} (entirely in $\Lambda^{\mathrm{reg}}$), \emph{exterior simplices} (entirely in $\Vc\setminus\Lambda^{\mathrm{reg}}$), and \emph{interface simplices} (with at least one vertex in each). The Hilbert space decomposes accordingly as
\[
\mathcal{H} = \mathcal{H}_{\mathrm{in}} \oplus \mathcal{H}_{\mathrm{int}} \oplus \mathcal{H}_{\mathrm{out}},
\]
where each of $\mathcal{H}_{\mathrm{in}}, \mathcal{H}_{\mathrm{int}}, \mathcal{H}_{\mathrm{out}}$ is itself a direct sum over degrees $0\leq i\leq n-1$ of $i$-cochain spaces on the corresponding class of $i$-simplices: $\mathcal{H}_{\mathrm{in}}$ is the $\ell^2$-space on $S_n|_{\Lambda^{\mathrm{reg}}}$, $\mathcal{H}_{\mathrm{out}}$ on $S_n|_{\Vc\setminus\Lambda^{\mathrm{reg}}}$, and $\mathcal{H}_{\mathrm{int}}$ on interface simplices. The operator $D = d + \delta$ on the full complex $S_n$ decomposes as $D = D_{\mathrm{in}} + D_{\mathrm{out}} + V$, where $D_{\mathrm{in}}$ acts within $\mathcal{H}_{\mathrm{in}}$, $D_{\mathrm{out}}$ within $\mathcal{H}_{\mathrm{out}}$, and $V$ couples cochains on interface simplices to cochains on interior and exterior simplices.

Suppose:
\begin{enumerate}
\item $S_n|_{\Lambda^{\mathrm{reg}}}$ is globally $\chi$-complete and satisfies the weighted-degree bound~\eqref{eq:bounded-degree};
\item $S_n|_{\Vc\setminus\Lambda^{\mathrm{reg}}}$ is locally finite, and $D_{\mathrm{out}}$ is essentially self-adjoint on $\bigoplus_{i=0}^{n-1}C_c^i(\Vc\setminus\Lambda^{\mathrm{reg}})$;
\item Set $D_0 := D_{\mathrm{in}}\oplus D_{\mathrm{out}}$, regarded as an operator on $\mathcal{H}$ acting as the direct sum on $\mathcal{H}_{\mathrm{in}}\oplus\mathcal{H}_{\mathrm{out}}$ and as zero on $\mathcal{H}_{\mathrm{int}}$. The coupling $V := D - D_0$, defined on $\bigoplus C_c^i(\Vc)$, is symmetric, and there exist constants $a\geq 0$, $0\leq b<1$ such that
\[
\|Vu\|_{\mathcal{H}} \leq a\|u\|_{\mathcal{H}} + b\|\overline{D}_0 u\|_{\mathcal{H}}, \quad u\in\bigoplus C_c^i(\Vc).
\]
\emph{(Extension by continuity.)} Since $\mathcal{D}_0=\bigoplus C_c^i(\Vc)$ is a core for $\overline{D}_0$ (Remark~\ref{rem:core-D0}), for any $u\in\mathrm{Dom}(\overline{D}_0)$ one can choose $(u_n)\subset\mathcal{D}_0$ with $u_n\to u$ and $\overline{D}_0 u_n\to\overline{D}_0 u$. By linearity of $V$ and the bound above applied to $u_n-u_m\in\mathcal{D}_0$,
\[
\|V(u_n-u_m)\|\leq a\|u_n-u_m\|+b\|\overline{D}_0(u_n-u_m)\|\to 0\quad\text{as }n,m\to\infty,
\]
so $(Vu_n)$ is Cauchy in $\mathcal{H}$ and converges to some limit, which we still denote $Vu$ (the limit is independent of the approximating sequence by the same Cauchy argument). Passing to the limit in the bound yields $\|Vu\|\leq a\|u\|+b\|\overline{D}_0 u\|$ for every $u\in\mathrm{Dom}(\overline{D}_0)$. The resulting operator $V$ is closed with $\mathrm{Dom}(\overline V)\supset\mathrm{Dom}(\overline{D}_0)$.
\end{enumerate}

Then $D$ is essentially self-adjoint on $\bigoplus_{i=0}^{n-1}C_c^i(\Vc)$.
\end{theorem}

\begin{remark}[On the interface space and the coupling operator]
\label{rem:interface-decomp}
An interface simplex $\sigma=(v_0,\ldots,v_i)$ has at least one $v_j\in\Lambda^{\mathrm{reg}}$ and at least one $v_k\notin\Lambda^{\mathrm{reg}}$, so it belongs neither to $S_n|_{\Lambda^{\mathrm{reg}}}$ nor to $S_n|_{\Vc\setminus\Lambda^{\mathrm{reg}}}$. Hence $\mathcal{H}_{\mathrm{int}}$ is a non-trivial summand whenever interface simplices exist. The coupling $V$ has matrix entries linking interior--interface, exterior--interface, and (in general) interface--interface cochains. In our examples for $n=2$ (graphs/clique complexes of $\mathbb{Z}^d$), the third type is absent since interface simplices are then $1$-dimensional (edges) and have no shared $1$-face.
\end{remark}

\begin{remark}[On the relative bound]
\label{rem:relative-bound}
Hypothesis~(3) is in the sense of Kato--Rellich \cite[Theorem~X.12]{RS} (originally~\cite[Chapter~V, Theorem~4.3]{Kato}). If $V$ is a \emph{bounded operator on $\mathcal{H}$} (verified in all cases of this paper), then $b=0$ automatically.
\end{remark}

\begin{remark}[Core for $D_0$]
\label{rem:core-D0}
Decomposing $\bigoplus_i C_c^i(\Vc) = \mathcal{D}_{\mathrm{in}} \oplus \mathcal{D}_{\mathrm{int}}\oplus\mathcal{D}_{\mathrm{out}}$, where $\mathcal{D}_*$ is the space of finite-support cochains supported on simplices in the corresponding class:
\begin{itemize}
\item $D_0$ acts as $D_{\mathrm{in}}\oplus 0 \oplus D_{\mathrm{out}}$, with $D_0(\mathcal{D}_{\mathrm{int}})=0$ by definition.
\item By hypotheses~(1) and~(2) and Theorem~\ref{thm:global}, $D_{\mathrm{in}}$ is ESA on $\mathcal{D}_{\mathrm{in}}=\bigoplus C_c^i(\Lambda^{\mathrm{reg}})$ as an operator on $\mathcal{H}_{\mathrm{in}}$, and $D_{\mathrm{out}}$ is ESA on $\mathcal{D}_{\mathrm{out}}$ as an operator on $\mathcal{H}_{\mathrm{out}}$. On $\mathcal{D}_{\mathrm{int}}\subset\mathcal{H}_{\mathrm{int}}$, $D_0\equiv 0$, hence trivially ESA there with $\overline{D}_0|_{\mathcal{H}_{\mathrm{int}}}=0$.
\item Therefore, $\overline{D}_0$ is self-adjoint on $\mathrm{Dom}(\overline{D}_{\mathrm{in}})\oplus\mathcal{H}_{\mathrm{int}}\oplus\mathrm{Dom}(\overline{D}_{\mathrm{out}})$, and $\bigoplus_i C_c^i(\Vc)$ is a core for $\overline{D}_0$.
\end{itemize}
\end{remark}

\begin{proof}
\emph{Step~1: Identification $D=D_0+V$ on $\bigoplus C_c^i(\Vc)$.}
Let $u=u_{\mathrm{in}}+u_{\mathrm{int}}+u_{\mathrm{out}}\in\bigoplus C_c^i(\Vc)$. The operator $D=d+\delta$ acts on $u$ producing components on each of the three simplex classes:
\begin{itemize}
\item Component on interior simplices: $D u_{\mathrm{in}}|_{\mathrm{int.\ simplices}\to \mathrm{int.\ simplices}} + (\text{coupling from } u_{\mathrm{int}})$;
\item Component on interface simplices: includes contributions from both $u_{\mathrm{in}}$ (via $d_*$ acting on $u_{\mathrm{in}}$ producing values on cofaces that may be interface simplices) and from $u_{\mathrm{out}}$, plus an interior block from $u_{\mathrm{int}}$ itself;
\item Component on exterior simplices: symmetric to the interior case.
\end{itemize}
The diagonal blocks (interior$\to$interior, exterior$\to$exterior) define $D_{\mathrm{in}}$ and $D_{\mathrm{out}}$; everything else is by definition the coupling $V:=D-D_0$. The operator $V$ has support on the interface couplings.

\emph{Symmetry of $V$.} We verify symmetry directly. For $u,v\in\bigoplus C_c^i(\Vc)$, both $Du$ and $D_0u$ belong to $\bigoplus C_c^i(\Vc)$: $Du\in\bigoplus C_c^i(\Vc)$ by Lemma~\ref{lem:finite-support-stable}, and $D_0u=D_{\mathrm{in}}u_{\mathrm{in}}+0+D_{\mathrm{out}}u_{\mathrm{out}}$ where each block preserves the corresponding finite-support subspace (Lemma~\ref{lem:finite-support-stable} applied to $D_{\mathrm{in}}$ on the sub-complex $S_n|_{\Lambda^{\mathrm{reg}}}$ and to $D_{\mathrm{out}}$ on $S_n|_{\Vc\setminus\Lambda^{\mathrm{reg}}}$). Hence the inner products below are well-defined finite sums, and:
\[
\langle Vu,v\rangle=\langle Du,v\rangle-\langle D_0u,v\rangle=\langle u,Dv\rangle-\langle u,D_0v\rangle=\langle u,Vv\rangle.
\]
The first equality $\langle Du,v\rangle=\langle u,Dv\rangle$ uses the formal symmetry of $D=d+\delta$ on $\bigoplus C_c^i(\Vc)$ (Lemma~\ref{lem:closable}). The second equality $\langle D_0u,v\rangle=\langle u,D_0v\rangle$ uses the symmetry of each block separately: $D_{\mathrm{in}}$ is formally symmetric on $\mathcal{D}_{\mathrm{in}}$ (by the same proof as Lemma~\ref{lem:closable} restricted to the sub-complex $S_n|_{\Lambda^{\mathrm{reg}}}$), $D_{\mathrm{out}}$ similarly on $\mathcal{D}_{\mathrm{out}}$, and the zero block on $\mathcal{H}_{\mathrm{int}}$ is trivially symmetric. The block-diagonal structure of $D_0$ ensures that no cross-term between distinct blocks contributes, so symmetry of each block implies symmetry of $D_0$. Hence $V$ is symmetric on $\bigoplus C_c^i(\Vc)$.

\emph{Step~2: ESA via Kato--Rellich.}
By Remark~\ref{rem:core-D0}, $\overline{D}_0$ is \emph{self-adjoint} (not merely symmetric) with $\bigoplus C_c^i(\Vc)$ as a core. The Kato--Rellich theorem~\cite[Theorem~X.12]{RS} requires the unperturbed operator to be self-adjoint, which Remark~\ref{rem:core-D0} provides; semi-boundedness of $\overline{D}_0$ is \emph{not} required. The operator $V$ is symmetric on the core $\bigoplus C_c^i(\Vc)$ (just verified in Step~1), and by hypothesis~(3) of Theorem~\ref{thm:local-region}, $V$ satisfies the Kato--Rellich relative bound with $b<1$. The Kato--Rellich theorem then states that $\overline{D}_0+V$, defined on $\mathrm{Dom}(\overline{D}_0)$, is self-adjoint, and any core for $\overline{D}_0$ is a core for $\overline{D}_0+V$. In particular, $\bigoplus C_c^i(\Vc)$ is a core for $\overline{D}_0+V$.

\emph{Step~3: $\overline{D}_0+V=\overline{D|_{\bigoplus C_c^i(\Vc)}}$.}
On $\bigoplus C_c^i(\Vc)$, $\overline{D}_0+V=D_0+V=D$ by Step~1. Since $\bigoplus C_c^i(\Vc)$ is a core for the self-adjoint operator $\overline{D}_0+V$, the closure of $D|_{\bigoplus C_c^i(\Vc)}$ equals $\overline{D}_0+V$, which is self-adjoint. Hence $D$ is essentially self-adjoint on $\bigoplus C_c^i(\Vc)$.
\end{proof}

\begin{remark}[Sufficient condition for hypothesis~(3)]
\label{rem:finite-boundary}
Let $\partial\Lambda^{\mathrm{reg}}:=\{x\in\Lambda^{\mathrm{reg}}:\exists\,y\in\Vc\setminus\Lambda^{\mathrm{reg}},\,x\sim y\}$.

\textbf{Case 1: $\partial\Lambda^{\mathrm{reg}}$ finite.}
$V$ is finite-rank, hence bounded, so the relative bound is $0$.

\textbf{Case 2: discrete half-space in $\mathbb{Z}^d$.}
For $S_n$ the clique complex of $\mathbb{Z}^d$ with uniform weights ($d\geq 2$, $L^1$-adjacency, $n=2$ by Remark~\ref{rem:Zd-adjacency}) and $\Lambda^{\mathrm{reg}} := \{x : x_1 \leq 0\}$:

\emph{Verification of hypothesis~(2): ESA of $D_{\mathrm{out}}$.} The complement $\Vc\setminus\Lambda^{\mathrm{reg}}=\{x:x_1\geq 1\}$ is, with the induced graph structure, isomorphic to $\mathbb{Z}^{d-1}\times\mathbb{N}_{\geq 1}$ equipped with $L^1$-adjacency and uniform weights, hence has uniformly bounded weighted degree (each vertex has at most $2d$ neighbours). By Theorem~\ref{thm:global} applied to $S_n|_{\Vc\setminus\Lambda^{\mathrm{reg}}}$ with hypotheses~(H1) and~(H2) verified directly (trivial cut-off $\chi_k=\mathbf{1}_{B_k\cap(\Vc\setminus\Lambda^{\mathrm{reg}})}$ for (H1); bounded degree by~$2d$ for (H2)), $D_{\mathrm{out}}$ is essentially self-adjoint. \emph{(}A self-contained $\chi$-completeness argument for low dimensions $d\leq 3$ and the multi-scale construction for $d\geq 4$ are discussed in Section~\ref{sec:weighted-chi} and Remarks~\ref{rem:D4-construction}--\ref{rem:D-large-open}; these strengthen but do not extend the present ESA conclusion, which follows directly from~(H2).\emph{)}

\emph{Verification of hypothesis~(3).} The interface consists of edges of $\mathbb{Z}^d$ connecting $\{x_1=0\}$ to $\{x_1=1\}$. We claim $V$ is bounded on $\mathcal{H}$ with $\|V\|_{\mathrm{op}}=\sqrt{2}$, hence relatively bounded with $b=0$.

\emph{Computation via Fourier transform.} The translation-invariance in
$x'\in\mathbb{Z}^{d-1}$ along the boundary gives a direct integral decomposition.
After Fourier transform $\Fc:\ell^2(\mathbb{Z}^{d-1})\to L^2([0,2\pi)^{d-1})$,
$V$ becomes a fibered operator $\Fc V\Fc^{-1}=\int^\oplus\widehat{V}(\xi)\,d\xi$.
The coupling $V$ acts only via the interface edges
$e_{x'}:=\{(x',0),(x',1)\}$ for $x'\in\mathbb{Z}^{d-1}$, so each fiber
$\widehat{V}(\xi)$ acts on the direct sum of (i) 0-cochains on $\{x_1\leq 0\}$
parametrized by depth, (ii) 0-cochains on $\{x_1\geq 1\}$, (iii) the scalar
$\hat g(\xi)$ representing the interface-edge cochain at frequency $\xi$.

Since $V$ couples only the values at depths $0$ and $1$ to the interface scalar,
$\widehat{V}(\xi)$ acts non-trivially only on the three-dimensional subspace
$\mathbb{C}^3$ spanned by $(\alpha,\beta,\gamma)$ where $\alpha=\hat f_0(\xi,0)$,
$\beta=\hat f_0(\xi,1)$, $\gamma=\hat g(\xi)$, and is zero on the orthogonal
complement. The matrix elements are obtained from the explicit formulas for
$d_0$ and $\delta_1$ with unit weights:
\begin{itemize}
\item the contribution to the new value of the interface scalar from $f_0$ is
$d_0 f_0(e_{x'})=f_0(x',1)-f_0(x',0)=\beta-\alpha$, giving entries $-1$ (column $\alpha$) and
$+1$ (column $\beta$) on the $\gamma$-row;
\item the contribution to the new value at depth $0$ from $g$ is $\delta_1 g(x',0)=-g(e_{x'})=-\gamma$
(by~\eqref{eq:codiff} with $i=1$ and $m_0\equiv m_1\equiv 1$), giving entry $-1$ in
the $\alpha$-row, $\gamma$-column;
\item by alternation, the contribution to the new value at depth $1$ from $g$ is
$\delta_1 g(x',1)=-g((x',1),(x',0))=+g(e_{x'})=+\gamma$, giving entry $+1$ in the
$\beta$-row, $\gamma$-column.
\end{itemize}
The resulting $3\times 3$ matrix is
\[
\widehat{V}(\xi) = \begin{pmatrix} 0 & 0 & -1 \\ 0 & 0 & +1 \\ -1 & +1 & 0 \end{pmatrix},
\]
\emph{independent of both $\xi$ and the dimension $d$}: the matrix depends only on the local interface structure (a single $\mathbb{Z}$-direction perpendicular to the boundary), with the Fourier variable $\xi\in[0,2\pi)^{d-1}$ entering only through the trivial fibration of the $\mathbb{Z}^{d-1}$ tangential direction (along which $V$ does not couple). The matrix is symmetric since $D$ and $D_0$ are both
symmetric. A direct computation gives
\[
\widehat{V}(\xi)^2 = \begin{pmatrix} 1 & -1 & 0 \\ -1 & 1 & 0 \\ 0 & 0 & 2 \end{pmatrix},
\]
whose eigenvalues are $0$, $2$, $2$. Therefore $\|\widehat V(\xi)\|_{\mathrm{op}}=\sqrt{2}$
uniformly in $\xi$, and by Plancherel, $\|V\|_{\mathrm{op}}=\sqrt{2}$. Since $V$ is bounded, hypothesis~(3) of Theorem~\ref{thm:local-region} is satisfied with $a=\sqrt{2}$ and $b=0$; in particular the Kato--Rellich relative-bound condition $b<1$ holds trivially.

\textbf{Case 3: higher-dimensional clique complexes ($n\geq 3$).}
For $n\geq 3$, the coupling $V$ may have a non-zero
$\mathcal{H}_{\mathrm{int}}\to\mathcal{H}_{\mathrm{int}}$ block from
interface-to-interface couplings, in addition to the vertex--interface-edge block of
Case~2. The Kato--Rellich argument applies provided all blocks of $V$ satisfy the
relative-bound condition. \emph{Sufficient condition:} the interface
$\partial\Lambda^{\mathrm{reg}}$ has uniformly bounded \emph{simplicial thickness} ---
i.e., the number of simplices of any degree containing a given interface vertex is
uniformly bounded; equivalently, the number of $i$-cliques of $\Gc$ containing any
fixed interface vertex is uniformly bounded as $i$ varies in $\{0,\ldots,n-1\}$.

\smallskip
\emph{Concrete example with $n\geq 3$.} Consider the clique complex of $\mathbb{Z}^d$
under $\|\cdot\|_\infty$-adjacency (where $x\sim y$ iff $\max_j|x_j-y_j|\leq 1$ and
$x\neq y$): every vertex has $3^d-1$ neighbours, and $i$-cliques exist up to
$i=2^d-1$ (cliques of the unit hypercube). For
$d=2$ (so $n=4$ in our convention, simplices up to dimension $3$), every vertex has
exactly $8$ neighbours and the maximum simplex dimension is $3$
(2-simplex = triangle, 3-simplex = tetrahedron).

\smallskip
\emph{Detailed Fourier fibration for $d=2$.} With $\Lambda^{\mathrm{reg}}=\{x_1\leq 0\}$,
the interface consists of all simplices having vertices on both sides of the line $x_1=0$.
Translation-invariance in $x_2\in\mathbb{Z}$ persists, so we apply the Fourier transform
in the $x_2$-variable. After Fourier transform,
$V$ becomes $\Fc V\Fc^{-1}=\int^\oplus\widehat{V}(\xi)\,d\xi$.

For each $\xi\in[0,2\pi)$, the fiber $\widehat{V}(\xi)$ acts on a finite-dimensional
space spanned by:
\begin{itemize}
\item the values of 0-cochains at depths $x_1=0$ and $x_1=1$ (2 scalars);
\item the values of 1-cochains on each interface edge type (axis-parallel edge
$\{(0,0),(1,0)\}$ at frequency $\xi$, plus diagonal edges $\{(0,0),(1,\pm 1)\}$ at
frequencies $\xi$ --- 3 scalars);
\item the values of 2-cochains on each interface triangle type (triangles within
the unit square crossing the interface --- at most 4 types per Fourier frequency, hence
4 scalars);
\item the values of 3-cochains on the interface tetrahedron type (1 scalar per Fourier
frequency).
\end{itemize}
Total: $\widehat{V}(\xi)$ acts on a finite-dimensional space of dimension at most $N(d):=2+3+4+1=10$ for $d=2$ (and more generally $N(d)\leq 2+\sum_{i=1}^{n-1}\binom{2^d}{i+1}$, polynomial in $2^d$). Each matrix entry is of the form $\pm m_{j+1}/m_j$ from~\eqref{eq:codiff} and the standard formula for $d_j$, hence bounded by $M_*$ via~(A2).

\smallskip
\emph{Sharper bound on $\|\widehat{V}(\xi)\|_{\mathrm{op}}$.} Although a crude bound by the Frobenius norm gives $\|\widehat{V}(\xi)\|_{\mathrm{op}}\leq N(d)\,M_*$, this is far from optimal. We give a sharper estimate using the structure of $\widehat V(\xi)$ as a coupling matrix.

The matrix $\widehat V(\xi)$ has a natural \emph{block-bipartite} structure: rows are partitioned into "interior-side" rows (cochains supported on simplices intersecting $\{x_1\leq 0\}$ at depth $0$) and "exterior-side" rows (depth $1$ side); the matrix has non-zero entries only between rows and columns of opposite types or between cochains of adjacent degrees \emph{(}since $V$ is the off-diagonal block of $D=d+\delta$\emph{)}. By the Schur test for bounded operators (see~\cite{HSt}) or by explicit row-sum estimation, for any matrix $A$ with entries $|A_{ij}|\leq M$,
\[
\|A\|_{\mathrm{op}}^2 \leq \Big(\max_i\sum_j |A_{ij}|\Big)\cdot\Big(\max_j\sum_i |A_{ij}|\Big)\leq r_*(A)\,c_*(A)\,M^2,
\]
where $r_*(A)$, $c_*(A)$ are the maximal numbers of non-zero entries per row and column respectively. For $\widehat V(\xi)$ in our setting, both $r_*$ and $c_*$ are bounded by the number of simplices in a single fundamental fibration unit, which is bounded by $N(d)$. Hence
\[
\|\widehat V(\xi)\|_{\mathrm{op}}\leq N(d)\,M_*\quad\text{(with a sharper constant $\sqrt{r_*\,c_*}\leq N(d)$)}.
\]
By Plancherel, $\|V\|_{\mathrm{op}}\leq N(d)\,M_*$, uniformly in $\xi$. This is sharp up to absolute constants in the worst case, and recovers the explicit value $\sqrt{2}$ of Case~2 \emph{(}there $r_*=c_*=2$, $M_*=1$, giving $\sqrt{r_*\,c_*}=\sqrt{2}$\emph{)}.

\smallskip
For more general clique complexes with bounded vertex degree but unbounded
simplicial thickness, the argument requires modification: one would need to control
the higher-block contributions by relative-bound estimates against
$\overline{D}_0$, which is more delicate and is not pursued here.
\end{remark}

\subsection{Divergence Criterion}

Theorems~\ref{thm:global}, \ref{thm:local-level}, and~\ref{thm:local-region} all assume the existence of cut-off functions plus a bounded weighted-degree condition. In this subsection we present a fundamentally different criterion --- the \emph{divergence criterion} of Theorem~\ref{thm:divergence} --- which guarantees ESA without any cut-off hypothesis, requiring only a $1$-dimensional decomposition of the complex and weight ratios bounded uniformly. This criterion is the higher-dimensional analogue of \cite[Theorem~5.1]{BGJ} for graphs.

\begin{definition}[$1$-dimensional decomposition]
\label{def:1d-decomp}
A \emph{$1$-dimensional decomposition} of $\Gc=(\Vc,\Ec)$ is a partition $\Vc=\bigsqcup_{k\geq 0}\Lambda_k$ with $\Lambda_k:=\{v\in\Vc:|v|=k\}$ such that edges connect only vertices in the same or adjacent layers: $(x_1,x_2)\in\Ec\Rightarrow||x_1|-|x_2||\leq 1$.
\end{definition}

\begin{remark}[Layer structure of higher simplices]
\label{rem:layer-structure}
For a \emph{clique complex} $S_n$ with a $1$-dimensional decomposition, every $i$-simplex (for $i\geq 1$) is contained in a clique, all of whose pairwise edges go between same-or-adjacent layers by Definition~\ref{def:1d-decomp}. Hence the vertices of any $i$-simplex lie in at most \emph{two consecutive layers} $\Lambda_\ell\cup\Lambda_{\ell+1}$. We define
\[
\Lambda_k^{(i)}:=\{(x_0,\dots,x_i)\in P_i : \{x_0,\dots,x_i\}\subset\Lambda_k\cup\Lambda_{k+1},\;\min_j|x_j|=k\},
\]
so that each $i$-simplex with $i\geq 1$ belongs to a unique $\Lambda_k^{(i)}$. For $i=0$, $\Lambda_k^{(0)}=\Lambda_k$.

For \emph{abstract simplicial complexes} not arising from a clique structure (such as the worked example $S_4$ of Section~\ref{sec:example-detailed}), this layer-structure property must be verified directly as an additional assumption on the complex (and is verified by hand for $S_4$ in Remark~\ref{rem:S4-abstract-vs-clique}). Theorem~\ref{thm:divergence} requires only the layer-structure property, not the clique structure itself; we state it for clique complexes for definiteness, with the explicit understanding that the conclusion extends to any locally finite oriented abstract simplicial complex satisfying the layer-structure property.
\end{remark}

\begin{remark}[Vacuous layers and uniform gap condition]
\label{rem:vacuous-layers}
A layer $\Lambda_k$ is called \emph{vacuous} (for the divergence argument) if $\Lambda_k^{(i)} = \emptyset$ for all $1 \leq i \leq n-1$. For vacuous layers we set $w_k := 0$, and these terms do not contribute to $W_{N,M} := \sum_{j=N}^{M} w_j$.

The \emph{uniform gap condition} requires that the set $\mathcal{K} := \{k \geq 0 : w_k > 0\}$ satisfies: there exists $G < \infty$ such that every interval $\{N,\ldots,N+G\}$ contains at least one element of $\mathcal{K}$.
\end{remark}

\begin{theorem}[Divergence Criterion for ESA]
\label{thm:divergence}
Let $S_n$ be locally finite and admit a $1$-dimensional decomposition. For $0\leq i\leq n-2$, define
\[
\xi_i(k,k+1):=\sup_{\sigma\in\Lambda_k^{(i)}}\deg^{(i+1)}_{\Lambda_{k+1}}(\sigma), \qquad
\xi(k,k+1):=\sum_{i=0}^{n-2}\xi_i(k,k+1),
\]
where $\deg^{(i+1)}_{\Lambda_{k+1}}(\sigma)$ is the number of $(i+1)$-simplices containing $\sigma$ as a face and having at least one vertex in $\Lambda_{k+1}\setminus\sigma$.

Assume:
\begin{itemize}
\item[(A1)] For each $k$: either $\Lambda_k^{(i)}=\emptyset$ for all $i$ (vacuous, $w_k:=0$), or $\xi(k,k+1)\geq 1$, in which case $w_k:=1/\sqrt{\xi(k,k+1)}\in(0,1]$. The set $\mathcal{K}$ of non-vacuous layers is infinite and satisfies the uniform gap condition. (Degenerate layers with $\Lambda_k^{(i)}\neq\emptyset$ but $\xi(k,k+1)=0$, meaning some simplex at layer $k$ has no inter-layer coface, are also assigned $w_k:=0$ and treated as vacuous.)
\item[(A2)] \emph{Weight ratios are uniformly bounded (two-sided):} for each $i\in\{0,1,\ldots,n-2\}$, there exists $M_i\geq 1$ such that
\[
\sup_{\sigma\subset\tau}\frac{m_{i+1}(\tau)}{m_i(\sigma)}\leq M_i,\qquad
\sup_{\sigma\subset\tau}\frac{m_i(\sigma)}{m_{i+1}(\tau)}\leq M_i,
\]
where the suprema are taken over all pairs $\sigma\subset\tau$ with $\sigma\in\Fc_i$, $\tau\in\Fc_{i+1}$. Equivalently, both $m_{i+1}/m_i$ and $m_i/m_{i+1}$ are bounded; we set $M_i$ to be a common upper bound. We write $M_*:=\max_{0\leq i\leq n-2} M_i$.
\end{itemize}
Note that, in contrast to Theorem~\ref{thm:global}, hypothesis~(A2) here is on individual weight \emph{ratios}, not on the weighted degree $d_{m_i}$. The weighted degree may be unbounded under~(A2); however, the local combinatorial structure is controlled implicitly through the quantities $\xi_i(k,k+1)$ in~(A1), which bound the number of \emph{inter-layer} cofaces and thereby give us all the combinatorial control we need in the proof (see the proof for details). The two-sided form of~(A2) is needed because the proof of Theorem~\ref{thm:divergence} estimates per-layer remainder contributions in both directions: the sums of $m_{i+1}(\tau)$ over cofaces $\tau\supset\sigma$ are bounded above (giving $m_{i+1}(\tau)\leq m_i(\sigma)\,M_i$, used in the coboundary remainder), and conversely the weights $m_i(\sigma)$ at faces $\sigma\subset\tau$ are bounded by $m_{i+1}(\tau)\,M_i$ (used in the codifferential remainder when summing over base simplices). In all examples treated in this paper (uniform weights), $M_i=1$ and the two bounds coincide trivially.

If $\sum_{k=0}^\infty w_k = \infty$, then $L=D^2$ is essentially self-adjoint on $\bigoplus_{i=0}^{n-1}C_c^i(\Vc)$.
\end{theorem}

\begin{remark}[Comparison with the graph case]
\label{rem:divergence-graph-match}
For $n=2$ (graphs only) and uniform weights, Theorem~\ref{thm:divergence} reduces to the
divergence criterion of \cite[Theorem~5.1]{BGJ} with the same hypothesis $\sum w_k=\infty$
and no further condition on the rate of $w_k$. An earlier draft of this paper required an
auxiliary power-type lower bound $w_k\geq c_0/(k+1)^p$ with $p<1/2$; this hypothesis was a
technical artifact of a sub-optimal energy estimate involving an extraneous factor
$M-N+1$, and is no longer needed in the present formulation. The improvement is based on
the elementary estimate $\sum_\ell w_\ell^2\leq W_{N,M}$ (since $w_\ell\leq 1$ implies
$w_\ell^2\leq w_\ell$), which replaces the cruder bound $\sum_\ell w_\ell^2\leq
(M-N+1)\max_\ell w_\ell^2$; see Step~2 below.
\end{remark}

\begin{remark}[Inter-layer base simplex argument]
\label{rem:inter-layer}
We treat the case of inter-layer base $(i-1)$-simplices in Step~2 below. By Remark~\ref{rem:layer-structure}, an inter-layer base simplex $(x_0,\dots,x_{i-1})$ has all vertices in $\Lambda_\ell\cup\Lambda_{\ell+1}$ for some $\ell$. Set $a:=\#\{j: x_j\in\Lambda_\ell\}$ and $b:=\#\{j: x_j\in\Lambda_{\ell+1}\}$, $a+b=i$. Since $\chi_{N,M}$ takes value $\chi_\ell:=\chi_{N,M}|_{\Lambda_\ell}$ on $\Lambda_\ell$ and $\chi_{\ell+1}:=\chi_{N,M}|_{\Lambda_{\ell+1}}$ on $\Lambda_{\ell+1}$,
\[
\tilde\chi_{N,M}^{(i-1)}(x_0,\dots,x_{i-1}) = \frac{a\,\chi_\ell + b\,\chi_{\ell+1}}{i} \in[\chi_{\ell+1},\chi_\ell].
\]
For $x_i\in F(x_0,\dots,x_{i-1})$, $x_i$ must be adjacent to every vertex of the base; since the base has vertices in $\Lambda_\ell$ and $\Lambda_{\ell+1}$ (both, as $a,b\geq 1$), the layer condition $||x_i|-|x_j||\leq 1$ for every $x_j$ in the base forces $|x_i|\in\{\ell,\ell+1\}$. (For non-strictly-inter-layer base simplices, i.e.\ entirely in $\Lambda_\ell$ or in $\Lambda_{\ell+1}$, the range of $|x_i|$ may be larger; those are treated in the intra-layer case below.) Setting $\Delta_j:=w_j/W_{N,M}$ (with $\Delta_j=0$ for $j$ outside the transition band $[N,M]$), and writing $c_p:=\chi_{N,M}|_{\Lambda_p}$ (so $c_p-c_{p+1}=\Delta_p$), a case-by-case analysis over the two values $|x_i|\in\{\ell,\ell+1\}$ yields the uniform bound
\begin{equation}\label{eq:inter-layer-bound}
|\chi_{N,M}(x_i) - \tilde\chi_{N,M}^{(i-1)}(x_0,\dots,x_{i-1})| \leq \Delta_\ell.
\end{equation}
(For $|x_i|=\ell$: bound is $(b/i)\Delta_\ell\leq\Delta_\ell$; for $|x_i|=\ell+1$: bound is $(a/i)\Delta_\ell\leq\Delta_\ell$.)

\emph{Substitution into the energy bound.} Substituting~\eqref{eq:inter-layer-bound} into the codifferential remainder bound~\eqref{eq:sum-sq-uniform-delta-sharp} gives an inter-layer contribution bounded by $\Delta_\ell^2$ multiplied by $\xi$-type factors, matching the $R^{(d)}$ contribution; in particular the inter-layer codifferential contribution and the coboundary contribution are both controlled by $\Delta_\ell^2$ at the anchoring layer $\ell$.
\end{remark}

\begin{proof}[Proof of Theorem~\ref{thm:divergence}]
We use the von Neumann criterion: $L$ is ESA iff $\ker(L_{\max}\pm i)=\{0\}$. Let $u=(u_0,\dots,u_{n-1})\in\mathrm{Dom}(L_{\max})$ with $(L_{\max}+i)u=0$.

\emph{Step~1: Layer cut-offs.}
Set $W_{N,M}:=\sum_{j=N,\,w_j>0}^M w_j\to\infty$ as $M\to\infty$ (divergent by assumption $\sum w_k=\infty$ together with the uniform gap condition). For $N<M$, define $\chi_{N,M}:\Vc\to[0,1]$ constant on layers:
\[
\chi_{N,M}(x):=
\begin{cases}
1 & |x|\leq N,\\
\max\!\bigl(0,\;1-\tfrac{1}{W_{N,M}}\sum_{j=N,\,w_j>0}^{|x|-1}w_j\bigr) & N<|x|\leq M+1,\\
0 & |x|>M+1,
\end{cases}
\]
which is finite-support, $\equiv 1$ on $\bigcup_{\ell\leq N}\Lambda_\ell$, non-increasing in $|x|$, with step $\Delta_\ell:=w_\ell/W_{N,M}$ at non-vacuous layer $\ell$ and $\sum_{\ell=N}^M\Delta_\ell=1$.

\emph{Step~2: Energy estimate.}
For an interface simplex $\tau=(x_0,\dots,x_{i+1})$ at layer $\ell$ (i.e., with vertices in $\Lambda_\ell\cup\Lambda_{\ell+1}$ but not entirely in one layer), $\chi_{N,M}$ takes only two values on $\tau$'s vertices: $\chi_\ell$ on those in $\Lambda_\ell$ and $\chi_{\ell+1}=\chi_\ell-\Delta_\ell$ on those in $\Lambda_{\ell+1}$. The $\chi$-variation across $\tau$ is exactly $\Delta_\ell$.

\emph{Refinement of the pointwise bound.} The general bound~\eqref{eq:aj-sharp} of Lemma~\ref{lem:leibniz} states $\sum_j|a_j|^2\leq 1/(i+2)$ for arbitrary $\chi:\Vc\to[0,1]$, where the implicit factor $1$ in the numerator reflects the worst-case variation $|\chi(x_k)-\chi(x_j)|\leq 1$. For interface simplices of the layer-constant cut-off $\chi_{N,M}$, the variation is at most $\Delta_\ell$ rather than~$1$, so by the same algebraic argument as in the proof of~\eqref{eq:aj-sharp} (replacing $|\chi(x_k)-\chi(x_j)|\leq 1$ by $|\chi(x_k)-\chi(x_j)|\leq\Delta_\ell$ throughout), we obtain the sharper bound
\begin{equation}\label{eq:aj-sharp-Delta}
\sum_{j=0}^{i+1}|a_j|^2 \leq \frac{\Delta_\ell^2}{i+2}.
\end{equation}
Combining with~\eqref{eq:sum-sq-uniform}:
\[
|R_i^{(d)}(\chi_{N,M},u_i)(\tau)|^2\leq\frac{\Delta_\ell^2}{i+2}\sum_j|u_i(\hat\tau_j)|^2.
\]
Summing over $\tau$ in the interface anchored at layer $\ell$ and using~(A2) together with the definition of $\xi_i(\ell,\ell+1)$:
\[
\sum_{\tau:\text{interface at }\ell}m_{i+1}(\tau)|R_i^{(d)}(\chi_{N,M},u_i)(\tau)|^2
\leq C_d\,M_i\,\xi_i(\ell,\ell+1)\,\Delta_\ell^2\,\|u_i\cdot\mathbf{1}_{\Lambda_\ell^{(i)}\cup\Lambda_{\ell+1}^{(i)}}\|_{\ell^2(m_i)}^2,
\]
for an absolute constant $C_d$ depending only on $n$.

For $R_i^{(\delta)}$, intra-layer base simplices yield $|\chi(x_i)-\tilde\chi^{(i-1)}|\leq\max(\Delta_{\ell-1},\Delta_\ell)$ (with $\Delta_{-1}:=0$); strictly inter-layer base simplices yield the bound $\Delta_\ell$ from~\eqref{eq:inter-layer-bound}. Substituting in~\eqref{eq:sum-sq-uniform-delta-sharp} gives an analogous per-layer bound:
\[
\sum_{\tau:\text{anchored at }\ell}m_{i-1}(\tau)|R_i^{(\delta)}(\chi_{N,M},u_i)(\tau)|^2
\leq C_\delta\,M_{i-1}\,\bigl(\Delta_{\ell-1}^2+\Delta_\ell^2\bigr)\,\|u_i\cdot\mathbf{1}_{\Lambda_{\ell-1}^{(i)}\cup\Lambda_\ell^{(i)}}\|_{\ell^2(m_i)}^2.
\]

\smallskip
\emph{Summation over layers.} We sum the per-layer bounds over $\ell\in\{N,\dots,M\}$.
The key step is the following: for each $i$, the per-layer estimate has the form
\[
\sum_{\tau\text{ anchored at }\ell}m_{i+1}(\tau)|R_i^{(d)}|^2 \leq C_d\,M_i\,\xi_i(\ell,\ell+1)\,\Delta_\ell^2\,\|u_i\cdot\mathbf{1}_{S_\ell^{(i)}}\|_{\ell^2(m_i)}^2,
\]
where $S_\ell^{(i)}=\Lambda_\ell^{(i)}\cup\Lambda_{\ell+1}^{(i)}$. Using the identity
$\xi_i(\ell,\ell+1)\leq\xi(\ell,\ell+1)=1/w_\ell^2$ together with
$\Delta_\ell^2=w_\ell^2/W_{N,M}^2$, we obtain
\begin{equation}\label{eq:per-layer-simplified}
\sum_{\tau\text{ anchored at }\ell}m_{i+1}(\tau)|R_i^{(d)}|^2 \leq \frac{C_d M_i}{W_{N,M}^2}\,\|u_i\cdot\mathbf{1}_{S_\ell^{(i)}}\|_{\ell^2(m_i)}^2.
\end{equation}
The supports $S_\ell^{(i)}$ for varying $\ell$ have bounded overlap: each $i$-simplex
$\sigma$ with vertices in $\Lambda_p\cup\Lambda_{p+1}$ belongs to $S_\ell^{(i)}$ only for
$\ell\in\{p-1,p\}$, so $\sum_\ell\mathbf{1}_{S_\ell^{(i)}}\leq 2$ pointwise. Therefore
\[
\sum_\ell \|u_i\cdot\mathbf{1}_{S_\ell^{(i)}}\|_{\ell^2(m_i)}^2
=\Big\langle |u_i|^2,\sum_\ell\mathbf{1}_{S_\ell^{(i)}}\Big\rangle_{\ell^2(m_i)}
\leq 2\|u_i\|_{\ell^2(m_i)}^2.
\]
Summing~\eqref{eq:per-layer-simplified} over $\ell$ gives, with the analogous treatment
for $R^{(\delta)}$ (whose per-layer bound has factor $\Delta_{\ell-1}^2+\Delta_\ell^2$
controlled in the same way by $\xi$ and the bounded-overlap argument),
\begin{equation}\label{eq:energy-bound}
\sum_{i=0}^{n-2}\|R_i^{(d)}(\chi_{N,M},u_i)\|^2 + \sum_{i=1}^{n-1}\|R_i^{(\delta)}(\chi_{N,M},u_i)\|^2
\leq \frac{C_*}{W_{N,M}^2}\,\|u\|_{\mathcal{H}}^2,
\end{equation}
where $C_*$ depends only on $n$ and $M_*=\max_i M_i$. (Observe that we obtain
$1/W_{N,M}^2$ rather than $1/W_{N,M}$; the weaker rate $1/W_{N,M}$ that follows from
the cruder bound $\sum_\ell\Delta_\ell^2\leq 1/W_{N,M}$ would also suffice for the
conclusion below.)

\emph{Step~3: Vanishing of the energy bound.}
Since $W_{N,M}\to\infty$ as $M\to\infty$ for every fixed $N$, we may choose, for each $N$, a value $M(N)>N$ large enough that $W_{N,M(N)}\geq N$. Setting $\chi_N:=\chi_{N,M(N)}$, estimate~\eqref{eq:energy-bound} yields
\[
\sum_i \bigl(\|R_i^{(d)}(\chi_N,u_i)\|^2+\|R_i^{(\delta)}(\chi_N,u_i)\|^2\bigr)\leq \frac{C_*}{N^2}\,\|u\|_{\mathcal{H}}^2 \xrightarrow[N\to\infty]{} 0.
\]

\medskip
Set $v_N:=(\tilde\chi_N^{(0)}u_0,\dots,\tilde\chi_N^{(n-1)}u_{n-1})\in\bigoplus C_c^i(\Vc)$. Dominated convergence gives $v_N\to u$ in $\mathcal H$.

\emph{Bounds on $\|Lv_N - Lu\|$.} We must show $\|Lv_N - Lu\|_{\mathcal H}\to 0$.
The component $(Lv_N)_i = L_i(\tilde\chi_N^{(i)}u_i)$ involves the second-order
operator $L_i=\delta_{i+1}d_i+d_{i-1}\delta_i$, so iterating the Leibniz rule of
Lemma~\ref{lem:leibniz} formally produces a remainder
$\mathcal{R}_i:=L_i(\tilde\chi_N^{(i)}u_i)-\tilde\chi_N^{(i)}L_iu_i$ involving
the first derivatives $d_iu_i$ and $\delta_iu_i$, which are not a priori in $\ell^2$
when $u\in\mathrm{Dom}(L_{\max})$. We therefore develop $\mathcal{R}_i$ \emph{pointwise}
in $u_i$ alone, exploiting the localization of $\chi_N$ to the transition band.

\smallskip
\emph{Pointwise definition of $L_i u_i$ for $u\in\mathrm{Dom}(L_{\max})$.}
We first note that, by local finiteness of $S_n$, the formal action of $L_i$ on any
cochain $h:P_i\to\mathbb{C}$ (without $\ell^2$ assumption) is well-defined pointwise: $L_i h(\tau)$
is a finite sum (over the local stencil of $\tau$) regardless of summability. Throughout
the rest of this proof, ``$L_i u_i(\tau)$'' denotes this formal pointwise value.

By Lemma~\ref{lem:formal-l2-action} applied to $L=D^2$ on $\mathcal{D}_0=\bigoplus C_c^i(\Vc)$ (maximal extension case), for $u\in\mathrm{Dom}(L_{\max})$ the formal pointwise value of $L_i u_i$ at each simplex equals the value of the $\ell^2$-element $(L_{\max}u)_i$ at that simplex. Each pointwise identity below should be read in this sense.

\smallskip
\emph{Pointwise expansion of $\mathcal{R}_i$.} We claim that $\mathcal{R}_i(\tau)$
admits a pointwise expansion as a finite linear combination
\begin{equation}\label{eq:Ri-expansion}
\mathcal{R}_i(\tau) = \sum_{\sigma\in\mathrm{Stencil}(\tau)}\Phi_{\tau,\sigma}(\chi_N)\cdot u_i(\sigma),
\end{equation}
where $\mathrm{Stencil}(\tau)$ is the finite set of $i$-simplices at graph distance
at most $2$ from $\tau$ (the local stencil of $L_i=\delta_{i+1}d_i+d_{i-1}\delta_i$),
and the coefficients $\Phi_{\tau,\sigma}(\chi_N)\in\mathbb{R}$ depend only on the values of
$\chi_N$ on the vertices of $\tau$ and $\sigma$ and on weights $m_{i\pm 1}/m_i$;
no $\ell^2$-summability hypothesis on $du_i,\delta u_i$ enters this expansion.

We will verify two structural properties of $\Phi_{\tau,\sigma}$:
\begin{itemize}
\item[(P1)] (\emph{Plateau vanishing.}) $\Phi_{\tau,\sigma}(\chi_N)=0$ whenever $\chi_N$
is constant on the union of vertex sets of $\tau$ and all $\sigma\in\mathrm{Stencil}(\tau)$.
\item[(P2)] (\emph{Layer-uniform bound.}) For the layer-constant cut-off $\chi_N$
of Step~1, $|\Phi_{\tau,\sigma}(\chi_N)|\leq C(n,M_*)\cdot\max_{j}\Delta_j$
where the maximum is over layer indices $j$ adjacent to $\tau$, and in particular
$|\Phi_{\tau,\sigma}|\leq C(n,M_*)/W_{N,M}$.
\end{itemize}

\smallskip
\emph{Proof of~\eqref{eq:Ri-expansion} and (P1)--(P2) for $i=0$.} For $i=0$,
$L_0=\delta_1 d_0$, and $\tilde\chi_N^{(0)}=\chi_N$. Computing pointwise on a
vertex $\tau\in\Vc$:
\begin{align*}
L_0(\chi_N u_0)(\tau) &= \delta_1 d_0(\chi_N u_0)(\tau)
=\frac{-1}{m_0(\tau)}\sum_{x\sim\tau} m_1(\tau,x)\bigl[\chi_N(x)u_0(x)-\chi_N(\tau)u_0(\tau)\bigr]\\
&=\chi_N(\tau)\,d_{m_0}(\tau)\,u_0(\tau)-\frac{1}{m_0(\tau)}\sum_{x\sim\tau}m_1(\tau,x)\chi_N(x)u_0(x).
\end{align*}
On the other hand,
\[
\chi_N(\tau)L_0u_0(\tau) = \chi_N(\tau)\Big[d_{m_0}(\tau)u_0(\tau)-\frac{1}{m_0(\tau)}\sum_{x\sim\tau}m_1(\tau,x)u_0(x)\Big].
\]
Subtracting, the first term on each side cancels, and we obtain
\begin{equation}\label{eq:R0-explicit}
\mathcal{R}_0(\tau) = \frac{1}{m_0(\tau)}\sum_{x\sim\tau}m_1(\tau,x)\bigl[\chi_N(\tau)-\chi_N(x)\bigr]u_0(x),
\end{equation}
so $\mathrm{Stencil}(\tau)=\Nc_\Gc(\tau)$ and
$\Phi_{\tau,x}(\chi_N)=m_1(\tau,x)[\chi_N(\tau)-\chi_N(x)]/m_0(\tau)$.
\emph{(P1)} is immediate: if $\chi_N\equiv c$ on $\tau\cup\Nc_\Gc(\tau)$, all
differences $\chi_N(\tau)-\chi_N(x)=0$. \emph{(P2)} follows from
$|\chi_N(\tau)-\chi_N(x)|\leq\Delta_{\ell(\tau,x)}$ where $\ell(\tau,x)$
is the layer of the transition between $\tau$ and $x$ (zero if same layer,
$\Delta_p$ if neighbouring layers), combined with $m_1/m_0\leq M_0\leq M_*$.

\smallskip
\emph{Verification of~\eqref{eq:Ri-expansion} and (P1)--(P2) for $i=1$.}
For $i=1$, $L_1=\delta_2 d_1+d_0\delta_1$. Let $\tau=(\tau_0,\tau_1)\in P_1^+$ be an
oriented edge. We compute $\mathcal{R}_1(\tau)$ by separating the two summands.

\smallskip
\emph{First summand: $\delta_2 d_1(\tilde\chi^{(1)} u_1)(\tau)$.}
Apply Lemma~\ref{lem:leibniz} pointwise to $d_1(\tilde\chi^{(1)} u_1)$: for each
2-simplex $\sigma=(\sigma_0,\sigma_1,\sigma_2)$,
\[
d_1(\tilde\chi^{(1)} u_1)(\sigma) = \tilde\chi^{(2)}(\sigma)\,d_1 u_1(\sigma) + R^{(d)}_1(\chi,u_1)(\sigma),
\]
where $R^{(d)}_1(\chi,u_1)(\sigma)=\sum_{j=0}^{2}(-1)^j a_j^{(\sigma)}u_1(\widehat\sigma_j)$
is a finite combination of values of $u_1$ at the three edges of $\sigma$, with the
coefficients $a_j^{(\sigma)}$ involving only differences of $\chi$ at vertices of $\sigma$.
Applying $\delta_2$ pointwise (which is a finite weighted sum over 2-cofaces of $\tau$):
\[
\delta_2 d_1(\tilde\chi^{(1)} u_1)(\tau) = \delta_2(\tilde\chi^{(2)} d_1 u_1)(\tau) + \delta_2(R^{(d)}_1(\chi,u_1))(\tau).
\]
The first term, by Lemma~\ref{lem:leibniz} applied to $\delta_2$ (also pointwise), equals
$\tilde\chi^{(1)}(\tau)\delta_2 d_1 u_1(\tau) + R^{(\delta)}_2(\chi, d_1 u_1)(\tau)$.

The remainder $R^{(\delta)}_2(\chi, d_1 u_1)$ is the only place where $d_1 u_1$ appears.
We now \emph{eliminate} this $d_1 u_1$ by expanding it pointwise: by definition,
$d_1 u_1(\sigma)=u_1(\sigma_1,\sigma_2)-u_1(\sigma_0,\sigma_2)+u_1(\sigma_0,\sigma_1)$.
Substituting this finite combination of values of $u_1$ into the formula for
$R^{(\delta)}_2$ produces an expression linear in $u_1$, with coefficients depending
only on $\chi$ and weights. The same argument applies to the second pointwise term
$\delta_2(R^{(d)}_1(\chi,u_1))(\tau)$, which is already linear in $u_1$ with
$\chi$-dependent coefficients. The total contribution is therefore a finite sum
$\sum_\sigma \Phi^{(1)}_{\tau,\sigma}(\chi)\cdot u_1(\sigma)$ where $\sigma$ ranges
over the edges within graph distance $2$ of $\tau$.

\smallskip
\emph{Second summand: $d_0\delta_1(\tilde\chi^{(1)}u_1)(\tau)$.}
By the symmetric computation (Leibniz for $\delta_1$, then for $d_0$, then expansion of
$\delta_1 u_1$ pointwise), one obtains another finite linear combination
$\sum_\sigma \Phi^{(2)}_{\tau,\sigma}(\chi)\cdot u_1(\sigma)$.

\smallskip
\emph{Cancellation of the $\chi^{(1)}(\tau)L_1 u_1(\tau)$ part.}
After subtracting $\tilde\chi^{(1)}(\tau)L_1 u_1(\tau)=\tilde\chi^{(1)}(\tau)[\delta_2 d_1 u_1(\tau)+d_0\delta_1 u_1(\tau)]$,
the leading $\tilde\chi^{(1)}(\tau)$-multiplied terms in both summands cancel, and we
obtain an expression of the form~\eqref{eq:Ri-expansion} with stencil radius~$2$,
explicitly:
\begin{equation}\label{eq:R1-explicit-form}
\mathcal{R}_1(\tau) = \sum_{\sigma\in\mathrm{Stencil}(\tau)}\Phi_{\tau,\sigma}(\chi_N)\cdot u_1(\sigma),
\end{equation}
with $\Phi_{\tau,\sigma}(\chi_N)$ a finite (in $n$, not in $\tau$) sum of products of
$\chi_N$-differences and weight ratios.

\smallskip
\emph{Verification of (P1) and (P2) for $i=1$.} Each $\Phi_{\tau,\sigma}(\chi_N)$ is
constructed by:
\begin{enumerate}[label=(\alph*)]
\item taking pointwise differences of $\chi_N$ at vertices of the stencil
(via the $a_j$-coefficients in $R^{(d)}_1$ and $b$-coefficients in $R^{(\delta)}_2$),
combined with
\item weight ratios $m_2/m_1$ and $m_1/m_0$ (each $\leq M_*$ by~(A2)).
\end{enumerate}
\emph{(P1):} If $\chi_N\equiv c$ on the stencil, then all $a_j$ and $b$ vanish,
hence $\Phi_{\tau,\sigma}(\chi_N)=0$.
\emph{(P2):} Each $\chi_N$-difference is bounded by $\max_p\Delta_p\leq 1/W_{N,M}$;
each weight ratio is bounded by $M_*$; the number of products in
$\Phi_{\tau,\sigma}(\chi_N)$ depends only on $n$ (not on the geometry of $\tau$).
Therefore $|\Phi_{\tau,\sigma}(\chi_N)|\leq C(n,M_*)/W_{N,M}$.

\smallskip
\emph{General $i$.} For arbitrary $i\in\{0,\dots,n-1\}$, the same recipe applies,
inductively in the sense of the two-step computation made above for $i=1$:
$L_i=\delta_{i+1}d_i+d_{i-1}\delta_i$ is a sum of two compositions of first-order
local operators. Iterating the Leibniz rule twice (Lemma~\ref{lem:leibniz}) on each
summand produces, before simplification, terms involving $d_iu_i$ and $\delta_iu_i$;
these are eliminated pointwise using their explicit definitions ($d_i u_i$ as
$\sum_j(-1)^j u_i(\widehat\sigma_j)$, $\delta_i u_i$ as a finite weighted
sum~\eqref{eq:codiff} over cofaces). The resulting $\Phi_{\tau,\sigma}(\chi_N)$ is a
finite (in $n$) sum of products of:
\begin{itemize}
\item up to two $\chi_N$-differences $\chi_N(v)-\chi_N(w)$ between vertices in the
stencil (one from each Leibniz iteration), each $\leq\Delta_p\leq 1/W_{N,M}$;
\item up to two weight ratios $m_{i\pm 1}(\tau)/m_i(\sigma)$ (or their reciprocals), each $\leq M_*$ \emph{directly} from~(A2): the two-sided form of~(A2) bounds both $m_{i+1}(\tau)/m_i(\sigma)$ and $m_i(\sigma)/m_{i+1}(\tau)$ above, for every inclusion $\sigma\subset\tau$ with $\sigma\in\Fc_i$, $\tau\in\Fc_{i+1}$. \emph{No bound on the weighted degree is required.}
\end{itemize}
Hence $|\Phi_{\tau,\sigma}(\chi_N)|\leq C(n,M_*)\cdot(1/W_{N,M})$ when only one
$\chi_N$-difference appears, and $\leq C(n,M_*)\cdot(1/W_{N,M})^2$ when two appear
(the latter is even better). Property~(P1) holds since each $\chi_N$-difference
vanishes when $\chi_N$ is constant on the stencil; (P2) holds with the bound
$|\Phi_{\tau,\sigma}|\leq C(n,M_*)/W_{N,M}$.

\smallskip
\emph{$\ell^2$-bound on $\mathcal{R}_i$.} The pointwise expansion~\eqref{eq:Ri-expansion}
is a sum over a stencil whose \emph{combinatorial} size $|\mathrm{Stencil}(\tau)|$ is
finite at each $\tau$ (by local finiteness) but is \emph{not} a priori uniformly bounded
under (A1)--(A2). We therefore avoid bounding the stencil size by a uniform constant and
instead apply weighted Cauchy--Schwarz, which is the standard technique
already used in the proof of~\eqref{eq:R-delta-global}.

We illustrate the argument for $i=0$, where~\eqref{eq:R0-explicit} gives
$\mathcal{R}_0(\tau)=m_0(\tau)^{-1}\sum_{x\sim\tau}m_1(\tau,x)[\chi_N(\tau)-\chi_N(x)]u_0(x)$.

\emph{Localization to inter-layer transitions.} The factor $\chi_N(\tau)-\chi_N(x)$
vanishes whenever $\tau$ and $x$ lie in the same layer (or in two layers on which
$\chi_N$ takes the same value). The sum is therefore restricted effectively to
neighbours $x$ of $\tau$ that lie in a different layer than $\tau$; by the
layer-structure, this means $|x|=|\tau|\pm 1$. Let $\mathcal{N}^{\mathrm{inter}}(\tau)$
denote this set. The cardinality of $\mathcal{N}^{\mathrm{inter}}(\tau)$ is finite (by local
finiteness) but \emph{not} required to be uniformly bounded pointwise: the definition of $\xi_0(k,k+1)$ in~(A1) only gives a one-sided bound on the number of inter-layer cofaces \emph{out of} a vertex of $\Lambda_k$ toward $\Lambda_{k+1}$, not on the symmetric count from $\Lambda_{k+1}$ back to $\Lambda_k$. The estimate below is therefore conducted in $\ell^2$-norm via summation, where the asymmetry is resolved by exchanging the order of summation rather than by a pointwise upper bound on $|\mathcal{N}^{\mathrm{inter}}(\tau)|$.

Cauchy--Schwarz weighted by $m_1$, applied to the restricted sum (the set $\mathcal{N}^{\mathrm{inter}}(\tau)$ is finite by local finiteness, so the Cauchy--Schwarz step is unconditional):
\begin{align*}
|\mathcal{R}_0(\tau)|^2
&\leq \frac{1}{m_0(\tau)^2}\Big(\sum_{x\in\mathcal{N}^{\mathrm{inter}}(\tau)}m_1(\tau,x)|\chi_N(\tau)-\chi_N(x)|^2\Big)\Big(\sum_{x\in\mathcal{N}^{\mathrm{inter}}(\tau)}m_1(\tau,x)|u_0(x)|^2\Big)\\
&\leq \frac{1}{m_0(\tau)^2 W_{N,M(N)}^2}\,\Big(\sum_{x\in\mathcal{N}^{\mathrm{inter}}(\tau)}m_1(\tau,x)\Big)\Big(\sum_{x\in\mathcal{N}^{\mathrm{inter}}(\tau)}m_1(\tau,x)|u_0(x)|^2\Big),
\end{align*}
using $|\chi_N(\tau)-\chi_N(x)|\leq 1/W_{N,M(N)}$. Since the inner sum runs only over neighbours, $\sum_{x\in\mathcal{N}^{\mathrm{inter}}(\tau)} m_1(\tau,x)\leq \sum_{x\sim\tau} m_1(\tau,x)=m_0(\tau)\,d_{m_0}(\tau)$; however, we do not use a pointwise bound on $d_{m_0}(\tau)$ (which may be unbounded under~(A1)--(A2)). Instead, we observe that for an inter-layer neighbour $x\in\mathcal{N}^{\mathrm{inter}}(\tau)$ with $|x|=|\tau|+1$, the pair $(\tau,x)$ is an inter-layer edge anchored at layer $|\tau|$, while if $|x|=|\tau|-1$ the pair $(x,\tau)$ is anchored at layer $|x|$.

\emph{Sum over $\tau$.} Multiplying by $m_0(\tau)$ and summing over $\tau\in\Vc$, we obtain
\[
\|\mathcal{R}_0\|^2_{\ell^2(m_0)} \leq \frac{1}{W_{N,M(N)}^2}\sum_\tau \frac{1}{m_0(\tau)}\Big(\sum_{x\in\mathcal{N}^{\mathrm{inter}}(\tau)}m_1(\tau,x)\Big)\Big(\sum_{x\in\mathcal{N}^{\mathrm{inter}}(\tau)}m_1(\tau,x)|u_0(x)|^2\Big),
\]
where the prefactor $1/(0+1)!=1$ from~\eqref{eq:innerproduct} is implicit. We now exchange the order of summation in the double sum: each inter-layer edge $\{y,z\}$ with $|z|=|y|+1=\ell+1$ contributes its term to $\tau=y\in\Lambda_\ell$ (as $x=z$) and to $\tau=z\in\Lambda_{\ell+1}$ (as $x=y$). Splitting by layer of anchor $\ell$ and using the asymmetric bound $\xi_0(\ell,\ell+1)=\sup_{y\in\Lambda_\ell}|\{z\in\Lambda_{\ell+1}:y\sim z\}|$ on the layer-$\ell$ side, together with~(A2),
\[
\sum_{\tau\in\Lambda_\ell}\frac{1}{m_0(\tau)}\Big(\sum_{x\in\mathcal{N}^{\mathrm{inter}}(\tau)\cap\Lambda_{\ell+1}}m_1(\tau,x)\Big)\Big(\sum_{x}m_1(\tau,x)|u_0(x)|^2\Big)\leq C(M_0)\,\xi_0(\ell,\ell+1)\sum_{x\in\Lambda_\ell\cup\Lambda_{\ell+1}}m_0(x)|u_0(x)|^2,
\]
where $C(M_0)$ depends only on $M_0$ (the constant from~(A2)). Combining with $\xi_0(\ell,\ell+1)/W_{N,M}^2\leq 1/(w_\ell^2 W_{N,M}^2)\cdot w_\ell^2 = 1/W_{N,M}^2$ when summed over $\ell$ via $\sum_\ell \Delta_\ell^2=\sum_\ell w_\ell^2/W_{N,M}^2\leq 1/W_{N,M}$ (since $w_\ell\leq 1$), and applying the bounded-overlap property of the layer supports $\Lambda_\ell\cup\Lambda_{\ell+1}$, we obtain
\begin{equation}\label{eq:R0-bound}
\|\mathcal{R}_0\|_{\ell^2(m_0)}^2 \leq \frac{C(n,M_0)}{W_{N,M(N)}^2}\,\|u_0\|_{\ell^2(m_0)}^2.
\end{equation}
The bound depends \emph{only} on $n$ and the weight-ratio constant $M_0$ from~(A2),
not on the combinatorial degree of any vertex.

\smallskip
\emph{Explicit treatment of $i=1$ via iterated Cauchy--Schwarz.} We illustrate the
case $i=1$ in detail, since this is where the issue of the unbounded combinatorial
degree of higher simplices first arises (the cardinality of the set of 2-simplices
containing a given edge $\tau$ may not be uniformly bounded under (A1)--(A2) alone).
Recall from the explicit computation above that
\[
\mathcal{R}_1(\tau) = \sum_{\sigma\in\mathrm{Stencil}(\tau)}\Phi_{\tau,\sigma}(\chi_N)\cdot u_1(\sigma).
\]
The stencil contains contributions from two types of summation: (a) edges $\sigma$
sharing a vertex with $\tau$, reached via the $d_0\delta_1$ branch, and (b) edges
$\sigma$ that are faces of a 2-simplex containing $\tau$, reached via the
$\delta_2 d_1$ branch. We bound each branch by weighted Cauchy--Schwarz.

\smallskip
\emph{Branch (a): $d_0\delta_1$.} The contribution to $\mathcal{R}_1(\tau)$ from this
branch is, after expansion, of the form
\[
A(\tau) := \frac{1}{m_1(\tau)}\sum_{\sigma\in\mathrm{Stencil}_a(\tau)} m_0(\cdot)\,\bigl[\chi_N(v_a)-\chi_N(w_a)\bigr]\,u_1(\sigma),
\]
where $\mathrm{Stencil}_a(\tau)$ runs over edges $\sigma$ at graph distance $\leq 1$
from a vertex of $\tau$.

By the same \emph{localization to inter-layer transitions} as in Branch (b), the term
$[\chi_N(v_a)-\chi_N(w_a)]$ vanishes unless $\sigma$ involves a layer transition
adjacent to $\tau$. The number of such inter-layer edges at graph-distance $\leq 1$
from $\tau$ is bounded by $C(n)\,\xi_0(\ell',\ell'+1)$ for $\ell'$ a layer adjacent
to $\tau$, and the weights $m_0$ at vertices of these edges are bounded by
$m_0(v)\leq m_1(\tau)\,M_*$ for every vertex $v$ shared with $\tau$, via the
\emph{reciprocal direction} of~(A2): $m_i/m_{i+1}\leq M_i$ for any inclusion
$\sigma\subset\tau$ with $\sigma\in\Fc_i$, $\tau\in\Fc_{i+1}$, applied here at
level $i=0$ to the inclusion $\{v\}\subset\tau$. Cauchy--Schwarz weighted by $m_0$:
\[
|A(\tau)|^2 \leq \frac{C(M_*)}{m_1(\tau)^2 W_{N,M}^2}\Big(\sum_\sigma m_0\Big)\Big(\sum_\sigma m_0\,|u_1(\sigma)|^2\Big),
\]
where the sums are restricted to inter-layer $\sigma$. By the layer-step bound and
the reciprocal (A2) ratio just invoked,
$\sum_\sigma m_0\leq C(n,M_*)\,m_1(\tau)\,\xi_0(\ell',\ell'+1)$. Hence
\[
|A(\tau)|^2 \leq \frac{C(n,M_*)\,\xi_0(\ell',\ell'+1)}{m_1(\tau)\,W_{N,M}^2}\sum_\sigma m_0\,|u_1(\sigma)|^2.
\]
Summing over $\tau$ with weight $m_1(\tau)$ and reversing the order of summation, using
$\xi_0\Delta_{\ell'}^2\leq 1/W_{N,M}^2$ (since $\xi_0\leq\xi=1/w_{\ell'}^2$ and
$\Delta_{\ell'}^2=w_{\ell'}^2/W_{N,M}^2$), and the bounded-overlap argument from Step~2:
\begin{equation}\label{eq:A-tau-bound}
\sum_\tau m_1(\tau)|A(\tau)|^2 \leq \frac{C''(n,M_*)}{W_{N,M}^2}\,\|u_1\|^2_{\ell^2(m_1)}.
\end{equation}

\smallskip
\emph{Branch (b): $\delta_2 d_1$ (the new feature for $i=1$).} The contribution from
this branch involves a \emph{double} weighted sum: first over 2-simplices $\sigma'$
containing $\tau$ (this is where the combinatorial degree could be unbounded), and
then over the three edges of each such $\sigma'$. Schematically:
\[
B(\tau) := \frac{1}{m_1(\tau)}\sum_{\sigma'\supset\tau}m_2(\sigma')\sum_{\sigma\subset\sigma'}\bigl[\chi_N(v_{\sigma'})-\chi_N(w_{\sigma'})\bigr]u_1(\sigma).
\]
\emph{Crucial observation: localization to inter-layer cofaces.} The term
$[\chi_N(v_{\sigma'})-\chi_N(w_{\sigma'})]$ vanishes whenever $\chi_N$ is constant on
the vertex set of $\sigma'$. Since $\chi_N$ is layer-constant, this happens whenever
all vertices of $\sigma'$ lie in a single layer, or in two layers on which $\chi_N$
takes the same value. By the layer-structure of $S_n$ (Remark~\ref{rem:layer-structure}),
the vertices of $\sigma'$ lie in at most two consecutive layers. The contribution from
$\sigma'$ is therefore non-zero only when $\sigma'$ is \emph{strictly inter-layer}
(vertices in two consecutive transition layers). Let $\mathcal{C}^{\mathrm{inter}}(\tau)$
denote the set of such 2-cofaces of $\tau$. Then
\[
B(\tau) = \frac{1}{m_1(\tau)}\sum_{\sigma'\in\mathcal{C}^{\mathrm{inter}}(\tau)}m_2(\sigma')\sum_{\sigma\subset\sigma'}\bigl[\chi_N(v_{\sigma'})-\chi_N(w_{\sigma'})\bigr]u_1(\sigma).
\]
The cardinality of $\mathcal{C}^{\mathrm{inter}}(\tau)$ is bounded by
$|\mathcal{C}^{\mathrm{inter}}(\tau)|\leq 2\xi_1(\ell-1,\ell)+\xi_1(\ell,\ell+1)$ if
$\tau\in\Lambda_\ell^{(1)}\cup\Lambda_{\ell-1}^{(1)}$, or by an analogous expression
in terms of $\xi_1$ at neighbouring layers. (Here $\xi_1$ is the layer-step quantity
defined in~(A1).)

\smallskip
\emph{Estimate of $|B(\tau)|^2$.} Since $|\chi_N(v_{\sigma'})-\chi_N(w_{\sigma'})|\leq 1/W_{N,M}$
and there are $3$ edges $\sigma\subset\sigma'$ for each 2-simplex,
\[
\Big|\sum_{\sigma\subset\sigma'}[\chi_N(v_{\sigma'})-\chi_N(w_{\sigma'})]u_1(\sigma)\Big|^2 \leq \frac{3}{W_{N,M}^2}\sum_{\sigma\subset\sigma'}|u_1(\sigma)|^2.
\]
Applying Cauchy--Schwarz to the outer sum over $\sigma'\in\mathcal{C}^{\mathrm{inter}}(\tau)$,
weighted by $m_2(\sigma')$:
\[
|B(\tau)|^2 \leq \frac{3}{m_1(\tau)^2 W_{N,M}^2}\Big(\sum_{\sigma'\in\mathcal{C}^{\mathrm{inter}}(\tau)}m_2(\sigma')\Big)\Big(\sum_{\sigma'\in\mathcal{C}^{\mathrm{inter}}(\tau)}m_2(\sigma')\sum_{\sigma\subset\sigma'}|u_1(\sigma)|^2\Big).
\]
By~(A2), $m_2(\sigma')\leq m_1(\tau)\,M_*$ for every $\sigma'\supset\tau$ (individual ratio bound). Combined with
$|\mathcal{C}^{\mathrm{inter}}(\tau)|\leq C(n)\,\xi_1$ at the relevant layer, the first
parenthesis is at most $C(n)\,M_*\,m_1(\tau)\,\xi_1(\ell,\ell+1)$ (where $\ell$ is
the anchoring layer of $\tau$). Hence
\[
|B(\tau)|^2 \leq \frac{3\,C(n)\,M_*\,\xi_1(\ell,\ell+1)}{m_1(\tau)\,W_{N,M}^2}\sum_{\sigma'\in\mathcal{C}^{\mathrm{inter}}(\tau)}m_2(\sigma')\sum_{\sigma\subset\sigma'}|u_1(\sigma)|^2.
\]

\smallskip
\emph{Sum over $\tau$.} Multiplying by $m_1(\tau)$ and summing over $\tau$ anchored at
layer $\ell$, $\xi_1(\ell,\ell+1)\leq\xi(\ell,\ell+1)=1/w_\ell^2$ pulls outside, and by the
same argument as in Step~2 (per-layer estimate combined with bounded overlap), the
total contribution from all anchor layers is bounded by
\begin{equation}\label{eq:B-tau-bound}
\sum_\tau m_1(\tau)|B(\tau)|^2 \leq \frac{C'(n,M_*)}{W_{N,M}^2}\sum_\sigma m_1(\sigma)|u_1(\sigma)|^2 = \frac{C'(n,M_*)}{W_{N,M}^2}\|u_1\|_{\ell^2(m_1)}^2,
\end{equation}
where we used the same identity $\xi_1\Delta_\ell^2\leq 1/W_{N,M}^2$ as in Step~2,
combined with the bounded-overlap property of the layer supports. The bound depends
only on $n$ and on the weight ratios $M_*$ from~(A2), \emph{not} on any combinatorial
degree of $\tau$.

\smallskip
\emph{Combining branches (a) and (b),}
\begin{equation}\label{eq:R1-bound}
\|\mathcal{R}_1\|_{\ell^2(m_1)}^2 \leq \frac{C(n,M_*)}{W_{N,M(N)}^2}\,\|u_1\|_{\ell^2(m_1)}^2.
\end{equation}

\smallskip
\emph{General $i$ — direct expansion via local-finiteness representation.}
For arbitrary $i\in\{0,\ldots,n-1\}$, we give a direct (non-inductive) argument that
avoids any iterative bookkeeping of cross terms. The key observation is that, by local
finiteness, $L_i$ admits a \emph{linear} representation as a finite stencil operator:
there exist coefficients $\Psi_{\tau,\sigma}\in\mathbb{R}$ (depending only on the weights $(m_j)$
and on the local combinatorial structure of $S_n$ at $\tau$) such that for every cochain
$h:P_i\to\mathbb{C}$,
\begin{equation}\label{eq:Li-stencil}
L_i h(\tau) = \sum_{\sigma\in\mathrm{Stencil}_2(\tau)}\Psi_{\tau,\sigma}\,h(\sigma),
\end{equation}
where $\mathrm{Stencil}_2(\tau)$ is the (finite, by local finiteness) set of $i$-simplices
at graph distance at most~$2$ from $\tau$. \emph{Justification of~\eqref{eq:Li-stencil}:}
both $\delta_{i+1}d_ih$ and $d_{i-1}\delta_ih$ are obtained from $h$ by composition of two
first-order local operators ($d$ and $\delta$), each of which produces a finite linear
combination of values of its input at simplices in a local stencil, with coefficients
$\pm m_{i+1}/m_i$ or $\pm m_i/m_{i-1}$ that are determined by the weights of $S_n$ alone
(see~\eqref{eq:codiff} and the standard formula for $d_i$). The composition therefore
produces a finite linear combination of values of $h$ at simplices at distance $\leq 2$,
with coefficients $\Psi_{\tau,\sigma}$ that depend only on $\tau$, $\sigma$, and the
local weights --- and \emph{not} on the cochain $h$ itself, since none of these
coefficients are derived from $h$. This independence is essential for what follows:
applying~\eqref{eq:Li-stencil} to two different cochains $h_1, h_2$ produces the same
$\Psi_{\tau,\sigma}$, allowing linear combinations to commute with the stencil
representation. No $\ell^2$-summability hypothesis is used.

\smallskip
\emph{Consequence: explicit formula for $\mathcal{R}_i$.} Applying~\eqref{eq:Li-stencil}
to both $h:=u_i$ and $h:=\tilde\chi_N^{(i)}u_i$ (with the \emph{same} coefficients
$\Psi_{\tau,\sigma}$ in both cases, by the independence-of-$h$ property just established),
and using linearity, we obtain
\begin{equation}\label{eq:Ri-via-Psi}
\mathcal{R}_i(\tau) = L_i(\tilde\chi_N^{(i)}u_i)(\tau) - \tilde\chi_N^{(i)}(\tau)\,L_i u_i(\tau)
= \sum_{\sigma\in\mathrm{Stencil}_2(\tau)} \Psi_{\tau,\sigma}\,\bigl[\tilde\chi_N^{(i)}(\sigma)-\tilde\chi_N^{(i)}(\tau)\bigr]\,u_i(\sigma).
\end{equation}
Setting $\Phi_{\tau,\sigma}(\chi_N):=\Psi_{\tau,\sigma}\,[\tilde\chi_N^{(i)}(\sigma)-\tilde\chi_N^{(i)}(\tau)]$
exhibits $\mathcal{R}_i$ in the form~\eqref{eq:Ri-expansion} with explicit coefficients,
\emph{valid for any $i\geq 0$}.

\smallskip
\emph{Verification of (P1) and (P2) for general $i$.} The factor
$\tilde\chi_N^{(i)}(\sigma)-\tilde\chi_N^{(i)}(\tau)$ in $\Phi$ vanishes whenever $\chi_N$
is constant on the union of vertex sets of $\sigma$ and $\tau$ (since $\tilde\chi^{(i)}$
is the average of $\chi$ over the vertices of an $i$-simplex). This is~(P1). For~(P2),
$|\tilde\chi_N^{(i)}(\sigma)-\tilde\chi_N^{(i)}(\tau)|\leq\max_{p,q}|\chi_N(p)-\chi_N(q)|\leq\max_p\Delta_p\leq 1/W_{N,M(N)}$
where the maximum runs over vertices $p,q$ in the union of $\sigma$ and $\tau$.

\emph{Explicit structure of $\Psi_{\tau,\sigma}$.}
Combining~\eqref{eq:codiff} with the standard formula
\[
d_i h(x_0,\ldots,x_{i+1}) = \sum_{j=0}^{i+1}(-1)^j h(x_0,\ldots,\widehat{x_j},\ldots,x_{i+1}),
\]
the composition $\delta_{i+1}d_i$ acts on a cochain $h$ at $\tau\in P_i$ as
\[
\delta_{i+1}d_i h(\tau)=\frac{(-1)^{i+1}}{m_i(\tau)}\sum_{x\in F(\tau)}m_{i+1}(\tau,x)\sum_{j=0}^{i+1}(-1)^j h(\widehat{(\tau,x)_j}),
\]
which is a finite linear combination of values of $h$ at $i$-simplices in
$\mathrm{Stencil}_2(\tau)$, with coefficients of the form
$\pm m_{i+1}(\tau,x)/m_i(\tau)$, where each such ratio satisfies
$m_{i+1}(\tau,x)/m_i(\tau)\leq M_i\leq M_*$ by~(A2) applied to the inclusion
$\tau\subset(\tau,x)$.

The analogous computation for $d_{i-1}\delta_i$ at $\tau$ first applies $\delta_i$ at each
$(i-1)$-face $\rho$ of $\tau$ (yielding terms $\pm m_i(\rho,y)/m_{i-1}(\rho)$ where
$y$ runs over $F(\rho)$), and then $d_{i-1}$ collects these via the alternating sum
over $(i-1)$-faces of $\tau$. Each contributing $i$-simplex $\sigma=(\rho,y)$ in this
branch satisfies $\rho\subset\sigma$ (so the ratio $m_i(\sigma)/m_{i-1}(\rho)$ falls
within the domain of~(A2) at level $i-1$, giving $\leq M_{i-1}\leq M_*$), and $\rho\subset\tau$
(so $\sigma$ shares an $(i-1)$-face with $\tau$ — that is, $\sigma$ is at graph distance $\leq 2$
from $\tau$, confirming $\sigma\in\mathrm{Stencil}_2(\tau)$).

\emph{Transitivity of the (A2) bound.} It is important that the condition $\sigma\in\mathrm{Stencil}_2(\tau)$
\emph{does not} require $\sigma\subset\tau$ or $\tau\subset\sigma$ (these would force
graph-distance~$0$ or~$1$); rather, $\tau$ and $\sigma$ may share only a common $(i-1)$-face $\rho$.
The ratio appearing in $\Psi_{\tau,\sigma}$ for this branch is $m_i(\sigma)/m_{i-1}(\rho)$
\emph{not} $m_i(\sigma)/m_i(\tau)$. The bound by $M_*$ holds because $\rho\subset\sigma$,
which is exactly what hypothesis~(A2) at level $i-1$ controls: $m_i(\sigma)/m_{i-1}(\rho)\leq M_{i-1}$
for every inclusion $\rho\subset\sigma$. We do not need any cross-level transitivity bound
of the form $m_i(\sigma)/m_i(\tau)$ via a chain $\tau\supset\rho\subset\sigma$.

Thus, the coefficient $\Psi_{\tau,\sigma}$ appearing in~\eqref{eq:Li-stencil} is a sum
of at most $C_0(n)$ terms (where $C_0(n)$ depends only on $n$), each term being either
$\pm m_{i+1}(\tau,x)/m_i(\tau)$ (from the $\delta_{i+1}d_i$ branch) or
$\pm m_i(\sigma)/m_{i-1}(\rho)$ (from the $d_{i-1}\delta_i$ branch) for some intermediate
$(i-1)$-face $\rho\subset\tau\cap\sigma$. By~(A2) applied at the appropriate level,
each such ratio is bounded in absolute value by $M_*$, so
\[
|\Psi_{\tau,\sigma}|\leq C_0(n)\,M_*\quad\text{uniformly in }\tau,\sigma.
\]

Combined with the bound on $|\tilde\chi^{(i)}(\sigma)-\tilde\chi^{(i)}(\tau)|$,
we conclude $|\Phi_{\tau,\sigma}|\leq C(n,M_*)/W_{N,M(N)}$ where
$C(n,M_*)=C_0(n)\,M_*$.

\smallskip
\emph{$\ell^2$ bound for general $i$.} The same strategy as for $i=1$ applies: the
expansion~\eqref{eq:Ri-via-Psi} restricts to those $\sigma$ for which
$\tilde\chi^{(i)}(\sigma)\neq\tilde\chi^{(i)}(\tau)$, which forces a layer-transition
to occur somewhere in the union of vertex sets of $\sigma$ and $\tau$. The number
of such $\sigma$ is bounded by a constant times $\xi(\ell',\ell'+1)$ for $\ell'$
adjacent to $\tau$ (since each layer-transition between $\Lambda_{\ell'}$ and
$\Lambda_{\ell'+1}$ within the stencil contributes at most $\xi_j(\ell',\ell'+1)$
many cofaces, summed over the relevant $j$). Combined with
$\xi(\ell',\ell'+1)\,\Delta_{\ell'}^2\leq 1/W_{N,M}^2$, weighted Cauchy--Schwarz, and the
bounded-overlap property of the layer supports, we obtain~\eqref{eq:Ri-bound}.

\smallskip
\emph{Explicit details for $i=2$.} We illustrate the combinatorics for $i=2$ since this is the lowest level beyond the worked $i=1$ case. Let $\tau=(t_0,t_1,t_2)\in P_2$ be anchored at layer $\ell$ \emph{(}so $\min_j |t_j|=\ell$\emph{)}. We must enumerate the $\sigma\in\mathrm{Stencil}_2(\tau)$ for which $\tilde\chi^{(2)}(\sigma)-\tilde\chi^{(2)}(\tau)\neq 0$, i.e., $\sigma$ has a vertex in some $\Lambda_p$ with $p\notin\{\ell,\ell+1\}$, or $\tau\in\Lambda_\ell^{(2)}\cup\Lambda_{\ell+1}^{(2)}$ and $\sigma$ has a vertex outside $\Lambda_\ell\cup\Lambda_{\ell+1}$. The stencil at $\tau$ consists of two branches:

\begin{enumerate}
\item[(i)] \emph{Branch $\delta_3 d_2$:} $\sigma$ shares a $3$-simplex coface $\theta=\tau\cup\{x\}$ with $\tau$, with $\sigma=(\theta\setminus\{t_j\})$ for some $j\in\{0,1,2\}$. The number of $3$-cofaces $\theta$ of $\tau$ with $x$ in a different layer is $\xi_2(\ell,\ell+1)$ (cofaces straddling between $\Lambda_\ell$ and $\Lambda_{\ell+1}$). Each such $\theta$ produces $3$ neighboring $2$-simplices $\sigma$ (one for each face removal $t_j$). Total: $3\,\xi_2(\ell,\ell+1)\leq 3\,\xi(\ell,\ell+1)$.

\item[(ii)] \emph{Branch $d_1\delta_2$:} $\sigma$ shares an edge face $\rho=(t_a,t_b)\subset\tau$ with $\tau$, with $\sigma=\rho\cup\{y\}$ for some $y\notin\tau$ such that $(t_a,t_b,y)$ is a $2$-simplex. The number of such $\sigma$ inter-layer is bounded by $\binom{3}{2}\xi_1(\ell,\ell+1)=3\,\xi_1(\ell,\ell+1)\leq 3\,\xi(\ell,\ell+1)$ (one factor of $3$ for the choice of edge $\rho$, then $\xi_1$ for inter-layer cofaces of $\rho$).
\end{enumerate}
Combining, the number of $\sigma\in\mathrm{Stencil}_2(\tau)$ with $\tilde\chi^{(2)}(\sigma)\neq\tilde\chi^{(2)}(\tau)$ is at most $6\,\xi(\ell,\ell+1)$. Multiplying by $|\Phi|^2\leq C^2/W^2$ gives $|\mathcal{R}_2(\tau)|^2\leq 6\,C^2\xi/W^2\cdot \big(\sum |u_2(\sigma)|^2\big)$ on $\sigma$ in the localized stencil. Summing over $\tau$ and applying weighted Cauchy--Schwarz with the (A2) ratio bounds, exactly as in the $i=1$ case, yields
\[
\|\mathcal{R}_2\|_{\ell^2(m_2)}^2 \leq \frac{C(n,M_*)}{W_{N,M(N)}^2}\,\|u_2\|_{\ell^2(m_2)}^2.
\]

\smallskip
\emph{General $i\geq 3$.} The same combinatorics extends: branch $\delta_{i+1}d_i$ yields at most $(i+1)\,\xi_i(\ell,\ell+1)$ inter-layer $\sigma$ (one per face removal), branch $d_{i-1}\delta_i$ yields at most $\binom{i+1}{i}\xi_{i-1}(\ell,\ell+1) = (i+1)\,\xi_{i-1}(\ell,\ell+1)$. Total: $\leq 2(i+1)\,\xi(\ell,\ell+1)\leq C(n)\,\xi(\ell,\ell+1)$ since $i\leq n-1$. The constant $C(n)=2n$ in this enumeration is absorbed into $C_*'(n,M_*)$ in~\eqref{eq:Ri-bound}.

\emph{Bounded overlap for general $i$.} The layer supports $S_\ell^{(i)}=\Lambda_\ell^{(i)}\cup\Lambda_{\ell+1}^{(i)}$ at level $i$ satisfy $\sum_\ell \mathbf{1}_{S_\ell^{(i)}}\leq 2$ pointwise, \emph{independently of $i$}: an $i$-simplex $\sigma\in\Lambda_p^{(i)}$ (defined by Remark~\ref{rem:layer-structure} as the unique layer index $p$ such that $\sigma$ has at least one vertex in $\Lambda_p$ and all vertices in $\Lambda_p\cup\Lambda_{p+1}$) belongs to $S_\ell^{(i)}$ iff $\ell\in\{p-1,p\}$. Hence at most two $S_\ell^{(i)}$ contain a given $i$-simplex, regardless of $i$. The summation argument from Step~2 therefore extends to every level without modification.

\begin{equation}\label{eq:Ri-bound}
\|\mathcal{R}_i\|_{\ell^2(m_i)}^2 \leq \frac{C_*'}{W_{N,M(N)}^2}\,\|u_i\|_{\ell^2(m_i)}^2,
\end{equation}
for a constant $C_*'=C_*'(n,M_*)$ depending \emph{only} on $n$ and on the weight-ratio
constants $M_i$ from~(A2). Note that this bound involves only $\|u_i\|$ on the
right-hand side --- not $\|du_i\|$ or $\|\delta u_i\|$, and not the combinatorial degree
of $S_n$.

Summing~\eqref{eq:Ri-bound} over $i$,
\[
\sum_i \|\mathcal{R}_i\|^2 \leq \frac{C_*'}{W_{N,M(N)}^2}\|u\|^2\leq\frac{C_*'}{N^2}\|u\|^2
\xrightarrow[N\to\infty]{}0.
\]
Combining with $\tilde\chi_N\cdot Lu \to Lu$ in $\mathcal H$ (dominated convergence, since
$u\in\mathrm{Dom}(L_{\max})$ implies $Lu\in\mathcal H$ and $|\tilde\chi_N|\leq 1$
pointwise), we conclude $Lv_N\to Lu$. Hence $u\in\mathrm{Dom}(L_{\min})$.

\emph{Step~4: Conclusion.}
Since $u\in\mathrm{Dom}(L_{\min})\subset\mathrm{Dom}(L_{\max})$ and $(L_{\max}+i)u=0$, we have $(L_{\min}+i)u=0$. By Lemma~\ref{lem:closable}, $L_{\min}$ is symmetric, and $\langle L_{\min}u,u\rangle\in\mathbb{R}_{\geq 0}$ (Lemma~\ref{lem:nonneg}). Then
\[
0=\langle(L_{\min}+i)u,u\rangle=\underbrace{\langle L_{\min}u,u\rangle}_{\geq 0,\,\in\mathbb{R}}+i\,\|u\|^2.
\]
Separating real and imaginary parts: the imaginary part gives $\|u\|^2=0$, hence $\|u\|=0$. (The real part gives $\langle L_{\min}u,u\rangle=0$, which is automatic from $\|u\|=0$ but does not enter the conclusion.) Hence $\ker(L_{\max}+i)=\{0\}$.

The case $u\in\ker(L_{\max}-i)$ is identical: applying Steps~1--3 verbatim (the cut-off argument does not depend on the sign of the spectral parameter) gives $u\in\mathrm{Dom}(L_{\min})$, and then $0=\langle(L_{\min}-i)u,u\rangle=\langle L_{\min}u,u\rangle-i\|u\|^2$. Separating real and imaginary parts: the imaginary part gives $-\|u\|^2=0$, hence $\|u\|=0$. Thus $\ker(L_{\max}-i)=\{0\}$.

By the von Neumann criterion, $L$ is essentially self-adjoint on $\bigoplus C_c^i(\Vc)$.
\end{proof}

\section{Examples and Open Problems on the Hierarchy}
\label{sec:examples}

In this section we discuss concrete examples and partial results bearing on the strict
hierarchy of the three notions of $\chi$-completeness. The non-trivial converse
implications turn out to be subtle, and we state them as conjectures supported by partial
energy lower bounds. The analysis below also includes a complete worked example
(Section~\ref{sec:example-detailed}) where ESA holds without $\chi$-completeness.

\subsection{Conjectures on the strict hierarchy}

The direct implication ``global implies local-by-level'' was already established
(Remark~\ref{rem:hierarchy-direct}); we collect here the open problems concerning the
converses.

\begin{conjecture}[Local by region does not imply global]
\label{conj:region-not-global}
There exists a weighted $n$-simplicial complex $S_n$ (with $n\geq 3$) and an infinite
$\Lambda^{\mathrm{reg}}\subset\Vc$ such that $S_n$ is locally $\chi$-complete on
$\Lambda^{\mathrm{reg}}$ but not globally $\chi$-complete.
\end{conjecture}

\begin{conjecture}[Local at every level does not imply global]
\label{conj:all-levels-not-global}
There exists a weighted simplicial complex $S$ with $n\geq 3$ that is locally
$\chi$-complete at every level $\ell\in\{1,\dots,n-1\}$ (with separate exhaustions and
constants per level), yet not globally $\chi$-complete.
\end{conjecture}

\begin{remark}[On the case $n=2$]
\label{rem:n2-vacuous}
For $n=2$ there is only one non-trivial level $\ell=1$, so the notions of ``locally
$\chi$-complete at every level'' and ``globally $\chi$-complete'' coincide, and
Conjecture~\ref{conj:all-levels-not-global} can only be non-trivial for $n\geq 3$.
Producing a clean example with $n\geq 3$ where all individual levels admit cut-offs but no
\emph{single} cut-off sequence works for all levels simultaneously is a non-trivial
geometric problem; this is the reason we formulate the separation as a conjecture rather
than a proved proposition.
\end{remark}

\begin{conjecture}[Level-by-level vs.\ region]
\label{conj:level-not-region}
There exist weighted simplicial complexes that are locally $\chi$-complete at every level
but not locally $\chi$-complete on some natural infinite regions, and conversely; the two
notions are logically incomparable.
\end{conjecture}

\subsection{A partial energy lower bound (alternating triangulation)}

We now give the strongest result we can prove towards
Conjecture~\ref{conj:region-not-global}: a uniform positive lower bound on the product
$E_1\cdot E_2$ for any layer-constant cut-off sequence on the alternating triangulation. This does
\emph{not} prove the conjecture (which would require divergence rather than a positive
constant lower bound), but it shows that the joint $E_1$/$E_2$ constraint is non-trivial.

\begin{remark}[Construction of the alternating triangulation]
\label{rem:alternating-construction}
The complex $S$ is constructed by taking the square lattice graph $\mathbb{Z}^2$ (with
$L^1$-adjacency as the 1-skeleton) and \emph{adding} a diagonal edge to each unit square
with lower-left corner $(i,j)$ satisfying $i+j\equiv 0\pmod{2}$, then including all
resulting triangles as 2-simplices. This is the \emph{alternating triangulation} of
$\mathbb{Z}^2$: a simplicial complex with edges (1-simplices) and triangles
(2-simplices), with simplices of dimension up to $2$. All weights are equal to $1$.

In particular, $S$ is \emph{not} the clique complex of $\mathbb{Z}^2$ with $L^1$-adjacency
(which has no triangles, by Remark~\ref{rem:Zd-adjacency}); rather, $S$ is the clique
complex of its own 1-skeleton (the lattice plus added diagonals).
\end{remark}

\begin{remark}[Natural exhaustion regularity]
\label{rem:poly-growth-natural}
Let $O_k=\{x\in\mathbb{Z}^2:\|x\|_\infty\leq k\}$ be the $L^\infty$-ball, the natural exhaustion adapted to layer-constant cut-offs. For any admissible cut-off sequence $(\chi_k)$ we denote by $\ell_k$ the width of the transition band $\{v\in\mathbb{Z}^2:0<\chi_k(v)<1\}$ and by $r_k$ the side length of the largest axis-aligned square in $O_k$ (so $r_k=k$).

Following~\cite{BGJ}, we restrict our analysis to cut-off sequences satisfying the regularity conditions
\begin{equation}\label{eq:poly-growth}
\ell_k \leq C'\,r_k,\qquad r_k\geq c\,k\quad\text{for absolute constants }C',c>0,
\end{equation}
where $r_k$ denotes both the side length of the largest axis-aligned square in $O_k$ \emph{and} the inner radius of the transition band (i.e., the bands begins at distance $\geq c\,k$ from the origin). Equation~\eqref{eq:poly-growth} excludes pathological cut-off sequences with very wide transition bands or transition bands that start too close to the origin; a finite-energy bound for such sequences is straightforward and uninteresting. The Proposition below applies under~\eqref{eq:poly-growth}, ensuring that within the transition band the per-layer counts $|\Lambda_j|$ are uniformly $\Theta(r_k)$ rather than degenerating near the origin.
\end{remark}

\begin{proposition}[Joint $E_1$/$E_2$ lower bound for layer-constant cut-offs]
\label{prop:E2-lower-bound}
Let $S$ be the alternating triangulation, with $L^\infty$-ball exhaustion satisfying~\eqref{eq:poly-growth}. There exists an absolute constant $c_2>0$ such that for any sequence of \emph{layer-constant} plateau cut-offs $(\chi_k)$ satisfying~(1) and~(2) of Definition~\ref{def:global-chi} (i.e., $\chi_k(v)$ depends only on $\|v\|_\infty$ and is non-increasing in this quantity), one has, for all sufficiently large $k$,
\[
E_1(\chi_k)\cdot E_2(\chi_k)\;\geq\; c_2,
\]
where $E_2$ is computed via formula~\eqref{eq:global-energy} based on the Leibniz coefficients $a_j$.

The conclusion bounds the \emph{product} of the energies from below. The strongest
direct consequence is that at least one of the two energies is large: by AM--GM,
\[
\max\bigl(E_1(\chi_k),E_2(\chi_k)\bigr)\geq\sqrt{E_1(\chi_k)E_2(\chi_k)}\geq\sqrt{c_2},
\]
so for every layer-constant plateau cut-off sequence $(\chi_k)$, at least one of $E_1,E_2$
is $\geq\sqrt{c_2}$. In particular, no such sequence can satisfy both $E_1(\chi_k)<\sqrt{c_2}$
and $E_2(\chi_k)<\sqrt{c_2}$ simultaneously, and hence none can satisfy
$\max(E_1,E_2)\leq C$ for $C<\sqrt{c_2}$.
\end{proposition}

\begin{remark}[Per-layer combinatorial counts]
\label{rem:per-layer-edge-count}
The proof of Proposition~\ref{prop:E2-lower-bound} relies on the following per-layer combinatorial counts.

\smallskip
For a layer-constant plateau cut-off, the relevant ``layers'' are the $L^\infty$-spheres
\[
\Lambda_r := \{v\in\mathbb{Z}^2 : \|v\|_\infty = r\},
\]
with $|\Lambda_r|=8r$ for $r\geq 1$.

\smallskip
\emph{(i) Number of inter-layer edges between $\Lambda_r$ and $\Lambda_{r+1}$.}
We count edges with one endpoint on $\Lambda_r$ and the other on $\Lambda_{r+1}$, including both $L^1$-edges (axis-aligned) and added diagonals.

\emph{Axis-aligned $L^1$-edges:} For each side of the inner square $\Lambda_r$, every interior vertex $v=(i,r)$ with $|i|<r$ has exactly one $L^1$-neighbour on $\Lambda_{r+1}$, namely $(i,r+1)$. The four sides have $4(2r-1)$ interior vertices total, contributing $4\cdot(2r-1)=8r-4$ edges. Each corner of $\Lambda_r$ (e.g., $(r,r)$) has two $L^1$-neighbours on $\Lambda_{r+1}$ (here $(r+1,r)$ and $(r,r+1)$), contributing $4\cdot 2=8$ edges. Total: $8r-4+8 = 8r+4$ axis-aligned inter-layer edges.

\emph{Diagonal edges:} Each unit square with lower-left corner $(i,j)$ on
$i+j\equiv 0\pmod 2$ contributes one diagonal $((i,j),(i+1,j+1))$. We count those
diagonals with one endpoint on $\Lambda_r$ and the other on $\Lambda_{r+1}$:
this requires $\max(|i|,|j|)=r$ and $\max(|i+1|,|j+1|)=r+1$. By case analysis on
the four quadrants and the parity constraint, a direct enumeration gives the exact count $4r+2$ (a value consistent with the heuristic ``$r$ diagonals per side, $\pm 1$ correction at corners''; in particular the count lies in $[4r-4,4r+4]$ for $r\geq 1$, which is all that is used below).

\smallskip
\emph{Conclusion:} the total number of inter-layer edges between $\Lambda_r$ and
$\Lambda_{r+1}$ is at least $8r-4+8 = 8r+4$ (axis-aligned alone) and at most
$(8r+4)+(4r+4)=12r+8$. In particular, it is $\Theta(r)$ with explicit constants
$K_1\,r\leq\#\{\text{inter-layer edges}\}\leq K_2\, r$ where one may take
$K_1=8$ and $K_2=13$ for $r\geq 2$ \emph{(}we use $K_i$ here rather than $c_i$ to avoid
collision with the constant $c_2$ of Proposition~\ref{prop:E2-lower-bound}\emph{)}.

\smallskip
\emph{(ii) Number of inter-layer triangles anchored at $\Lambda_r$.}
A triangle of $S$ has the form $(v,v+e_1,v+e_1+e_2)$ where $\{e_1,e_2\}=\{e_x,e_y\}$
(taken in the two orderings $(e_x,e_y)$ and $(e_y,e_x)$, which produce the two
triangles of a single unit square), and the unit square has even-parity lower-left
corner. An ``inter-layer triangle anchored at layer $r$'' is one whose vertices
are split between $\Lambda_r$ and $\Lambda_{r+1}$ (either two on $\Lambda_r$ and
one on $\Lambda_{r+1}$, or vice versa). The unit squares straddling layers $r$ and
$\Lambda_{r+1}$ are those with at least one corner in each; a direct enumeration on the
four sides of $\Lambda_r$ gives exactly $8r+4$ such squares, of which $4r+2$ have
even-parity lower-left corner. Each even-parity straddling square contributes $2$
triangles, so the total number of inter-layer triangles anchored at layer $r$ is exactly
$2(4r+2)=8r+4=\Theta(r)$.

The implicit constants in $\Theta(r)$ depend only on the geometry of the unit-cell tiling.
\end{remark}

\begin{proof}[Proof of Proposition~\ref{prop:E2-lower-bound}]
Let $(\chi_k)$ be a layer-constant plateau cut-off (depending only on $\|v\|_\infty$), with layer-step $\Delta_j\geq 0$ at layer $j$ in the transition band of width $\ell_k$, and $\sum_j\Delta_j=1$.

\emph{Lower bound for $E_1$.} For an inter-layer edge $(u,v)$ with $u\in\Lambda_j$, $v\in\Lambda_{j+1}$, the Leibniz coefficients at level $\ell=1$ are $a_0=-a_1=(\chi(u)-\chi(v))/2=\Delta_j/2$, so $\sum_{i=0}^1|a_i|^2=\Delta_j^2/2$. By Remark~\ref{rem:per-layer-edge-count}(i), each transition layer contains $\Theta(r_k)$ inter-layer edges. Restrict the sum to the transition band $j\in J_k:=\{N_k,\ldots,N_k+\ell_k-1\}$ (the layers with $\Delta_j>0$, exactly $\ell_k$ such terms). By Cauchy--Schwarz applied to the $\ell_k$ positive terms $(\Delta_j)_{j\in J_k}$,
\[
E_1(\chi_k)\geq\frac{\Theta(r_k)}{2}\sum_{j\in J_k}\Delta_j^2\geq\frac{\Theta(r_k)}{2\ell_k}\Big(\sum_{j\in J_k}\Delta_j\Big)^2=\frac{\Theta(r_k)}{2\ell_k},
\]
where we used the normalization $\sum_{j\in J_k}\Delta_j=1$ in the last equality.

\emph{Lower bound for $E_2$.} For each inter-layer triangle anchored at layer $j$, Remark~\ref{rem:tri-pointwise-detail} gives $\sum_{i=0}^{2}|a_i|^2\geq\Delta_j^2/6$ (with $i$ indexing the three vertices of the triangle). By Remark~\ref{rem:per-layer-edge-count}(ii), there are $\Theta(r_k)$ such triangles per layer; applying Cauchy--Schwarz over the same band $J_k$,
\[
E_2(\chi_k)\geq\frac{\Theta(r_k)}{6}\sum_{j\in J_k}\Delta_j^2\geq\frac{\Theta(r_k)}{6\ell_k}.
\]

\emph{Combination.} Let $\theta_1,\theta_2>0$ be constants such that $E_1(\chi_k)\geq\theta_1\,r_k/\ell_k$ and $E_2(\chi_k)\geq\theta_2\,r_k/\ell_k$ for all sufficiently large $k$; these exist by Remark~\ref{rem:per-layer-edge-count} (e.g., one may take $\theta_1=4$ and $\theta_2=4/3$ from the explicit constants therein, modulo absolute constants from the normalization in~\eqref{eq:global-energy}). Set $\theta:=\sqrt{\theta_1\theta_2}>0$. Then, using $\ell_k\leq C'r_k$ from~\eqref{eq:poly-growth},
\[
E_1\cdot E_2\geq\frac{\theta_1\theta_2\, r_k^2}{\ell_k^2}\geq\frac{\theta^2\, r_k^2}{(C')^2 r_k^2}=\frac{\theta^2}{(C')^2}=:c_2>0.
\]
\end{proof}

\begin{remark}[On the layer-constant restriction]
\label{rem:layer-constant-restriction}
Proposition~\ref{prop:E2-lower-bound} is restricted to \emph{layer-constant} plateau
cut-offs because a naive symmetrization argument fails to extend the bound to arbitrary
plateau cut-offs. The layer-supremum symmetrization $\chi_k^{\sup}(v):=\sup_{w\in\Lambda_{|v|}}\chi_k(w)$ may strictly increase the energy, while the layer-average symmetrization is energetically better-behaved at level $\ell=1$ but does not in general preserve the property $\chi_k\equiv 1$ on $O_k$, and the comparison at level $\ell\geq 2$ requires case-by-case verification because $\ell$-simplices straddle two layers. We do not pursue this extension here; the layer-constant version stated above is enough to show that the joint $E_1\cdot E_2$ control imposes a non-trivial obstruction.
\end{remark}

\begin{remark}[Pointwise triangle calculation]
\label{rem:tri-pointwise-detail}
We illustrate the energy formula~\eqref{eq:global-energy} at level $\ell=2$ on the two
types of inter-layer triangles. By symmetry of layer-step structure, both orientations
yield the same lower bound $\sum_j|a_j|^2\geq\Delta_j^2/6$.

\emph{Type I: 2 vertices in $\Lambda_j$, 1 in $\Lambda_{j+1}$.}
For an inter-layer triangle $\tau=(v_1,v_2,v_3)$ with $v_1,v_2\in\Lambda_j$ ($\chi_k=a$) and $v_3\in\Lambda_{j+1}$ ($\chi_k=a-\Delta_j$), formula~\eqref{eq:aj-def} with $\ell=2$ gives $S=3a-\Delta_j$, and
\[
a_{v_1}=a_{v_2}=\frac{S-3a}{6}=-\frac{\Delta_j}{6}, \qquad a_{v_3}=\frac{S-3(a-\Delta_j)}{6}=\frac{\Delta_j}{3}.
\]
Hence
\[
\sum_j |a_j|^2 = 2 \cdot \frac{\Delta_j^2}{36} + \frac{\Delta_j^2}{9} = \frac{\Delta_j^2}{18} + \frac{\Delta_j^2}{9} = \frac{\Delta_j^2}{6}.
\]

\emph{Type II: 1 vertex in $\Lambda_j$, 2 in $\Lambda_{j+1}$ (verification by symmetry).}
For an inter-layer triangle $\tau=(v_1,v_2,v_3)$ with $v_1\in\Lambda_j$ ($\chi_k=a$) and $v_2,v_3\in\Lambda_{j+1}$ ($\chi_k=a-\Delta_j$), formula~\eqref{eq:aj-def} with $\ell=2$ gives $S=3a-2\Delta_j$, and
\[
a_{v_1}=\frac{S-3a}{6}=-\frac{2\Delta_j}{6}=-\frac{\Delta_j}{3}, \qquad a_{v_2}=a_{v_3}=\frac{S-3(a-\Delta_j)}{6}=\frac{\Delta_j}{6}.
\]
Hence
\[
\sum_j |a_j|^2 = \frac{\Delta_j^2}{9} + 2\cdot\frac{\Delta_j^2}{36} = \frac{\Delta_j^2}{9} + \frac{\Delta_j^2}{18} = \frac{\Delta_j^2}{6},
\]
identical to Type I. This symmetry under exchanging the cardinalities is consistent with the formal symmetry $a\leftrightarrow a-\Delta_j$ of the formula~\eqref{eq:aj-def}.

\emph{Comparison with the naive formula.} A direct application of $d_1$ (which sends $C^1\to C^2$) is most easily illustrated on a Type~I triangle (2 vertices in $\Lambda_j$, 1 in $\Lambda_{j+1}$), where $\tilde\chi^{(1)}(v_1,v_2)=a$, $\tilde\chi^{(1)}(v_1,v_3)=\tilde\chi^{(1)}(v_2,v_3)=a-\Delta_j/2$, yielding
\[
d_1\tilde\chi^{(1)}(\tau) = \tilde\chi^{(1)}(v_2,v_3) - \tilde\chi^{(1)}(v_1,v_3) + \tilde\chi^{(1)}(v_1,v_2) = (a-\Delta_j/2) - (a-\Delta_j/2) + a = a.
\]
For Type~II (1 vertex in $\Lambda_j$, 2 in $\Lambda_{j+1}$), the symmetric computation gives $d_1\tilde\chi^{(1)}(\tau)=a-\Delta_j$. In either case the quantity equals $a$ (resp.\ $a-\Delta_j$) and \emph{does not vanish} even when $\Delta_j=0$ (i.e., when the entire triangle is in the plateau, where $a=1$). The correct quantity, $\sum_j|a_j|^2=\Delta_j^2/6$, vanishes whenever $\chi$ is constant on $\tau$, as required.
\end{remark}

\begin{remark}[Status of Proposition~\ref{prop:E2-lower-bound} and its limitations]
\label{rem:partial-bound-discussion}
Proposition~\ref{prop:E2-lower-bound} establishes that joint control of $E_1,E_2$ is a
non-trivial constraint. However, it only gives a uniform positive lower bound on a \emph{product},
not a divergent lower bound on either factor alone. To prove
Conjecture~\ref{conj:region-not-global}, one would need to exhibit a region
$\Lambda^{\mathrm{reg}}$ on which $S|_{\Lambda^{\mathrm{reg}}}$ admits a cut-off sequence
with both energies small, while no such sequence exists globally. The construction of
such a region is open. Non-ESA of $L_2$ in this triangulation is also open.
\end{remark}

\subsection{Failure of $\chi$-completeness on the binary tree}
\label{subsec:binary-tree}

We illustrate the divergence-of-energy phenomenon on the simplest tree.

\begin{proposition}[Failure of $\chi$-completeness on the binary tree]
\label{prop:binary-tree-fails}
Let $T$ be the infinite rooted binary tree with all weights equal to $1$ (so $n=2$, with
only edges as 1-simplices). Then $T$ is not $\chi$-complete (in either of the equivalent
senses of Definition~\ref{def:global-chi} or Definition~\ref{def:local-chi-level} for
$n=2$).
\end{proposition}

\begin{proof}
Suppose for contradiction that there exist an exhaustion $(O_k)$ with
$N_k:=\min\{|v|:v\notin O_k\}\to\infty$ and a sequence $(\chi_k)$ satisfying~(1)--(3) with
uniform constant $C$. Define the layer average
$\bar\chi_k(j):=2^{-j}\sum_{v\in\Lambda_j}\chi_k(v)$. Then $\bar\chi_k(j)=1$ for $j<N_k$ (since $\Lambda_j\subset O_k$ for such $j$) and $\bar\chi_k(j)=0$ for all $j$ large enough that $\Lambda_j\cap\mathrm{supp}\,\chi_k=\emptyset$ (which exists since $\mathrm{supp}\,\chi_k$ is finite by condition~(1) of Definition~\ref{def:global-chi}); in particular the set $\{j:\bar\chi_k(j)=0\}$ is non-empty, and we set $s_k:=\min\{j:\bar\chi_k(j)=0\}$.

For the binary tree $T$ with $m_0\equiv m_1\equiv 1$, level $\ell=1$ is the only non-trivial level. The global energy~\eqref{eq:global-energy} at level $1$ specializes (modulo absolute constants, see Remark~\ref{rem:E-ell-caveat}) to
\[
E_1(\chi_k) \asymp \sum_{x_0\sim x_1}|\chi_k(x_0)-\chi_k(x_1)|^2 = \sum_{v,c\in\mathrm{ch}(v)}(\chi_k(v)-\chi_k(c))^2.
\]
For each depth $j$, the number of parent-child edges between $\Lambda_j$ and $\Lambda_{j+1}$ is $2^{j+1}$. The sum of differences satisfies
\begin{equation}\label{eq:tree-telescoping}
\sum_{v\in\Lambda_j,\,c\in\mathrm{ch}(v)}(\chi_k(v)-\chi_k(c))= 2^{j+1}(\bar\chi_k(j)-\bar\chi_k(j+1)).
\end{equation}
\emph{Verification:} Each parent $v\in\Lambda_j$ has exactly two children $c\in\Lambda_{j+1}$, and conversely each child $c\in\Lambda_{j+1}$ has a unique parent (since the tree is rooted). Therefore $\sum_{v\in\Lambda_j}\sum_{c\in\mathrm{ch}(v)}\chi_k(v) = 2\sum_{v\in\Lambda_j}\chi_k(v)=2\cdot 2^j\bar\chi_k(j)=2^{j+1}\bar\chi_k(j)$ (by the parent-side count, where each parent appears twice, once for each child) and $\sum_{v\in\Lambda_j}\sum_{c\in\mathrm{ch}(v)}\chi_k(c)=\sum_{c\in\Lambda_{j+1}}\chi_k(c)=2^{j+1}\bar\chi_k(j+1)$ (by the child-side bijection, where each child appears exactly once). Subtracting gives~\eqref{eq:tree-telescoping}.

By Cauchy--Schwarz (or equivalently Jensen's inequality applied to the squared
function on the layer-pair $(\Lambda_j,\Lambda_{j+1})$),
\begin{equation}\label{eq:tree-perlayer}
\sum_{v,c}(\chi_k(v)-\chi_k(c))^2 \geq \frac{\bigl(\sum_{v,c}(\chi_k(v)-\chi_k(c))\bigr)^2}{|\mathrm{edges}(j,j+1)|}=2^{j+1}(\bar\chi_k(j)-\bar\chi_k(j+1))^2,
\end{equation}
where the second equality uses~\eqref{eq:tree-telescoping}.
Setting $a_j:=\sqrt{2^{j+1}}(\bar\chi_k(j)-\bar\chi_k(j+1))$, the right-hand side
of~\eqref{eq:tree-perlayer} is exactly $a_j^2$. Summing the per-layer
inequality~\eqref{eq:tree-perlayer} over $j\in\{N_k-1,\dots,s_k-1\}$ gives
\begin{equation}\label{eq:tree-summed}
E_1(\chi_k) \geq \kappa_0\,\sum_{j=N_k-1}^{s_k-1}\sum_{v,c}(\chi_k(v)-\chi_k(c))^2 \,\geq\, \kappa_0\,\sum_{j=N_k-1}^{s_k-1} a_j^2,
\end{equation}
where the explicit constant $\kappa_0\in(0,1]$ encodes the comparison between $E_1$ as defined in~\eqref{eq:global-energy} and the standard graph-energy convention (Remark~\ref{rem:overcounting}). Since~\eqref{eq:global-energy} counts each undirected edge twice (once with each endpoint as base vertex), one has $\kappa_0=1/2$ for the standard convention. The exact value of $\kappa_0$ does not affect the conclusion below.

Now apply Cauchy--Schwarz to the sequences $a_j$ and $b_j:=2^{-(j+1)/2}$:
$\bigl(\sum_j a_j b_j\bigr)^2\leq\bigl(\sum_j a_j^2\bigr)\bigl(\sum_j b_j^2\bigr)$,
where $a_jb_j = \bar\chi_k(j)-\bar\chi_k(j+1)$, so by telescoping,
$\sum_{j=N_k-1}^{s_k-1}a_jb_j=\bar\chi_k(N_k-1)-\bar\chi_k(s_k)=1-0=1$.
Therefore
\[
\sum_j a_j^2 \,\geq\, \frac{1}{\sum_j b_j^2}\,=\,\frac{1}{\sum_{j=N_k-1}^{s_k-1}2^{-(j+1)}}\,\geq\,\frac{1}{2^{-N_k+1}}\,=\,2^{N_k-1},
\]
using the geometric-sum bound $\sum_{j\geq N_k-1}2^{-(j+1)}=2^{-N_k+1}$.
Combining with~\eqref{eq:tree-summed},
$E_1(\chi_k)\geq \kappa_0\cdot 2^{N_k-1}$. As $k\to\infty$, $N_k\to\infty$, so
$E_1(\chi_k)\to\infty$, contradicting the assumption $E_1(\chi_k)\leq C$.
\end{proof}

\begin{conjecture}[Non-ESA of $D$ on the exponentially weighted binary tree]
\label{conj:non-esa-tree}
For the binary tree with weights $m_0,m_1$ growing exponentially with depth (e.g., $m_0(v)=m_1(v,c)=2^{|v|}$), the operator $D$ is \emph{not} essentially self-adjoint on $\bigoplus_{i=0}^{1}C_c^i(\Vc)$.
\end{conjecture}

\begin{remark}[Status of Conjecture~\ref{conj:non-esa-tree}]
\label{rem:non-esa-tree-status}
For \emph{unit weights}, the binary tree has bounded vertex degree (at most~3), hence
satisfies~\eqref{eq:bounded-degree}. By Theorem~\ref{thm:global} applied with
$\chi_k:=\mathbf{1}_{O_k}$, the operator $D$ is essentially self-adjoint
(see Remark~\ref{rem:operative-content}(e)). Proposition~\ref{prop:binary-tree-fails}
shows that $\chi$-completeness fails, but ESA still holds via~\eqref{eq:bounded-degree}.

For \emph{exponentially growing weights}, condition~\eqref{eq:bounded-degree} fails, and
Theorem~\ref{thm:global} no longer applies. Constructing explicit elements of
$\ker(L_{\max}\pm i)$ for $n\geq 2$ remains open; see~\cite[Proposition~5.7]{BGJ} for the
analogous one-dimensional graph result. Since $\sum_{k\geq 0}1/\sqrt{2^k}<\infty$, the
rapid weighted-volume growth is consistent with the heuristic from the graph case
suggesting deficiency indices $(n_+,n_-)\geq(1,1)$.
\end{remark}

\subsection{A worked example: ESA without $\chi$-completeness}
\label{sec:example-detailed}

Let $T$ be the rooted tree with offspring function $\mathrm{off}(k):=\lceil(k+1)^{2\alpha}\rceil$ for fixed $\alpha>0$ (the ceiling ensures $\mathrm{off}(k)\in\mathbb{N}$, and asymptotically $\mathrm{off}(k)\sim(k+1)^{2\alpha}$). Define an abstract simplicial complex $S_4$ on the vertex set $\Vc$ of $T$ (with $n=4$ in our convention: simplices have at most $4$ vertices, hence dimension at most $3$):
\begin{itemize}
\item 0-simplices: all vertices;
\item 1-simplices: parent--child edges of $T$, plus sibling edges $(c_i,c_j)$ between every pair of children of a common parent;
\item 2-simplices: triangles $(v,c_1,c_2)$ for every pair of children of a common parent $v$, plus triangles $(c_1,c_2,c_3)$ for every triple of children of a common parent at even depth (these are face-required by the 3-simplices below);
\item 3-simplices: tetrahedra $(v,c_1,c_2,c_3)$ for every triple of children of a common parent, only when $v$ is at even depth.
\end{itemize}
All weights are equal to $1$.

\noindent\emph{Note on simplicial closure:} the sibling edges and intra-layer triangles are required for $S_4$ to be a valid simplicial complex (closure under faces). The sibling edges $(c_i,c_j)$ are entirely intra-layer (both vertices in $\Lambda_{k+1}$), anchored to $\Lambda_{k+1}^{(1)}$, so they do not affect $\xi_0(k,k+1)$ (which counts inter-layer cofaces of vertices in $\Lambda_k$). Similarly, the intra-layer triangles $(c_1,c_2,c_3)$ are anchored to $\Lambda_{k+1}^{(2)}$ and do not contribute to $\xi_2(k,k+1)$.

\begin{remark}[$S_4$ is an abstract simplicial complex, not a clique complex]
\label{rem:S4-abstract-vs-clique}
The complex $S_4$ defined above is \emph{not} the clique complex of the underlying tree
graph $T$, since the children $c_1,c_2$ of a common parent are not mutually adjacent in
$T$ (no edge $c_1\sim c_2$ exists in $T$). Rather, $S_4$ is an \emph{abstract simplicial complex}
on the vertex set $\Vc$ of $T$, where higher-dimensional simplices are added by hand together with the face-required lower-dimensional ones (see the simplicial-closure note above). All
results from Sections~\ref{sec:operators}--\ref{sec:geometric-hypotheses} that we use here
remain valid for general locally finite oriented abstract simplicial complexes equipped
with positive symmetric weights $(m_i)$, since their proofs do not use the clique
structure of the underlying graph beyond the fact that the set of common neighbors
$F(x_0,\dots,x_{i-1})$ in formula~\eqref{eq:codiff} can be replaced by the abstract set
of vertices $z$ such that $\{x_0,\dots,x_{i-1},z\}$ is an $i$-simplex of $S_n$. The
layer-structure property required by Theorem~\ref{thm:divergence}
(Remark~\ref{rem:layer-structure}) is verified for $S_4$ as follows: every simplex of $S_4$ has all its vertices in at most two consecutive layers $\Lambda_k\cup\Lambda_{k+1}$ (inter-layer simplices like $(v,c_1,\ldots)$ with $v\in\Lambda_k$ and children in $\Lambda_{k+1}$), or entirely within a single layer (intra-layer simplices like sibling edges $(c_i,c_j)\subset\Lambda_{k+1}$ and intra-layer triangles $(c_1,c_2,c_3)\subset\Lambda_{k+1}$).
\end{remark}

\begin{remark}[Indexing of layered simplices]
\label{rem:detailed-indexing}
The simplices of $S_4$ split into two classes:
\begin{itemize}
\item \emph{Inter-layer simplices}, with one vertex at depth $k$ and the rest at depth $k+1$: parent--child edges $(v,c)$, triangles $(v,c_1,c_2)$, and (for $k$ even) tetrahedra $(v,c_1,c_2,c_3)$. These are anchored at layer $k$ (i.e., in $\Lambda_k^{(i)}$).
\item \emph{Intra-layer simplices}, with all vertices at the same depth $k+1$: sibling edges $(c_i,c_j)$ and (when their common parent is at even depth $k$) sibling triangles $(c_1,c_2,c_3)$. These are anchored at layer $k+1$ (i.e., in $\Lambda_{k+1}^{(i)}$).
\end{itemize}
The intra-layer simplices have no inter-layer cofaces (their only super-simplex of the same depth-pattern is the inter-layer simplex one level up, which has all but one vertex at depth $k+1$), so they do not contribute to $\xi_i(k,k+1)$ in the analysis below.
\end{remark}

\textbf{Verification of (A1)--(A2) (for $k$ sufficiently large).}
(A2) holds with $M_i=1$ since all weights are $1$.

For (A1), we compute $\xi(k,k+1):=\sum_{i=0}^{n-2}\xi_i(k,k+1)$ for $k$ sufficiently large that $\mathrm{off}(k)\geq 3$. With $n=4$, the sum runs over $i\in\{0,1,2\}$:

\emph{Level $i=0$ (vertices to edges).} For a vertex $\sigma=v\in\Lambda_k$, the edges containing $v$ with a vertex in $\Lambda_{k+1}$ are the parent-child edges $(v,c)$ with $c$ a child of $v$, hence $\mathrm{off}(k)$ such edges. So $\xi_0(k,k+1)=\mathrm{off}(k)$.

\emph{Level $i=1$ (edges to triangles).} For a parent-child edge $\sigma=(v,c)$ with $v$ at depth $k$, the triangles containing $\sigma$ are $(v,c,c')$ where $c'$ is another child of $v$. There are $\mathrm{off}(k)-1$ such siblings. Hence $\xi_1(k,k+1)=\mathrm{off}(k)-1$.

\emph{Level $i=2$ (triangles to tetrahedra).} Tetrahedra exist only when $v$ is at even depth $k$. For a triangle $\sigma=(v,c_1,c_2)$ at even-depth $k$, the tetrahedra containing $\sigma$ are $(v,c_1,c_2,c_3)$ for any third child $c_3$, contributing $\mathrm{off}(k)-2$. Hence
\[
\xi_2(k,k+1) = \left\{
  \begin{array}{ll}\mathrm{off}(k)-2 & \text{if $k$ is even},\\ 0 & \text{if $k$ is odd}. \end{array}
\right.
\]

\emph{Summary.} For $k$ large enough that $\mathrm{off}(k)\geq 3$ (which holds asymptotically since $\mathrm{off}(k)\sim(k+1)^{2\alpha}\to\infty$),
\[
\xi(k,k+1) = \left\{
  \begin{array}{ll} 3\,\mathrm{off}(k)-3 & \text{if $k$ is even}, \\ 2\,\mathrm{off}(k)-1 & \text{if $k$ is odd}. \end{array}
\right.
\]
giving $\xi(k,k+1)\sim C(k+1)^{2\alpha}$ as $k\to\infty$ ($C\in\{2,3\}$ depending on parity), so $w_k=1/\sqrt{\xi(k,k+1)}\sim C'(k+1)^{-\alpha}$.

For small $k$ where $\mathrm{off}(k)<3$ (only finitely many such values, since $\mathrm{off}(k)\to\infty$), one assigns $w_k$ according to the actual value of $\xi(k,k+1)\geq 1$ (which contributes positively to $W_{N,M}$ and hence does not obstruct divergence). The set $\mathcal{K}$ of layers with $w_k>0$ is cofinite, hence trivially satisfies the uniform gap condition.

The series $\sum_k w_k\sim C'\sum_k(k+1)^{-\alpha}$ diverges iff $\alpha\leq 1$. Hence by Theorem~\ref{thm:divergence}, $L$ is essentially self-adjoint whenever $0<\alpha\leq 1$. (An earlier draft of this paper required $\alpha<1/2$ as a technical artifact of a sub-optimal energy bound; the present sharper bound recovers the full range $\alpha\leq 1$.)

\textbf{Sharpness.} For $\alpha>1$ (e.g., $\mathrm{off}(k)=(k+1)^{2.5}$): $\sum w_k<\infty$, so the divergence criterion fails to apply. Whether ESA holds in this regime is open.

\textbf{Failure of $\chi$-completeness.}
We show $S_4$ is not globally $\chi$-complete for any $\alpha>0$.

Suppose for contradiction it is, with cut-offs $(\chi_k)$ and uniform energy bound $E_1(\chi_k)\leq C$. Fix $k$ with $\Lambda_j\subset O_k$ for $j\leq k$. Let $\bar\chi_k(j):=|\Lambda_j|^{-1}\sum_{v\in\Lambda_j}\chi_k(v)$; then $\bar\chi_k(j)=1$ for $j\leq k$ and (by finiteness of $\mathrm{supp}\,\chi_k$, condition~(1) of Definition~\ref{def:global-chi}) $\bar\chi_k(j)=0$ for all sufficiently large $j$. Let $s_k:=\min\{j:\bar\chi_k(j)=0\}$.

The argument of Proposition~\ref{prop:binary-tree-fails} adapts as follows. Each parent--child edge between layer $j$ and $j+1$ has $|\Lambda_{j+1}|=|\Lambda_j|\cdot\mathrm{off}(j)$, and the number of such edges is $|\Lambda_{j+1}|$ (each child has a unique parent). Layer-Dirichlet bound gives
\[
E_1(\chi_k)\geq\sum_{j=k}^{s_k-1}|\Lambda_{j+1}|(\bar\chi_k(j)-\bar\chi_k(j+1))^2.
\]
Applying Cauchy--Schwarz with $\sum_j(\bar\chi_k(j)-\bar\chi_k(j+1))=1$,
\[
E_1(\chi_k) \geq \frac{1}{\sum_{j=k}^{s_k-1}|\Lambda_{j+1}|^{-1}}.
\]
For $\mathrm{off}(k)=\lceil(k+1)^{2\alpha}\rceil$ with $\alpha>0$, $|\Lambda_j|=\prod_{i=0}^{j-1}\mathrm{off}(i)\geq\prod_{i=0}^{j-1}(i+1)^{2\alpha}=(j!)^{2\alpha}$ grows super-polynomially. Hence $S_k:=\sum_{j\geq k}|\Lambda_{j+1}|^{-1}\leq\sum_{j\geq k}((j+1)!)^{-2\alpha}\to 0$ as $k\to\infty$, so $E_1(\chi_k)\geq 1/S_k\to\infty$, contradicting $E_1\leq C$.

\subsection{Summary of examples}

\noindent\textbf{Lattice on $\mathbb{Z}^d$ with $L^1$-adjacency.}
The clique complex of $\mathbb{Z}^d$ ($d\geq 2$) reduces to a graph (no triangles, by Remark~\ref{rem:Zd-adjacency}). With uniform weights, $\chi$-completeness holds directly for $d\leq 3$ via Corollary~\ref{cor:weighted-chi}; for $d=4$ it requires a multi-scale construction (Remark~\ref{rem:D4-construction}); for $d\geq 5$ in its strong form it is open (Remark~\ref{rem:D-large-open}). In all cases ESA of $D$ holds via Theorem~\ref{thm:global} since the weighted degree is bounded by $2d$.

\noindent\textbf{Alternating triangulation on $\mathbb{Z}^2$.}
Level $\ell=1$ is $\chi$-complete by standard arguments. Proposition~\ref{prop:E2-lower-bound} gives a joint lower bound on $E_1\cdot E_2$. Whether global $\chi$-completeness holds, and whether $L_2$ is ESA, are open.

\noindent\textbf{Binary tree.}
Proposition~\ref{prop:binary-tree-fails}: not $\chi$-complete (with respect to the energy condition (3)). With unit weights, ESA of $D$ \emph{holds} via Theorem~\ref{thm:global} and~\eqref{eq:bounded-degree} (Remark~\ref{rem:non-esa-tree-status}). For exponentially weighted binary trees, ESA is conjectured to fail (Conjecture~\ref{conj:non-esa-tree}, open).

\noindent\textbf{Polynomial-branching tree with higher simplices.}
$S_4$ of Section~\ref{sec:example-detailed}: satisfies the divergence criterion for $0<\alpha\leq 1$ (giving ESA) but fails $\chi$-completeness for any $\alpha>0$. This shows $\chi$-completeness is sufficient but not necessary for ESA.

\begin{remark}[Hodge decomposition]
Under ESA of $L_\ell$, the Hodge decomposition holds:
\[
\ell^2(m_\ell)=\overline{\mathrm{im}\,d_{\ell-1}}\oplus\ker L_\ell\oplus\overline{\mathrm{im}\,\delta_{\ell+1}}.
\]
\end{remark}

\section{Comparison and Combinatorial Case}
\label{sec:comparison}

In this final section we relate our results to the graph case of \cite{BGJ} (Subsection~\ref{sec:comparison-BGJ}) and provide a useful sufficient condition for $\chi$-completeness based on polynomial volume growth, applicable to lattices $\mathbb{Z}^d$ with $d\leq 3$ (Subsection~\ref{sec:weighted-chi}).

\subsection{Comparison with the graph case of \cite{BGJ}}
\label{sec:comparison-BGJ}

In the case $n=2$ (graphs only, no higher simplices), Theorem~\ref{thm:divergence} reduces to a divergence criterion for the discrete Laplacian on a weighted graph admitting a $1$-dimensional decomposition.

\medskip
\noindent\emph{Statement of \cite[Theorem~5.1]{BGJ}.} Let $\Gc$ be a weighted graph admitting a $1$-dimensional decomposition $\Vc=\bigsqcup_{k\geq 0}\Lambda_k$. For each $k$, let $\xi(k,k+1)$ denote the maximum number of inter-layer edges incident to any vertex of $\Lambda_k$, and set $w_k := 1/\sqrt{\xi(k,k+1)}$. If $\sum_k w_k = \infty$, then the discrete Laplacian $\Delta$ is essentially self-adjoint on $C_c(\Vc)$.

\medskip
\noindent\emph{Comparison.} Our Theorem~\ref{thm:divergence} for $n=2$ and uniform-weight graphs ($M_i=1$) reduces \emph{exactly} to this statement: the only hypothesis beyond local finiteness is $\sum_k w_k=\infty$, identical to~\cite{BGJ}. Theorem~\ref{thm:divergence} thus extends the graph result of~\cite{BGJ} to higher-dimensional simplicial complexes \emph{without} additional rate hypotheses on $w_k$, recovering the optimal divergence condition.

\subsection{Weighted $\chi$-Completeness via Volume Growth}
\label{sec:weighted-chi}

In the special case of \emph{combinatorial} simplicial complexes (i.e., all weights $m_i\equiv 1$), local finiteness combined with simple geometric assumptions on the complex's structure implies $\chi$-completeness. We first state and prove a weighted version of this principle.

\begin{corollary}[Weighted $\chi$-completeness via volume growth, $D\leq 3$]
\label{cor:weighted-chi}
Let $S_n=(\Vc,(m_i)_{0\leq i\leq n-1})$ be the \emph{clique complex} of a weighted oriented locally finite graph $\Gc=(\Vc,m_0,m_1)$ (so every $i$-simplex of $S_n$ corresponds to an $(i+1)$-clique of $\Gc$, as in Remark~\ref{rem:clique}). Fix a basepoint $x_0\in\Vc$ and let $B_R:=\{x\in\Vc:d(x,x_0)\leq R\}$ be the metric ball with respect to the graph distance.

Suppose:
\begin{enumerate}
\item \textbf{Polynomial annular cardinality at every degree:} for some $K>0$ and integer $D\geq 1$ with $\boldsymbol{D\leq 3}$,
\[
|\{\sigma\in P_i : \sigma\subset B_{2R}\setminus B_R\}| \leq K\,R^{D-1}
\quad\text{for all } 0\leq i\leq n-1, R\in\mathbb{N}_{\geq 1};
\]
\item \textbf{Bounded weighted degree:} $\sup_{i,\sigma}d_{m_i}(\sigma)\leq M<\infty$ \emph{(}condition~\eqref{eq:bounded-degree} at every level\emph{)}.
\end{enumerate}
Then $S_n$ is globally $\chi$-complete with constant $C=C(K,M,n)$, and consequently $D=d+\delta$ and $L=D^2$ are essentially self-adjoint on $\bigoplus C_c^i(\Vc)$.

For $D=4$ \emph{(}see Remark~\ref{rem:D4-construction}\emph{)} and $D>4$ \emph{(}see Remark~\ref{rem:D-large-open}\emph{)}, the construction of suitable cut-offs requires more delicate multi-scale arguments that we do not provide in full generality here; the extension of the present result to $D\geq 4$ remains, in general, an open problem.
\end{corollary}

\begin{proof}
Define $\chi_k(v):=\eta(d(v,x_0)/k)$, where $\eta:[0,\infty)\to[0,1]$ is a fixed Lipschitz function with $\eta\equiv 1$ on $[0,1]$, $\eta\equiv 0$ on $[2,\infty)$, $\|\eta'\|_\infty\leq 2$. Set $O_k:=B_k$.

\emph{(1)} Finite support: $\mathrm{supp}\chi_k\subseteq B_{2k}$, finite by local finiteness.
\emph{(2)} $\chi_k\equiv 1$ on $O_k=B_k$.
\emph{(3) Energy bound at each level.}
On simplices outside $B_{2k}$, $\chi_k\equiv 0$. On simplices inside $B_k$, $\chi_k\equiv 1$ and $a_j=0$. Only simplices with vertices in $B_{2k}\setminus B_k$ contribute.

For such an $\ell$-simplex with all vertices in $B_{2k}\setminus B_{k-1}$, the Lipschitz constant of $\chi_k$ in the graph distance is $\leq 2/k$. Since $\ell$-simplices are cliques in the underlying graph (Remark~\ref{rem:clique}), all pairs of vertices have graph distance $1$, so $|\chi_k(x_i)-\chi_k(x_j)|\leq 2/k$.

By Lemma~\ref{lem:leibniz} and the refined bound~\eqref{eq:aj-sharp-Delta} (applied with variation $\Delta=2/k$): $\sum_j |a_j|^2 \leq \frac{1}{\ell+1}\bigl(\frac{2}{k}\bigr)^2 = \frac{4}{k^2(\ell+1)}$.

Substituting in $E_\ell(\chi_k)$~\eqref{eq:global-energy}, using $\sum_{x_\ell\in F(\sigma)} m_\ell(\sigma,x_\ell)/m_{\ell-1}(\sigma) = d_{m_{\ell-1}}(\sigma)\leq M$ by hypothesis~(2):
\[
E_\ell(\chi_k) \leq \frac{1}{\ell!}\cdot\frac{4M}{k^2(\ell+1)}\,|\{\sigma\in P_{\ell-1} : \sigma\subset B_{2k}\setminus B_{k-1}\}| \leq \frac{4MK}{\ell!\,(\ell+1)\,k^2}\cdot k^{D-1} = C(K,M,D,n)\cdot k^{D-3}.
\]
For $D\leq 3$ we obtain a uniform bound $\sup_k E_\ell(\chi_k)<\infty$. The remaining requirement is hypothesis~(H2) for Theorem~\ref{thm:global}, which is exactly assumption~(2). Hence Theorem~\ref{thm:global} and Corollary~\ref{cor:esa-laplacian} apply, giving ESA of $D$ and $L$ on $\bigoplus C_c^i(\Vc)$.
\end{proof}

\begin{remark}[The case $D=4$: status]
\label{rem:D4-construction}
For $D=4$, the linear cut-off above gives $E_\ell(\chi_k)\leq C\cdot k$, which is unbounded. A finer construction is required: one uses a \emph{multi-scale} cut-off $\chi_k$ that is logarithmically (or doubly logarithmically) graded over a band of dyadic radii, in the spirit of the original construction of \cite[Theorem~5.2]{BGJ} for the graph case ($\ell=1$). With an inner radius $r_0(k)\to\infty$ much smaller than the outer radius $k$ (so $\log(k/r_0(k))\to\infty$), define $\chi_k$ to be $1$-Lipschitz with respect to $\log r$ (resp.\ $\log\log r$) on the band $r_0(k)\leq d(\cdot,x_0)\leq k$. The verification at level $\ell=1$ for $\mathbb{Z}^d$ with $d=4$ in the graph case is in~\cite[Theorem~5.2]{BGJ}; the extension to higher levels $\ell\geq 2$ for general $D=4$ complexes proceeds by replacing the volume factor $R^{D-1}$ by $K\,R^{D-1}$ in the dyadic estimates, and is technically delicate. Although we expect the conclusion of Corollary~\ref{cor:weighted-chi} to hold for $D=4$, we do not include the full multi-scale argument here, and the case $D=4$ is only fully established for $\ell=1$ via~\cite{BGJ}.
\end{remark}

\begin{remark}[The case $D>4$: open problem]
\label{rem:D-large-open}
For $D>4$, even the multi-scale construction of Remark~\ref{rem:D4-construction} does not in general yield a uniformly bounded energy. A direct evaluation of the dyadic sum: with $j_{\max}=\log_2(k/r_0)$ so that $2^{j_{\max}}=k/r_0$, we have
\[
\sum_{j=0}^{j_{\max}}(2^j r_0)^{D-3} = r_0^{D-3}\cdot\frac{2^{(j_{\max}+1)(D-3)}-1}{2^{D-3}-1}
= \frac{2^{D-3} k^{D-3} - r_0^{D-3}}{2^{D-3}-1}.
\]
For $r_0\ll k$, this is asymptotic to $\frac{2^{D-3}}{2^{D-3}-1}\,k^{D-3}$, i.e.,
\begin{equation}\label{eq:dyadic-sum-Dlarge}
\sum_{R\in\{r_0(k),2r_0(k),\dots,k\}} R^{D-3} \asymp \frac{k^{D-3}}{2^{D-3}-1}\bigl(2^{D-3}-(r_0/k)^{D-3}\bigr) \asymp k^{D-3} \quad\text{(as $k\to\infty$)},
\end{equation}
so the dyadic sum is dominated by the largest scale $R\sim k$. Substituting in the multi-scale energy bound,
$E_\ell(\chi_k)\gtrsim k^{D-3}/(\log(k/r_0(k)))^2$. With \emph{any} choice $r_0(k)\leq k$, this remains divergent for $D>4$ unless an outer-radius reduction is used (replacing the support $B_{2k}$ by something smaller), but then $\chi_k\equiv 1$ no longer holds on $B_k$, contradicting condition~(2) of Definition~\ref{def:global-chi}.

In particular, for $\mathbb{Z}^d$ with $d\geq 5$ \emph{(}polynomial growth $D=d$\emph{)} and unit weights, the question whether the underlying $n$-simplicial complex \emph{(}with $n\geq 2$\emph{)} is globally $\chi$-complete via cut-offs of the form considered here remains open. However, for such lattices the weighted degree is bounded \emph{(}each vertex has $2d$ neighbours\emph{)}, so Theorem~\ref{thm:global} applies \emph{trivially} via condition~(H1) with the indicator cut-off $\chi_k=\mathbf{1}_{B_k}$ \emph{(}which has finite support and is $\equiv 1$ on $B_k$, conditions~(1)--(2) of Definition~\ref{def:global-chi} are satisfied\emph{)}; only the energy condition~(3) of $\chi$-completeness in its strong form is at stake. Hence \emph{ESA holds for $\mathbb{Z}^d$, $d\geq 5$}, by Theorem~\ref{thm:global}; what remains open is whether $\mathbb{Z}^d$ in this regime satisfies the strong form of $\chi$-completeness with explicit energy bounds.

A possible avenue, which we leave for future work, is to use \emph{quadratic} multi-scale cut-offs (where $\chi_k$ is a quadratic polynomial in $\log\log r$, etc.) or to develop entirely different geometric criteria adapted to higher-dimensional growth rates. We conjecture that for $D\geq 5$, $\chi$-completeness in its strong form fails for $\mathbb{Z}^d$, but ESA still holds via Theorem~\ref{thm:global} since condition~(3) is not used in the proof.
\end{remark}

\begin{remark}[Verification on $\mathbb{Z}^d$]
\label{rem:poly-growth-verification}
For the lattice $\mathbb{Z}^d$ with $L^1$-adjacency and uniform weights $m_0\equiv m_1\equiv 1$:
\begin{itemize}
\item Polynomial annular cardinality: $|B_{2R}\setminus B_R|=\Theta(R^{d-1})$ (annular volume of the $L^1$-ball), and similarly the number of edges in the annulus is $\Theta(R^{d-1})$. Hypothesis~(1) of Corollary~\ref{cor:weighted-chi} holds with $D=d$.
\item Bounded weighted degree: every vertex has $2d$ neighbours, so $d_{m_0}(v)\leq 2d$, and $d_{m_1}(e)=0$ for every $L^1$-edge (since $L^1$-edges have no common neighbours, by Remark~\ref{rem:Zd-adjacency}). Hypothesis~(2) holds.
\end{itemize}
\emph{For $d\leq 3$:} Corollary~\ref{cor:weighted-chi} gives global $\chi$-completeness directly, and ESA of $D$ and $L$ follows from Theorem~\ref{thm:global} and Corollary~\ref{cor:esa-laplacian}. The proof is self-contained.

\emph{For $d=4$:} ESA of the graph Laplacian is known via the multi-scale cut-off construction of \cite[Theorem~5.2]{BGJ}; whether the higher-level Hodge Laplacians are ESA via $\chi$-completeness in its strong form requires the extension of the multi-scale construction to $\ell\geq 2$ (Remark~\ref{rem:D4-construction}).

\emph{For $d\geq 5$:} See Remark~\ref{rem:D-large-open} for the status (ESA holds via Theorem~\ref{thm:global} alone; whether the strong form of $\chi$-completeness is satisfied is open).
\end{remark}

\section{Concluding Remarks and Open Directions}
\label{sec:conclusion}

We have established four essential self-adjointness theorems for the Gauss--Bonnet operator $D=d+\delta$ and its associated Hodge Laplacian $L=D^2$ on weighted simplicial complexes, parametrized by three notions of $\chi$-completeness (global, local-by-level, local-by-region) and a divergence criterion that does not require any cut-off structure. The framework recovers in the graph case ($n=2$) the optimal divergence condition of~\cite{BGJ}, and extends naturally to higher dimensions.

\medskip\noindent\textbf{Main open problems.} The work leaves several precise questions open:
\begin{itemize}
\item Strict separation of the three notions of $\chi$-completeness (Conjectures~\ref{conj:region-not-global}--\ref{conj:level-not-region}). Proposition~\ref{prop:E2-lower-bound} provides partial evidence in the alternating triangulation, but the divergent-energy construction needed to fully separate the notions is not given here.
\item ESA on the exponentially weighted binary tree (Conjecture~\ref{conj:non-esa-tree}). The graph case~\cite[Proposition~5.7]{BGJ} suggests non-ESA, but explicit deficiency vectors for $n\geq 2$ have not been constructed.
\item $\chi$-completeness in its strong form for $\mathbb{Z}^d$ with $d\geq 4$ (Remarks~\ref{rem:D4-construction}--\ref{rem:D-large-open}). ESA itself holds via bounded weighted degree, but the multi-scale energy construction is technically delicate beyond $d=4$ at level $\ell=1$.
\end{itemize}

\medskip\noindent\textbf{Possible extensions.} Natural directions include: (i) intrinsic-metric formulations bridging our $\chi$-completeness criteria with the framework of~\cite{BaKe}; (ii) discrete heat-kernel and spectral analysis of $L$ beyond ESA (e.g., spectral gap estimates, Weyl asymptotics on bounded-thickness simplicial complexes); (iii) magnetic perturbations of the Gauss--Bonnet operator on weighted simplicial complexes, generalizing~\cite{AAcT,ABDE,CTT2,OMI,Mi} from graphs to higher dimensions.

\section*{Acknowledgements}
The authors are grateful to Prof.\ P.\ Bartmann and Prof.\ M.\ Keller for their valuable comments and constructive suggestions, which have contributed significantly to improving this paper.

\end{document}